\documentclass[11pt]{article}
\usepackage{amsmath,amssymb,amsfonts,times}
\usepackage{bm}
\usepackage{graphicx,epsfig}
\usepackage{caption,subfigure}
\usepackage{multirow,booktabs}
\usepackage{enumerate}
\usepackage{color}

\newtheorem{theorem}{Theorem}

\newtheorem{remark}{Remark}

\newcommand\red{\textcolor[rgb]{0.00,0.00,0.00}}

\allowdisplaybreaks \numberwithin{equation}{section}

\title{An extended sampling-ensemble Kalman filter approach for partial data inverse elastic problems\thanks{The work of the first and third authors was partially supported by the NNSF of China (National Natural Science Foundation of China)[grant number 11771068].}
}

\author{
\red{Zhaoxing Li
}\thanks{
School of Mathematical Sciences, University of Electronic Science and Technology of China, Chengdu, 611731, China (lzx130682@163.com).
}
\quad\red{
Jiguang Sun
}\thanks{
Department of Mathematical Sciences, Michigan Technological University, Houghton, MI 49931 (jiguangs@mtu.edu).
}
\quad\red{
Liwei Xu
}\thanks{
School of Mathematical Sciences, University of Electronic Science and Technology of China, Chengdu, 611731, China (xul@uestc.edu.cn).
}}
\date{}

\begin{document}

\maketitle

\begin{quote}
\small
{\bf Abstract:}\;
Inverse problems are more challenging when only partial data are available in general. In this paper, we propose a two-step approach combining the extended sampling method and the ensemble Kalman filter to reconstruct an elastic rigid obstacle using partial data. In the first step, the approximate location of the unknown obstacle is obtained by the extended sampling method. In the second step, the ensemble Kalman filter is employed to reconstruct the shape. The location obtained in the first step guides the construction of the initial particles of the ensemble Kalman filter, which is critical to the performance of the second step. Both steps are based on the same physical model and use the same scattering data. Numerical examples are shown to illustrate the effectiveness of the proposed method.


\end{quote}

\section{Introduction }

Inverse scattering theory is an active research area in mathematics and engineering.
It has many important applications such as non-destructive testing, seismology, and geological exploration.
In this paper, we consider the inverse elastic scattering problem to
determine the location and shape of the obstacle from the measured displacement field in the frequency domain.
Due to how much data is available, such problems are usually divided into full aperture and limited aperture problems.
The full aperture problems have been extensively studied, and many existing methods can achieve satisfactory reconstructions
(see, e.g., \cite {Bao2018,Charalambopoulos2006,Hu2012,Kress1996}).
However, for applications such as underground explorations, full aperture data are not available.
It is desirable to develop effective numerical methods for limited aperture data.

Compared to the full aperture case, limited aperture problems are more challenging in general.
One approach is to first recover the full aperture data and then apply the existing methods.
However, it is a severely ill-posed problem to recover the full aperture data using analytic continuation or optimization.
We refer the readers to \cite{Ahn2014,Bao2003,Ikehata2012,Zinn1989,Ochs1987,LiuSun2019} for more discussions.

Another approach for the limited aperture problems is to take the advantages of different inversion methods by combining them in a suitable way \cite{Li2020, LiEtal2020}.
In this paper, we continue the investigation along this direction and propose a two-step method combining the extended sampling method (ESM)
and the ensemble Kalman filter (EnKF).
The ESM is a qualitative method which was originally proposed in \cite{Liu2018} for the acoustic inverse scattering problem to reconstruct the approximate location and size of the scatterer
using the scattering data due to one incident wave. It was extended to the inverse elastic problem in \cite{Liu2019}.
The key ingredient of the ESM is a new far-field equation, whose regularized solution is used to define an indicator for the unknown scatterer.
For classical sampling methods such as the linear sampling method, the kernel of the far field equation is the measured full aperture scattering data
of all incident and observation directions \cite{Colton2013}. In contrast, the integral kernel in the ESM is the full aperture far-field data of a known rigid disc.
The measured data is moved to the right hand side of the integral equations. This arrangement enables the ESM to treat limited aperture data flexibly.
As the first step of the proposed method, we modify the ESM to reconstruct the approximate location of the obstacle using limited aperture data.

In the second step, the EnKF is employed to refine the location and construct the shape of the unknown obstacle using the same data.
The EnKF can be regarded as a Monte Carlo variation of the standard Kalman filter (KF).
The mean and covariance are approximated by an ensemble of particles, and the propagation of these particles are encoded in an iteration process
through the standard Kalman update formula. Due to its robustness, ease of implementation, and accuracy for state estimation of partially observed dynamical systems,
the EnKF has applications in many areas such as oceanography, meteorology, and geophysics. For inverse problems, the EnKF only uses the forward operator and
the Fr\'{e}chet derivative is not needed. For more details about the EnKF and its applications to inverse problems,
one can refer to \cite{Chada2018,Iglesias2016,Kaipio2005,Schillings2017} and references therein.
The initial ensemble of the particles in the EnKF is critical to its performance.
In \cite{Iglesias2013}, it is proved that the inversion solution generated by the EnKF lies in the linear span of the initial ensemble.
The approximate location obtained by the ESM in the first step is used to construct the initial ensemble of particles.
Both steps use the same physical model and measured data. Numerical experiments show that this approach
inherits the merits of the two methods and can effectively recover the obstacle using limited aperture data.
We refer the readers to \cite{Li2020} and \cite{LiEtal2020} for the applications of the combined approach to an inverse scattering problem and an inverse acoustic source problems, respectively.

The rest of the paper is arranged as follows. In Section 2, we introduce the limited aperture inverse obstacle scattering problem
and an equivalent form of the elastic equation based on the Helmholtz decomposition. In Section 3, we develop a modified ESM to find the approximate location of the unknown obstacle.
In Section 4, the EnKF is employed to recover the shape of the obstacle. In Section 5, numerical experiments are presented to demonstrate the effectiveness of the proposed method.
Finally, we draw some conclusions in Section 6.

\section{Direct and Inverse Elastic Scattering Problems}

For $\bm{x}=(x_1,x_2)^\top\in \mathbb{R}^2$, let $\hat{\bm{x}}:=\bm{x}/|\bm{x}|\in \mathbb{S}$, where $\mathbb{S}=\{\bm{\hat{x}}\in \mathbb{R}^2,|\bm{\hat{x}}|=1\}$ is the unit circle.
Let $\hat{\bm{x}}^\bot\in \mathbb{S}$ be the vector obtained by rotating $\hat{\bm{x}}$ counterclockwise $\pi/2$ .
For a scaler function $v$ and a vector function $\bm{v}=(v_1,v_2)^\top$, define the curl operators
$\textbf{curl}\,v=\left(\partial_{x_2}v,-\partial_{x_1}v\right)^\top$ and $\text{curl}\,\bm{v}=\partial_{x_1} v_2-\partial_{x_2} v_1$, respectively.
Denote by $\Omega\subset \mathbb{R}^2$ a rigid obstacle with $C^2$-boundary and assume that $\mathbb{R}^2\backslash \overline{\Omega}$ is occupied by isotropic homogeneous elastic solid.
Denote by $\bm{\tau}=(\tau_1,\tau_2)^\top$ the unit tangential and by $\bm{\nu}=(\nu_1,\nu_2)^\top$ the unit outward normal vector on $\partial \Omega$, respectively,
where $\tau_1=-\nu_2$ and $\tau_2=\nu_1$.

The time-harmonic elastic scattering problem is to find  $\bm{v}$ satisfying the Navier equation
\begin{equation}\label{Navier}
\mu\Delta \bm{v}+(\lambda+\mu)\nabla\nabla\cdot \bm{v}+\omega^2\bm{v}=0\quad \text{in}\; \mathbb{R}^2\backslash \overline{\Omega},
\end{equation}
where $\mu$ and $\lambda$ are the L\'{a}me constants such that $\mu>0,\;\lambda+\mu>0$, $\omega>0$ is the angular frequency.
In \eqref{Navier}, $\bm{v}=\bm{u}^{inc}+\bm{u}$ is the total displacement field, $\bm{u}^{inc}$ is the incident field, and $\bm{u}$ is the scattered field.

The incident field is the plane wave given by
\begin{equation}\label{uinc}
\bm{u}^{inc}(\bm{x})=\bm{u}^{inc}(\bm{x};\bm{d},\alpha_p,\alpha_s)=\alpha_p \bm{d}e^{i k_p \bm{x}\cdot \bm{d}}+\alpha_s \bm{d}^{\bot}e^{i k_s \bm{x}\cdot \bm{d}},\; \alpha_p,\alpha_s \in  \mathbb{C},
\end{equation}
where $\bm{d}\in \mathbb{S}$ is the incident direction, $k_p=\omega/\sqrt{\lambda+2\mu}$ and $k_s=\omega/\sqrt{\mu}$ are the compressional and shear wave numbers, respectively.
For a rigid obstacle $\Omega$, the total field $\bm{v}$ satisfies the boundary condition
\begin{equation}\label{bmvBC}
\bm{v}={\boldsymbol 0}\;\; \text{on}\; \partial \Omega.
\end{equation}
Hence $\bm{u}=\bm{v}-\bm{u}^{inc}$ satisfies the following boundary value problem
\begin{equation}\label{usc}
\left\{\begin{array}{ll}
\mu\Delta \bm{u}+(\lambda+\mu)\nabla\nabla\cdot \bm{u}+\omega^2\bm{u}=0\; &\text{in}\; \mathbb{R}^2\backslash \overline{\Omega},\\
\bm{u}=-\bm{u}^{inc}\; &\text{on}\; \partial \Omega.
\end{array}\right.
\end{equation}
The solution $\bm{u}$ can be decomposed as $\bm{u}=\bm{u}_p+\bm{u}_s$, where the compressional wave $\bm{u}_p$ and the shear wave $\bm{u}_s$ are given by
\begin{equation*}
\bm{u}_p=-\frac{1}{k_p^2}\nabla\nabla\cdot \bm{u},\quad \bm{u}_s=\frac{1}{k_s^2}\textbf{curl}\,
\text{curl}\,\bm{u}.
\end{equation*}
In addition, $\bm{u}$ satisfies the Kupradze radiation condition
\begin{equation}\label{Kupradze}
\lim\limits_{r\to\infty}\sqrt{r}\left(\partial_r \bm{u}_t-i k_t \bm{u}_t\right)=0,\quad t=p \text{ or }s,\;r=|\bm{x}|.
\end{equation}

The solution $\bm{u}$ to \eqref{usc}-\eqref{Kupradze} has the following asymptotic expansion \cite{Arens2001}
\begin{equation}\label{asymptotic}
\bm{u}(\bm{x})=\frac{e^{ik_p|\bm{x}|}}{\sqrt{|\bm{x}|}}
{u}^{\infty}_p(\bm{\hat{x}})\bm{\hat{x}}+
\frac{e^{ik_s|\bm{x}|}}{\sqrt{|\bm{x}|}}
{u}^{\infty}_s(\bm{\hat{x}})\bm{\hat{x}}^{\bot}+
\mathcal{O}\left({|\bm{x}|^{-3/2}}\right),
\quad |\bm{x}|\to\infty,
\end{equation}
uniformly in all direction $\bm{\hat{x}}=\bm{x}/|\bm{x}|\in \mathbb{S}$, where $u^{\infty}_p$ and $u^{\infty}_s$ defined on $\mathbb{S}$
are the compressional and shear far-field pattern of $\bm{u}$, respectively.

Let $\gamma^o \subset \mathbb{S}$ and $\gamma^i \subset \mathbb{S}$ be the observation aperture and the incident aperture, respectively.
The inverse obstacle scattering problems (IOSP) considered in this paper are as follows:

\begin{itemize}
\item \textbf{IOSP-P:} Determine $\partial \Omega$ from ${u}_p^{\infty}(\bm{\hat{x},d})$, $(\bm{\hat{x}},\bm{d})\in \gamma^o\times \gamma^i$, $\gamma^o \times \gamma^i \subsetneq \mathbb{S} \times \mathbb{S}$.

\item \textbf{IOSP-S:} Determine $\partial \Omega$ from ${u}_s^{\infty}(\bm{\hat{x},d})$, $(\bm{\hat{x}},\bm{d})\in \gamma^o\times \gamma^i$, $\gamma^o \times \gamma^i \subsetneq \mathbb{S} \times \mathbb{S}$.

\item \textbf{IOSP-F:} Determine $\partial \Omega$ from $\bm{u}_{\infty}(\bm{\hat{x},d})=({u}_p^{\infty};{u}_s^{\infty})$, $(\bm{\hat{x}},\bm{d})\in \gamma^o\times \gamma^i$, $\gamma^o,\gamma^i \subsetneq \mathbb{S}$.
\end{itemize}

We end this section by introducing an equivalent form of the Navier equation \eqref{usc}-\eqref{Kupradze} (see, e.g., \cite{Li2015,Dong2019}). For a solution $\bm{u}$ of \eqref{usc}-\eqref{Kupradze}, the Helmholtz decomposition holds
\begin{equation}\label{HelmDecomp}
\bm{u}=\nabla \phi+\textbf{curl}\,\psi,
\end{equation}
where $\phi$ and $\psi$ are two scalar functions. Using \eqref{HelmDecomp} and \eqref{usc}-\eqref{Kupradze}, $\phi$ and $\psi$ satisfies
\begin{equation}\label{phipsi}
\left\{\begin{array}{ll}
\Delta \phi+k_p^2 \phi=0,\quad &\text{in}\; \mathbb{R}^2\backslash \overline{\Omega},\\
\Delta \psi+k_s^2 \psi=0,\quad &\text{in}\; \mathbb{R}^2\backslash \overline{\Omega},\\
\frac{\partial \phi}{\partial \bm{\nu}}+\frac{\partial \psi}{\partial \bm{\tau}} =g_1,\quad &\text{on}\; \partial \Omega,\\
\frac{\partial \phi}{\partial \bm{\tau}}-\frac{\partial \psi}{\partial \bm{\nu}} =g_2,\quad &\text{on}\; \partial \Omega,\\
\lim\limits_{r\to\infty}\sqrt{r}\left(\frac{\partial \phi}{\partial r}-i k_p \phi\right)=0,\quad &r=|\bm{x}|,\\
\lim\limits_{r\to\infty}\sqrt{r}\left(\frac{\partial \psi}{\partial r}-i k_s \psi\right)=0,\quad &r=|\bm{x}|,
\end{array}\right.
\end{equation}
where $g_1=-\bm{\nu}\cdot \bm{u}^{inc}$ and $g_2=-\bm{\tau}\cdot \bm{u}^{inc}$.
The relation between the solutions of \eqref{usc}-\eqref{Kupradze} and \eqref{phipsi} is stated in the following theorem.

\begin{theorem}\label{thm1}\cite{Liu2019,Dong2019}
Let $\bm{u}$ be the solution of \eqref{usc}-\eqref{Kupradze}. Then
\begin{equation}\label{relation_uphi}
\phi=-\frac{1}{k_p^2}\nabla\cdot \bm{u} \;\;\; \text{and} \;\;\; \psi=\frac{1}{k_s^2}\text{curl}\,\bm{u},
\end{equation}
are the solution of the Helmholtz equation \eqref{phipsi}. Moreover, the far field patterns $\phi_\infty$ and $\psi_\infty$ of $\phi$ and $\psi$
satisfy
\begin{equation}\label{relation_uinfphiinf}
\phi_\infty(\bm{\hat{x}})=\frac{1}{i k_p}u_p^\infty(\bm{\hat{x}}) \;\;\; \text{and} \;\;\; \psi_\infty(\bm{\hat{x}})=\frac{1}{i k_s}u_s^\infty(\bm{\hat{x}}).
\end{equation}
\end{theorem}

\section{Extended Sampling Method}

As the first step of the combined approach, we consider the problem of finding the approximate location of the obstacle $\Omega$.
In this section, we modify the extended sampling method (ESM) in \cite{Liu2019} for the limited aperture inverse elastic obstacle problem.

\subsection{ESM for IOSP-P and IOSP-S}

We first consider the case when the measured data is the far-field pattern $u^\infty_t(\hat{\bm{x}},\bm{d}_0)$, $t=p,s$ of $\Omega$ due to one incident direction $\bm{d}_0$
and of all observation directions $\hat{\bm{x}}\in \mathbb{S}$. Denote by $B_{\bm{z}}\subset \mathbb{R}^2$ a rigid disc centered at $\bm{z}$ with radius large enough.
Let $u^{B_{\bm{z}}}(\hat{\bm{x}};k_t,\bm{d})$, $t=p,s$, be the solution of
\begin{equation}\label{ubz}
\left\{\begin{array}{ll}
\Delta u+k_t^2 u=0,\quad &\text{in}\; \mathbb{R}^2\backslash \overline{B}_{\bm{z}},\\
u=-e^{i k_t \bm{x}\cdot \bm{d}}, \quad &\text{on}\; \partial B_{\bm{z}},\\
\lim\limits_{r\to\infty}\sqrt{r}\left(\frac{\partial u}{\partial r}-i k_t u\right)=0,\quad &r=|\bm{x}|.
\end{array}\right.
\end{equation}
Let $u^{B_{\bm{z}}}_\infty(\bm{\hat{x}};k_t,\bm{d})$ be the far field pattern of $u^{B_{\bm{z}}}(\bm{\hat{x}};k_t,\bm{d})$.
Define the far-field operator $\mathcal{F}_{\bm{z}}: L^2(\mathbb{S})\rightarrow L^2(\mathbb{S})$ such that
\begin{equation}\label{Fzg}
\mathcal{F}_{\bm{z}} g(\bm{\hat{x}})=\int_{\mathbb{S}} u^{B_{\bm{z}}}_\infty(\bm{\hat{x}};k_t,\bm{d}) g(\bm{d}) d s(\bm{d}),\;\; \bm{\hat{x}}\in \mathbb{S}.
\end{equation}
Using $\mathcal{F}_{\bm{z}}$, for the far-field data $u_t^\infty(\hat{\bm{x}},\bm{d}_0)$, $ \hat{\bm{x}}\in \mathbb{S}$, we set up a far-field equation
\begin{equation}\label{farfieldeqn1}
  (\mathcal{F}_{\bm{z}} g)(\bm{\hat{x}})=\frac{1}{i k_t} u^\infty_t(\hat{\bm{x}},\bm{d}_0),\;\; \hat{\bm{x}}\in \mathbb{S}.
\end{equation}
The approximate location of $\Omega$ can be obtained using the solutions of \eqref{farfieldeqn1}.
Let $V$ be a domain such that $\Omega \subset V$.
For a point $\bm{z}\in V$, let $g_{\bm{z}}^\epsilon$ be the regularized solution of \eqref{farfieldeqn1}.
The norm $\|g_{\bm{z}}^\epsilon\|_{L^2(\mathbb{S})}$ is relatively small when $\Omega$ is inside $B_{\bm{z}}$ and relatively large when $\Omega$ is outside $B_{\bm{z}}$ (see Theorem 3.3 in \cite{Liu2019}).
Therefore, $\|g_{\bm{z}}^\epsilon\|_{L^2(\mathbb{S})}$ can be used to characterize the location of $\Omega$.

In contrast to the classical linear sampling method \cite{Colton2013},
the kernel of $\mathcal{F}_{\bm{z}}$ in \eqref{Fzg} is the far-field pattern of $B_{\bm{z}}$ with all observation directions.
The right hand side of \eqref{farfieldeqn1}  is the measured far-field data.
This arrangement makes it possible to treat the limited aperture data.
For a fixed incident direction $\bm{d}_0$, the far-field equation \eqref{farfieldeqn1} for observation aperture $\gamma^o$ is
\begin{equation}\label{farfieldeqn2}
(\mathcal{F}_{\bm{z}} g)(\bm{\hat{x}})=\frac{1}{i k_t} u^\infty_t(\hat{\bm{x}},\bm{d}_0),\;\; \hat{\bm{x}}\in \gamma^o.
\end{equation}
Define the indicator function
\begin{equation}\label{indicator1}
I_{\bm{z}}(\bm{d}_0)=\|g_{\bm{z}}^\epsilon(\bm{d}_0)\|_{L^2(\mathbb{S})},\;\; \bm{z}\in V,
\end{equation}
where $g_{\bm{z}}^\epsilon(\bm{d}_0)$ is the regularized solution of \eqref{farfieldeqn2}.

For $u^\infty_t(\hat{\bm{x}},\bm{d})$, $(\hat{\bm{x}},\bm{d})\in \gamma^o\times \gamma^i$, the indicator is defined as
\begin{equation}\label{indicator2}
I_{\bm{z}}=\int_{\gamma^i}I_{\bm{z}}(\bm{d})d s(\bm{d}),\;\;\; \bm{z}\in V.
\end{equation}
In practice, the measured data are usually discrete
\begin{equation*}
u^\infty_t(\hat{\bm{x}}_i,\bm{d}_j),\;\; \hat{\bm{x}}_i\in\{\hat{\bm{x}}_1,\cdots,\hat{\bm{x}}_I\} \subset\mathbb{S},\;\; {\bm{d}}_j\in\{\bm{d}_1,\cdots,\bm{d}_J\} \subset\mathbb{S}.
\end{equation*}
For each $j$, let $g_{\bm{z}}^\epsilon(\bm{d}_j)$ be the solution of the far field equation
\[
(\mathcal{F}_{\bm{z}} g)(\bm{\hat{x}}_i,\bm{d}_j)=\frac{1}{i k_t} u^\infty_t(\hat{\bm{x}}_i,\bm{d}_j).
\]
Consequently, the discrete indicator is defined as
\begin{equation}\label{indicator3}
I_{\bm{z}}=\sum_{j=1}^{J}\|g_{\bm{z}}^\epsilon(\bm{d}_j)\|_{L^2(\mathbb{S})}, \;\;\; \bm{z}\in V.
\end{equation}

\subsection{ESM for IOSP-F}

Let $\mathbb{L}^2:=L^2(\mathbb{S})^2$. For $\bm{g}:=(g_p;g_s)$, $\bm{h}:=(h_p;h_s)\in \mathbb{L}^2$, define the inner product
\begin{equation*}
\langle\bm{g},\bm{h}\rangle:= \frac{\omega}{k_p}\int_{\mathbb{S}}g_p(\hat{\bm{x}}) \overline{h_p(\hat{\bm{x}}}) d s(\hat{\bm{x}})+ \frac{\omega}{k_s}\int_{\mathbb{S}}g_s(\hat{\bm{x}}) \overline{h_s(\hat{\bm{x}}}) d s(\hat{\bm{x}}).
\end{equation*}
Denote by $\bm{u}^{B_{\bm{z}}}(\hat{\bm{x}})$ the solution of \eqref{usc}-\eqref{Kupradze} with $\Omega$ replaced by $B_{\bm{z}}$ and
$\bm{u}^{B_{\bm{z}}}_\infty=({u}^{B_{\bm{z}}}_{p,\infty};{u}^{B_{\bm{z}}}_{s,\infty})$ the far-field pattern of $\bm{u}^{B_{\bm{z}}}(\hat{\bm{x}})$.
According to the Helmholtz decomposition \eqref{HelmDecomp}, $\bm{u}^{B_{\bm{z}}}=\nabla \phi+\textbf{curl}\,\psi$, where $(\phi,\psi)$ is the solution of
\begin{equation*}
\left\{\begin{array}{ll}
\Delta \phi+k_p^2 \phi=0,\quad &\text{in}\; \mathbb{R}^2\backslash \overline{B}_{\bm{z}},\\
\Delta \psi+k_s^2 \psi=0,\quad &\text{in}\; \mathbb{R}^2\backslash \overline{B}_{\bm{z}},\\
\frac{\partial \phi}{\partial \bm{\nu}}+\frac{\partial \psi}{\partial \bm{\tau}} =-\bm{\nu}\cdot \bm{u}^{inc},\quad &\text{on}\; \partial B_{\bm{z}},\\
\frac{\partial \phi}{\partial \bm{\tau}}-\frac{\partial \psi}{\partial \bm{\nu}} =-\bm{\tau}\cdot \bm{u}^{inc},\quad &\text{on}\; \partial B_{\bm{z}},\\
\lim\limits_{r\to\infty}\sqrt{r}\left(\frac{\partial \phi}{\partial r}-i k_p \phi\right)=0,\quad &r=|\bm{x}|,\\
\lim\limits_{r\to\infty}\sqrt{r}\left(\frac{\partial \psi}{\partial r}-i k_s \psi\right)=0,\quad &r=|\bm{x}|.
\end{array}\right.
\end{equation*}
Define the far-field operator $\tilde{\mathcal{F}}_{\bm{z}}: \mathbb{L}^2\rightarrow \mathbb{L}^2$ as in \cite{Liu2019}
\begin{equation}
\begin{split}
\tilde{\mathcal{F}}_{\bm{z}} \bm{g}(\bm{\hat{x}})&=\int_{\mathbb{S}} \left\{\sqrt{\frac{k_p}{\omega}} \bm{u}^{B_{\bm{z}}}_\infty(\bm{\hat{x}};\bm{d},1,0) g_p(\bm{d})+ \sqrt{\frac{k_s}{\omega}} \bm{u}^{B_{\bm{z}}}_\infty(\bm{\hat{x}};\bm{d},0,1) g_s(\bm{d})\right\} d s(\bm{d})\\
&=\int_{\mathbb{S}}\left(
\begin{array}{cc}
\sqrt{\frac{k_p}{\omega}}u^{B_{\bm{z}}}_{p,\infty}(\bm{\hat{x}};\bm{d},1,0) & \sqrt{\frac{k_s}{\omega}}u^{B_{\bm{z}}}_{p,\infty}(\bm{\hat{x}};\bm{d},0,1)\\
\sqrt{\frac{k_p}{\omega}}u^{B_{\bm{z}}}_{s,\infty}(\bm{\hat{x}};\bm{d},1,0) &\sqrt{\frac{k_s}{\omega}}u^{B_{\bm{z}}}_{s,\infty}(\bm{\hat{x}};\bm{d},0,1)
\end{array}\right)
\left(
\begin{array}{c}
g_p(\bm{d})\\
g_s(\bm{d})
\end{array}
\right)d s(\bm{d}),
\end{split}
\end{equation}
where $\bm{u}^{B_{\bm{z}}}_\infty(\hat{\bm{x}};\bm{d},\alpha_p,\alpha_s)$ denotes the far-field pattern of $B_{\bm{z}}$ due to the incident plane wave \eqref{uinc}.

For the far-field data $\bm{u}_\infty(\hat{\bm{x}},\bm{d})$, $(\hat{\bm{x}},\bm{d})\in \gamma^o\times \gamma^i$, we introduce
\begin{equation}\label{farfieldeqn3}
(\tilde{\mathcal{F}}_{\bm{z}} \bm{g})(\bm{\hat{x}})=\bm{u}_\infty(\bm{\hat{x}},\bm{d}), \;\;(\hat{\bm{x}},\bm{d})\in \gamma^o\times \gamma^i,
\end{equation}
where $\bm{g}\in \mathbb{L}^2$. From Theorem 4.2 of \cite{Liu2019}, the solution of \eqref{farfieldeqn3} has the same property as that of \eqref{farfieldeqn1}.
Similar to Section 3.1, define an indicator function
\begin{equation}\label{indicator4}
I_{\bm{z}}=\sum_{j=1}^{J}\|\bm{g}_{\bm{z}}^\epsilon(\bm{d}_j)\|_{\mathbb{L}^2}, \;\;\; \bm{z}\in V,
\end{equation}
where $\bm{g}_{\bm{z}}^\epsilon$ is the regularized solution of \eqref{farfieldeqn3}.

Since $B_{\bm{z}}$ as a disc with radius $R$ centered at $\bm{z}$,
the far-field pattern $u^{B_{\bm{z}}}_\infty(\hat{\bm{x}},k_t,\bm{d})$ and $\bm{u}^{B_{\bm{z}}}_\infty=(u^{B_{\bm{z}}}_{p,\infty};u^{B_{\bm{z}}}_{s,\infty})$
have series expansions (see, e.g., \cite{Colton2013,Liu2019}).
Given $\bm{u}_\infty(\hat{\bm{x}},\bm{d})$, $(\hat{\bm{x}},\bm{d})\in \gamma^o\times \gamma^i$,
the approximate location of $\Omega$ can be reconstructed by the ESM as follows.
\begin{enumerate}
\item For a domain $V$ such that $\Omega\subset V$, generate a set $T$ of sampling points for $V$.
\item For each $\bm{z}\in T$, calculate $u^{B_{\bm{z}}}_\infty(\hat{\bm{x}},k_t,\bm{d})$ (or $\bm{u}^{B_{\bm{z}}}_\infty=(u^{B_{\bm{z}}}_{p,\infty};u^{B_{\bm{z}}}_{s,\infty})$) for all $\hat{\bm{x}}\in \mathbb{S}$ and ${\bm{d}}\in \mathbb{S}$.
\item For each $\bm{d}_j$, solve the far-field equation \eqref{farfieldeqn2} (or \eqref{farfieldeqn3}) to obtain $g_{\bm{z}}^\epsilon(\bm{d}_j)$ (or $\bm{g}_{\bm{z}}^\epsilon(\bm{d}_j)$).

\item Calculate the indicator function ${I}_{ESM}(\bm{z})={I_{\bm{z}}}(\bm{z})/{\max_{{\bm{z}}\in T} I_{\bm{z}}(\bm{z})}$. The global minimum point $\bm{z}^*\in T$ for $I_{ESM}(\bm{z})$ is the location of $\Omega$.
\end{enumerate}

\begin{remark}
The ESM only provides the approximate location of $\Omega$.
One can use a multilevel technique to set a suitable radius of $B_{\bm{z}}$ and thus find the approximate size of $\,\Omega$.
Since the construction of initial particles proposed in Section 4 just needs an approximate location of $\,\Omega$, the above ESM is enough for the purpose of this paper.
\end{remark}

\section{Ensemble Kalman Filter}

The inverse obstacle scattering problem can be written as the statistical model to seek $\partial \Omega$ such that
\begin{equation}\label{yGOmega}
\bm{y}=\mathcal{G}(\Omega)+\bm{\eta},
\end{equation}
where $\bm{y}$ is the measured far-field data, $\mathcal{G}$ is the scattering operator and $\bm{\eta}$ is the noise.
Assume that $\bm{\eta}$ is Gaussian $\bm{\eta}\sim \mathcal{N}(0,C)$, where $C$ is the covariance matrix.
Let $\Omega$ be a starlike domain such that the boundary $\partial\Omega$ can be written as
\begin{equation}\label{Omegapz}
  \partial\Omega=r(\theta)(\cos\theta,\sin\theta)+\bm{z}= \exp(p(\theta))(\cos\theta,\sin\theta)+\bm{z},\;\; \theta\in (0,2\pi],
\end{equation}
where $p(\theta ) = \ln r(\theta )$, $0 < r(\theta ) < r_{\max}$, and $\bm{z}$ is the location of $\Omega$.

In particular, we assume that $p(\theta)$ has the following form \cite{Stuart2010,Li2020}
\begin{equation}\label{ptheta}
{p}(\theta) =\frac{{a}_0}{\sqrt{2\pi}}+\sum^{M}_{m=1}
\frac{{a}_m}{m^{s}}
\frac{\cos (m\theta)}{\sqrt{\pi}}+
\frac{{b}_m}{m^{s}}
\frac{\sin (m\theta)}{\sqrt{\pi}},
\end{equation}
where $s$ is a smoothing parameter.
Let $\bm{q}:=(a_0,a_1,a_2,\cdots,a_m,b_m)^\top$.
The inverse problem is to determine  $\bm{\xi}$ from $\bm{y}$ such that
\begin{equation}\label{yGxi}
\bm{y}=\mathcal{G}(\bm{\xi})+\bm{\eta},\;\;\; \bm{\eta}\sim \mathcal{N}(0,C),
\end{equation}
where $\bm{\xi}:=(\bm{q},\bm{z})^\top =(a_0,a_1,a_2,\cdots,a_m,b_m,z_1,z_2)^\top\in \mathbb{R}^{2m+3}$.

To solve the inverse problem by the Kalman filter (KF),
we construct an artificial dynamic system as follows.
Let $Z:=\mathbb{R}^{2M+3}\times \mathbb{C}^N$ and $\bm{\phi}=(\bm{\xi}, \bm{\omega})^{\top}\in Z$. We define $\Psi: Z\rightarrow Z$ by
\begin{equation*}
\Psi:
\begin{pmatrix}
\bm{\xi} \\ \bm{\omega}
\end{pmatrix}
\rightarrow
\begin{pmatrix}
\bm{\xi} \\ \mathcal{G}(\bm{\xi})
\end{pmatrix}.
\end{equation*}
Define $H: Z\rightarrow \mathbb{C}^N$ such that $H=(0,I)$. Introduce the artificial dynamic system
\begin{equation}\label{dynamicsystem}
\begin{array}{l}
  \bm{\phi}_{n+1}=\Psi(\bm{\phi}_{n}),  \\
  \bm{y}_{n+1}=H \bm{\phi}_{n+1} +\bm{\eta}_{n+1},
\end{array}
\end{equation}
where $\{\bm{\eta}_n\}_{n\in \mathbb{Z}^{+}}$ is an i.i.d.(independent and identically distributed) Gaussian sequence, i.e., $\bm{\eta}_n \sim \mathcal{N}(0,C)$.

The formulation of the Kalman filter can be interpreted in the framework of either optimization or Bayesian inference.
We shall briefly discuss the Bayesian perspective (see, e.g., \cite{Kaipio2005}) and refer the readers to \cite{Iglesias2013,Chada2020} for the optimization perspective.
In the Bayesian framework, all variables in \eqref{dynamicsystem} are treated as random variables.
The target of the filter is to extract information from the distribution of $\bm{\phi}_n$ conditioned on the data $\bm{\mathfrak{g}}_n:=\{\bm{y}_n\}_{n=1}^{N}$, i.e., $\pi(\bm{\phi}_n|\bm{\mathfrak{g}}_n)$.
This can be done by using a sequential procedure consisting the following two steps.
The first step is prediction. Given $\mathfrak{g}_{n-1}$, one computes the distribution $\pi(\bm{\phi}_n|\bm{\mathfrak{g}}_{n-1})$ according to
\begin{equation*}
\pi(\bm{\phi}_n|\bm{\mathfrak{g}}_{n-1})=\int \pi(\bm{\phi}_n|\bm{\phi}_{n-1}) \pi(\bm{\phi}_{n-1}|\bm{\mathfrak{g}}_{n-1}) d \bm{\phi}_n.
\end{equation*}
The second step is analysis. For the new observations $\bm{y}_n$, the distribution $\pi(\bm{\phi}_n|\bm{\mathfrak{g}}_{n})$ is
\begin{equation*}
\pi(\bm{\phi}_n|\bm{\mathfrak{g}}_{n})= \frac{\pi(\bm{y}_n|\bm{\phi}_{n}) \pi(\bm{\phi}_n|\bm{\mathfrak{g}}_{n-1})} {\pi(\bm{y}_n|\bm{\mathfrak{g}}_{n-1})},
\end{equation*}
where
\begin{equation*}
\pi(\bm{y}_n|\bm{\mathfrak{g}}_{n-1})=\int \pi(\bm{y}_n|\bm{\phi}_{n}) \pi(\bm{\phi}_{n}|\bm{\mathfrak{g}}_{n-1}) d \bm{\phi}_n.
\end{equation*}

When the system \eqref{dynamicsystem} is linear, the filtered distribution $\pi(\bm{\phi}_n|\bm{\mathfrak{g}}_{n})$ is Gaussian.
The mean and covariance are given by the Kalman equations (Theorem 4.3 of \cite{Kaipio2005})
\begin{equation}\label{KFupdate}
\begin{split}
&\bm{\phi}_{n|n-1}=F \bm{\phi}_{n-1|n-1}, \\
&\bm{\phi}_{n|n}=\bm{\phi}_{n|n-1}+\Xi_{n} (\bm{y}_n-H\bm{\phi}_{n|n-1}),
\end{split}\;\;\;\;\;\;\;\;\;\;
\begin{split}
& \Gamma_{n|n-1}=F \Gamma_{n-1|n-1} F^T,\\ &\Gamma_{n|n}=(1-\Xi_{n}H)\Gamma_{n|n-1},
\end{split}
\end{equation}
where $F$ denotes the matrix for the mapping $\Psi: Z\rightarrow Z$, $\bm{\phi}_{n|l}=\mathbb{E}(\bm{\phi}_n|\bm{\mathfrak{g}}_l)$, $\Gamma_{n|l}=\text{Cov}(\bm{\phi}_n|\bm{\mathfrak{g}}_l)$, and $\Xi_n$ is the Kalman gain matrix given by
\begin{equation}\label{Kalmangain}
\Xi_n=\Gamma_{n|n-1}H^T(H\Gamma_{n|n-1}H^T+C)^{-1}.
\end{equation}

For nonlinear systems, the filtering distribution $\pi(\bm{\phi}_n|\bm{\mathfrak{g}}_{n})$ is no longer Gaussian.
However, the framework can be generalized by approximating the distribution $\pi(\bm{\phi}_n|\bm{\mathfrak{g}}_{n})$ through its Gaussian approximation $\pi_G(\bm{\phi}_n|\bm{\mathfrak{g}}_{n})$.
The ensemble Kalman filter (EnKF) is a powerful tool to deal with both linear and nonlinear systems. Compared with the extended Kalman filter (EKF), the EnKF does not need to compute the Fr\'{e}chet derivative of the forward operator.

For the EnKF, the true mean and covariance appearing in the KF are estimated by an ensemble of particles $\{\bm{\phi}_n^{(j)}\}_{j=1}^{J}$, and the propagation of these particles $\{\bm{\phi}_n^{(j)}\}_{j=1}^{J}$ follows the standard Kalman equations \eqref{KFupdate}.
The initial ensemble $\{\bm{\phi}_0^{(j)}\}_{j=1}^{J}$ is
\begin{equation}\label{initialensemble}
\bm{\phi}_0^{(j)}=
\begin{pmatrix}
\bm{\xi}_0^{(j)}\\
\mathcal{G}(\bm{\xi}_0^{(j)})
\end{pmatrix},
\end{equation}
where $\bm{\xi}_0^{(j)}$ is generated according to the prior distribution.
Assume $a_0,a_m,b_m\sim \mathcal{N}(0,1)$, $z_i\sim \mathcal{N}(z_i^*,1)$, $m=1,\cdots, M$, $i=1,2$, where $z_i^*$ is the approximate location of $\Omega$ obtained by the ESM.
Let $\mathcal{A}:= \text{span}\{{\bm{\xi}_0^{(j)}}\}_{j=1}^{J}$.
From the invariance subspace property of the EnKF \cite{Iglesias2013}, the inversion solution of \eqref{yGxi} still lies in the subspace of $\mathcal{A}$.
Given the initial ensemble, the procedure of the propagation of each ensemble particle $\{\bm{\phi}_n^{(j)}\}_{j=1}^{J}$ is as follows.

(1) Prediction step. Map forward the current ensemble $\{\bm{\phi}_n^{(j)}\}_{j=1}^{J}$ according to the artificial dynamic system
\begin{equation*}
\hat{\bm{\phi}}_{n+1}^{(j)}=\Psi({\bm{\phi}}_n^{(j)}),
\end{equation*}
and calculate the sample mean and covariance
\begin{equation*}
\bar{\bm{\phi}}_{n+1}=\frac{1}{J} \sum_{j=1}^{J}\hat{\bm{\phi}}_{n+1}^{(j)},\;\;\;
{\Gamma}_{n+1}=\frac{1}{J-1}\sum_{j=1}^{J}(\hat{\bm{\phi}}_{n+1}^{(j)} -\bar{\bm{\phi}}_{n+1})\otimes (\hat{\bm{\phi}}_{n+1}^{(j)} -\bar{\bm{\phi}}_{n+1}),
\end{equation*}
where $J$ is the number of particles. The mean $\bar{\bm{\phi}}_{n+1}$ and covariance  ${\Gamma}_{n+1}$ have the following block structures
\begin{equation*}
\bar{\bm{\phi}}_{n+1}=
\begin{pmatrix}
\bar{\bm{\xi}}_{n+1} \\
\bar{\bm{\omega}}_{n+1}
\end{pmatrix},\;\;\;\;\;
{\Gamma}_{n+1}=
\begin{pmatrix}
{\Gamma}_{n+1}^{{\xi}{\xi}} & {\Gamma}_{n+1}^{{\xi}{\omega}} \\
({\Gamma}_{n+1}^{{\xi}{\omega}})^T & {\Gamma}_{n+1}^{{\omega}\bm{\omega}}
\end{pmatrix}.
\end{equation*}
In the above equation,
\begin{equation*}
\bar{\bm{\xi}}_{n+1}=\frac{1}{J} \sum_{j=1}^{J}\hat{\bm{\xi}}_{n+1}^{(j)}=\frac{1}{J} \sum_{j=1}^{J}{\bm{\xi}}_{n}^{(j)},\;\;\; \bar{\bm{\omega}}_{n+1}=\frac{1}{J} \sum_{j=1}^{J}\hat{\bm{\omega}}_{n+1}^{(j)}=\frac{1}{J} \sum_{j=1}^{J}\mathcal{G}({\bm{\xi}}_{n}^{(j)}),
\end{equation*}
and
\begin{equation*}
\begin{aligned}
&{\Gamma}_{n+1}^{{\xi}{\xi}}=\frac{1}{J-1}\sum_{j=1}^{J}(\hat{\bm{\xi}}_{n+1}^{(j)} -\bar{\bm{\xi}}_{n+1})\otimes (\hat{\bm{\xi}}_{n+1}^{(j)} -\bar{\bm{\xi}}_{n+1}),\\[-1ex]
&{\Gamma}_{n+1}^{{\xi}{\omega}}=\frac{1}{J-1}\sum_{j=1}^{J}(\hat{\bm{\xi}}_{n+1}^{(j)} -\bar{\bm{\xi}}_{n+1})\otimes (\hat{\bm{\omega}}_{n+1}^{(j)} -\bar{\bm{\omega}}_{n+1}),\\[-1ex]
&{\Gamma}_{n+1}^{{\omega}{\omega}}=\frac{1}{J-1}\sum_{j=1}^{J}(\hat{\bm{\omega}}_{n+1}^{(j)} -\bar{\bm{\omega}}_{n+1})\otimes (\hat{\bm{\omega}}_{n+1}^{(j)} -\bar{\bm{\omega}}_{n+1}).
\end{aligned}
\end{equation*}
(2) Analysis step.
Calculate the Kalman gain matrix
\begin{equation}\label{Kalmangain}
\Xi_{n+1}=\Gamma_{n+1} H^T (H \Gamma_{n+1} H^T +C)^{-1},
\end{equation}
and update each ensemble
\begin{equation}\label{Kalmanupdate}
\bm{\phi}_{n+1}^{(j)}=I\hat{\bm{\phi}}_{n+1}^{(j)} +\Xi_{n+1} (\bm{y}_{n+1}^{(j)}-H \hat{\bm{\phi}}_{n+1}^{(j)}).
\end{equation}
Due to the structure of $H$, \eqref{Kalmanupdate} is equivalent to
\begin{equation}\label{updatexiomega}
\begin{aligned}
\bm{\xi}_{n+1}^{(j)}&={\bm{\xi}}_{n}^{(j)}+\Gamma_{n+1}^{\xi\omega} (\Gamma_{n+1}^{\omega\omega}+C)^{-1}\left(\bm{y}_{n+1}^{(j)}- \mathcal{G}(\bm{\xi}_{n}^{(j)})\right),\\
\bm{\omega}_{n+1}^{(j)}&=\mathcal{G}({\bm{\xi}}_{n}^{(j)})+ \Gamma_{n+1}^{\omega\omega}(\Gamma_{n+1}^{\omega\omega}+C)^{-1}
\left(\bm{y}_{n+1}^{(j)}- \mathcal{G}(\bm{\xi}_{n}^{(j)})\right).
\end{aligned}
\end{equation}
The EnKF estimator of the inverse problem is obtained by averaging over the particles $\{\bm{\xi}_{n+1}^{(j)}\}_{j=1}^J$
\begin{equation}\label{EnKFestimator}
\bm{\xi}_{n+1}=\frac{1}{J}\sum_{j=1}^{J} \bm{\xi}_{n+1}^{(j)}.
\end{equation}

In the numerical experiments, one can use a blocking strategy by updating the two components $\bm{q}_n^{(j)}$ and $\bm{z}_n^{(j)}$ of $\bm{\xi}_{n}^{(j)}$ separately:
\begin{equation}\label{updataqz}
\begin{aligned}
\bm{q}_{n+1}^{(j)}&={\bm{q}}_{n}^{(j)}+\Gamma_{n+1}^{q\omega} (\Gamma_{n+1}^{\omega\omega}+C)^{-1}\left(\bm{y}_{n+1}^{(j)}- \mathcal{G}(\bm{\xi}_{n}^{(j)})\right), \\
\bm{z}_{n+1}^{(j)}&=\bm{z}_{n}^{(j)}+\Gamma_{n+1}^{z\omega} (\Gamma_{n+1}^{\omega\omega}+C)^{-1}\left(\bm{y}_{n+1}^{(j)}- \mathcal{G}(\bm{\xi}_{n}^{(j)})\right),
\end{aligned}
\end{equation}
where
\begin{equation*}
\begin{aligned}
&{\Gamma}_{n+1}^{{q}{\omega}}= \frac{1}{J-1}\sum_{j=1}^{J}(\hat{\bm{q}}_{n+1}^{(j)} -\bar{\bm{q}}_{n+1})\otimes (\hat{\bm{\omega}}_{n+1}^{(j)} -\bar{\bm{\omega}}_{n+1}),\\[-1ex]
&{\Gamma}_{n+1}^{{z}{\omega}}=
\frac{1}{J-1}\sum_{j=1}^{J}(\hat{\bm{z}}_{n+1}^{(j)} -\bar{\bm{z}}_{n+1})\otimes (\hat{\bm{\omega}}_{n+1}^{(j)} -\bar{\bm{\omega}}_{n+1}).
\end{aligned}
\end{equation*}
Consequently, the EnKF estimator of the inverse problem is given by
\begin{equation}\label{EnKFestimator2}
\bm{q}_{n+1}=\frac{1}{J}\sum_{j=1}^{J} \bm{q}_{n+1}^{(j)},\;\;\;\;\;
\bm{z}_{n+1}=\frac{1}{J}\sum_{j=1}^{J} \bm{z}_{n+1}^{(j)}.
\end{equation}

\section{Numerical Experiments}

We present some examples to show the performance of the proposed approach.
Let $\omega=\pi$, $\lambda=2$, $\mu=1$. The incident field is the plane compressional wave
$
\bm{u}^{inc}(\bm{x})=\bm{d}e^{i k_p \bm{x}\cdot \bm{d}}.
$
The forward problem is solved by the Nystr\"{o}m method \cite{Dong2019} on a finer mesh (128 equidistant points on $\partial\Omega$).
Then 3\% relative error is added to the computed far-field data, which is the simulated measured data.
In the inversion stage, a coarser mesh is used (64 equidistant points on $\partial \Omega$).

In the first step, we set $V=[-5,5]\times[-5,5]$ and the sampling points are given by
\begin{equation*}
T:=\{(-5+0.1 k,-5+0.1 l),\;\; k,l=0,1,\cdots,100\}.
\end{equation*}
The radius of the reference disc is $R = 1$. For each $\bm{z}\in T$, the far-field equations \eqref{farfieldeqn2} and \eqref{farfieldeqn3}
are solved by the Tikhonov regularization with a fixed regularization parameter $10^{-5}$.
In the EnKF, we set $s=1.2$ in \eqref{ptheta} and $M = 6$ in the Fourier expansion \eqref{ptheta}.
The maximum number of iterations and particle size are set to 30 and 500, respectively.

Five different observation apertures are
\begin{eqnarray*}
&& \gamma_1^o=\{(\cos \theta, \sin \theta) | \; \theta \in [0,2\pi]\}, \\
&& \gamma_2^o=\{(\cos \theta, \sin \theta) | \;\theta \in [0,\pi]\}, \\
&& \gamma_3^o=\{(\cos \theta, \sin \theta) | \;\theta \in [0,\pi/2]\},\\
&& \gamma_4^o=\{ (\cos \theta, \sin \theta) | \;\theta \in [0,\pi/2]\cup[\pi,3\pi/2]\}, \\
&& \gamma_5^o=\{(\cos \theta, \sin \theta) | \;\theta \in [0,\pi/4]\cup [\pi,5\pi/4]\}.
\end{eqnarray*}
The incident apertures are
\begin{eqnarray*}
&& \gamma_1^i = \{ (1/2,\sqrt{3}/2) \}, \\
&& \gamma_2^i = \{(\cos\theta,\sin \theta) |\; \theta=\{0,\pi/8,\pi/4,3\pi/8,\pi/2\}\}.
\end{eqnarray*}

For the difference between the exact and reconstructed boundaries, we use the Hausdorff distance defined as
\begin{equation*}
d_H(\partial \Omega_{\text{1}},\partial \Omega_{\text{2}}):=\max \Big\{\sup_{\bm{x}\in \partial\Omega_{\text{1}}} \inf_{\bm{y}\in \partial\Omega_{\text{2}}}|\bm{x}-\bm{y}|,\; \sup_{\bm{y}\in \partial\Omega_{\text{2}}} \inf_{\bm{x}\in \partial\Omega_{\text{1}}} |\bm{y}-\bm{x}| \Big\}.
\end{equation*}

\begin{figure}[!h]
\centering
\begin{minipage}{0.32\linewidth}
\centering
\includegraphics[height=40mm]{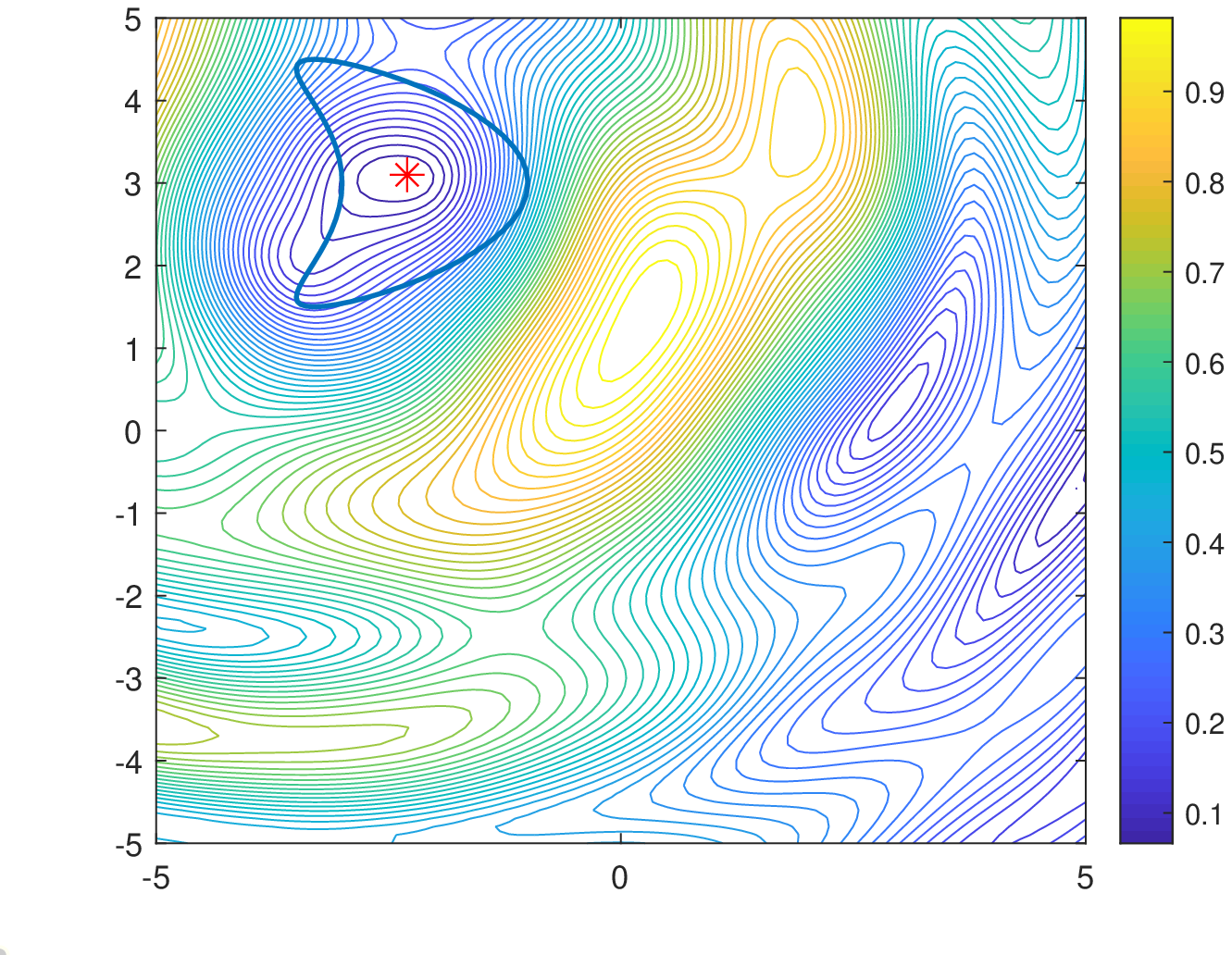}
\end{minipage}
\hfill
\begin{minipage}{0.32\linewidth}
\centering
\includegraphics[height=40mm]{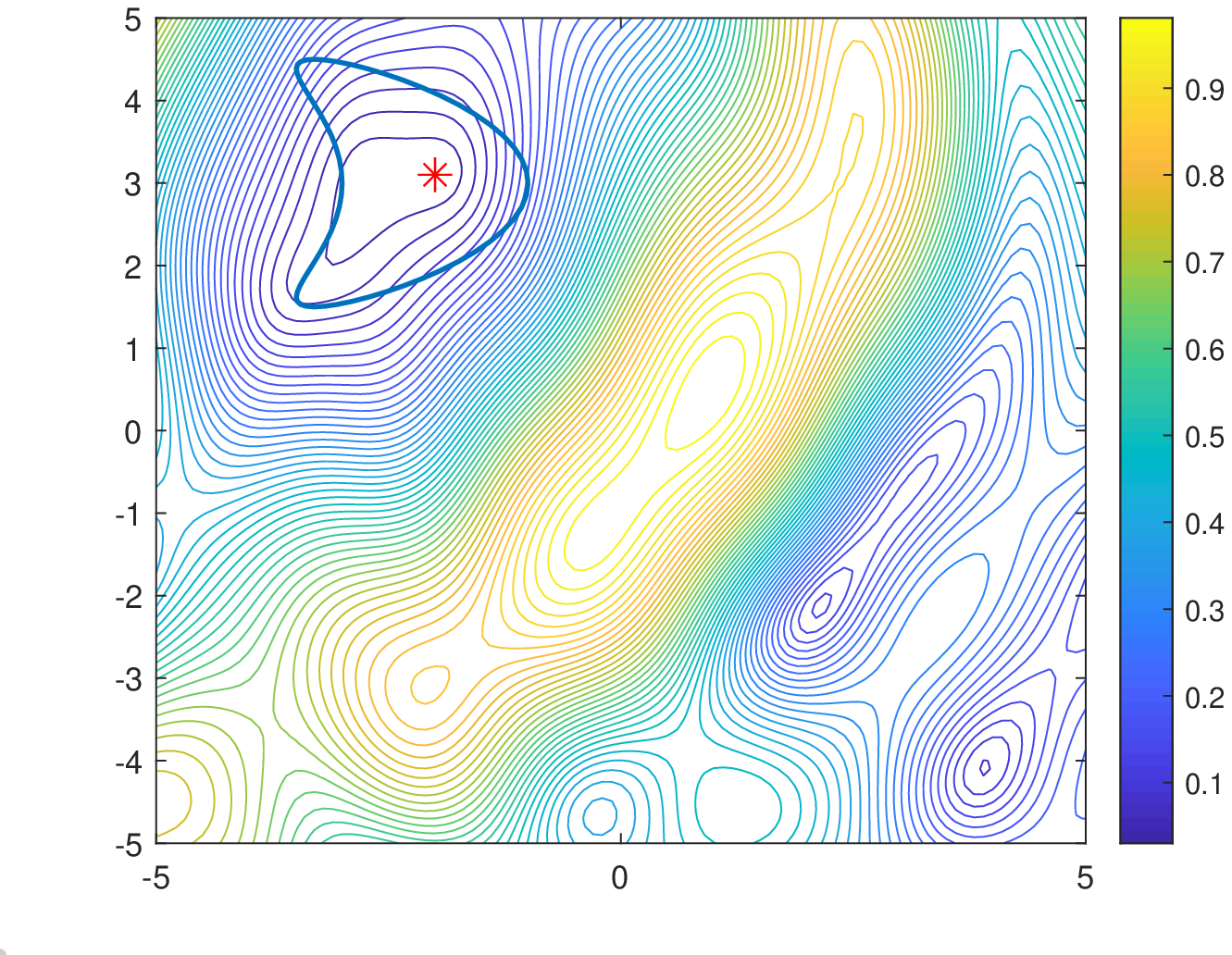}
\end{minipage}
\hfill
\begin{minipage}{0.32\linewidth}
\centering
\includegraphics[height=40mm]{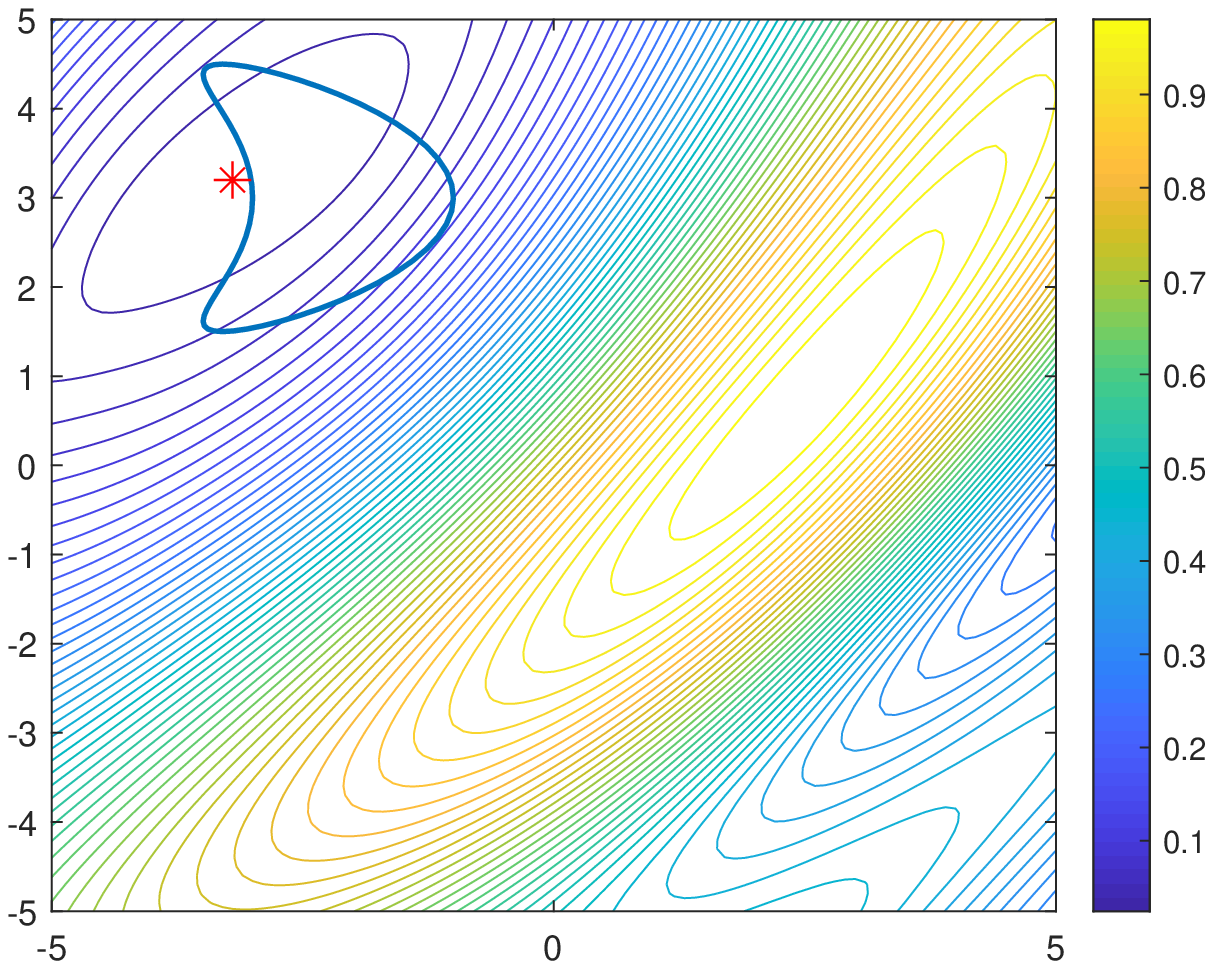}
\end{minipage}

\begin{minipage}{0.32\linewidth}
\centering
\includegraphics[height=40mm]{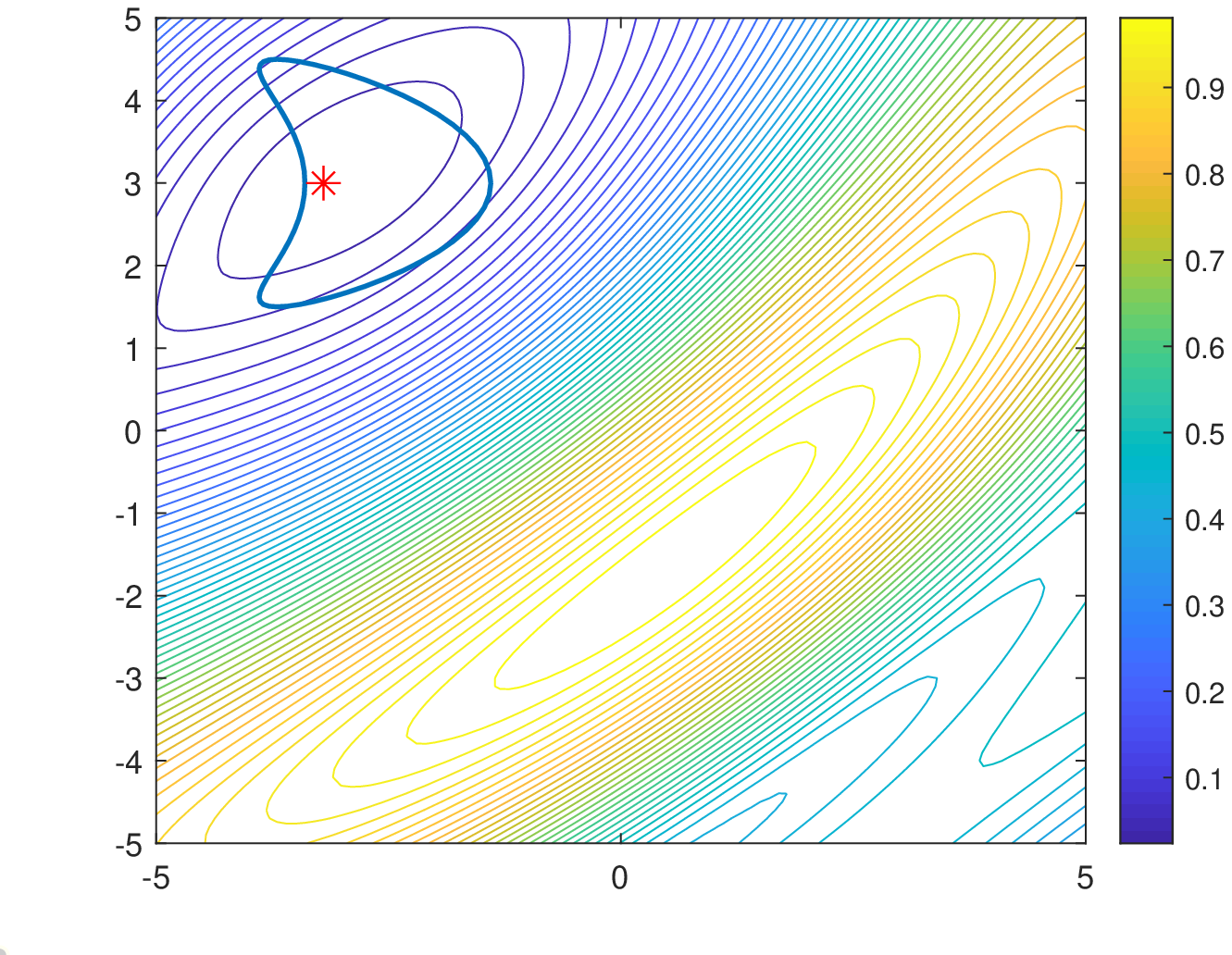}
\end{minipage}
\hfill
\begin{minipage}{0.32\linewidth}
\centering
\includegraphics[height=40mm]{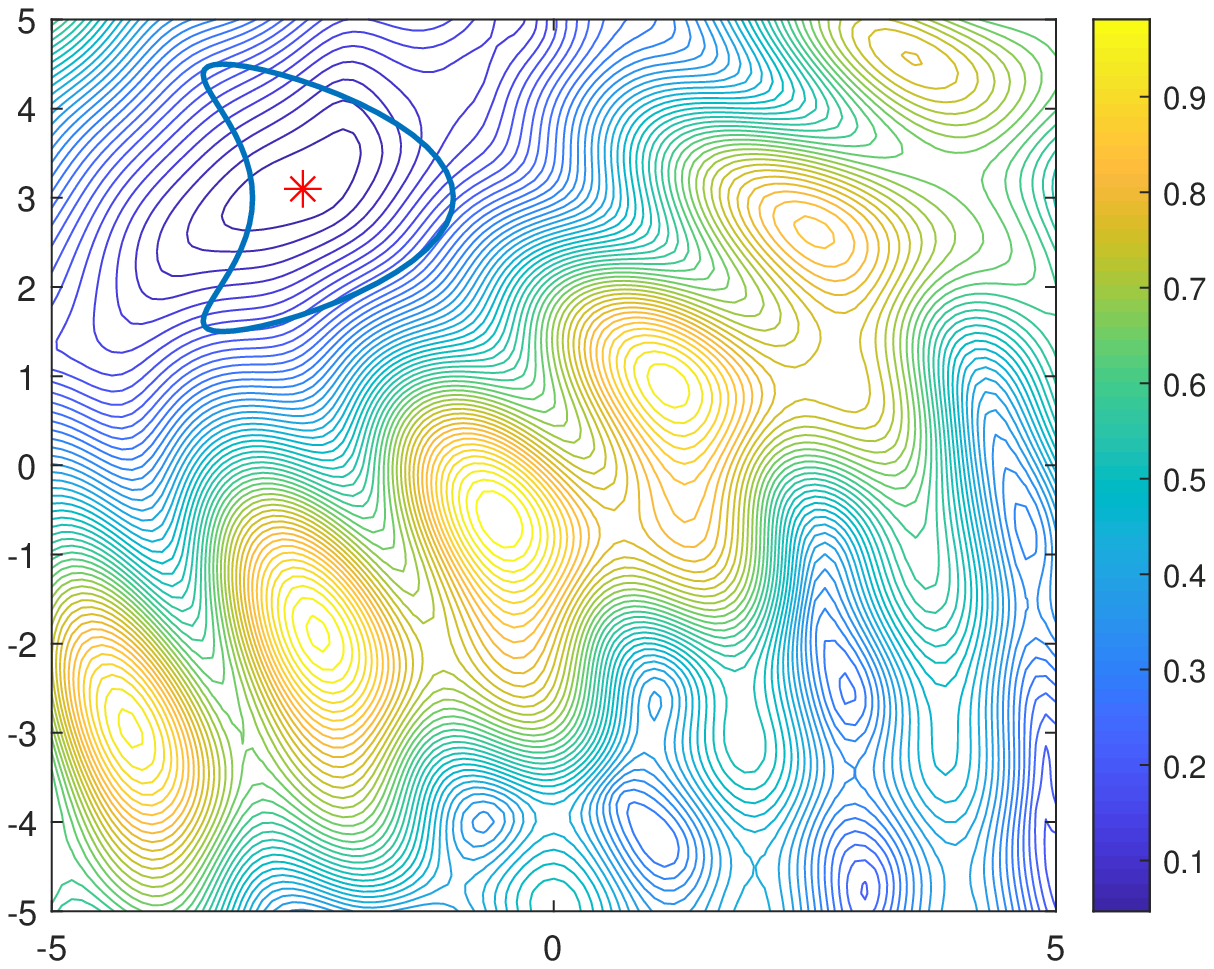}
\end{minipage}
\hfill
\begin{minipage}{0.32\linewidth}
\centering
\includegraphics[height=40mm]{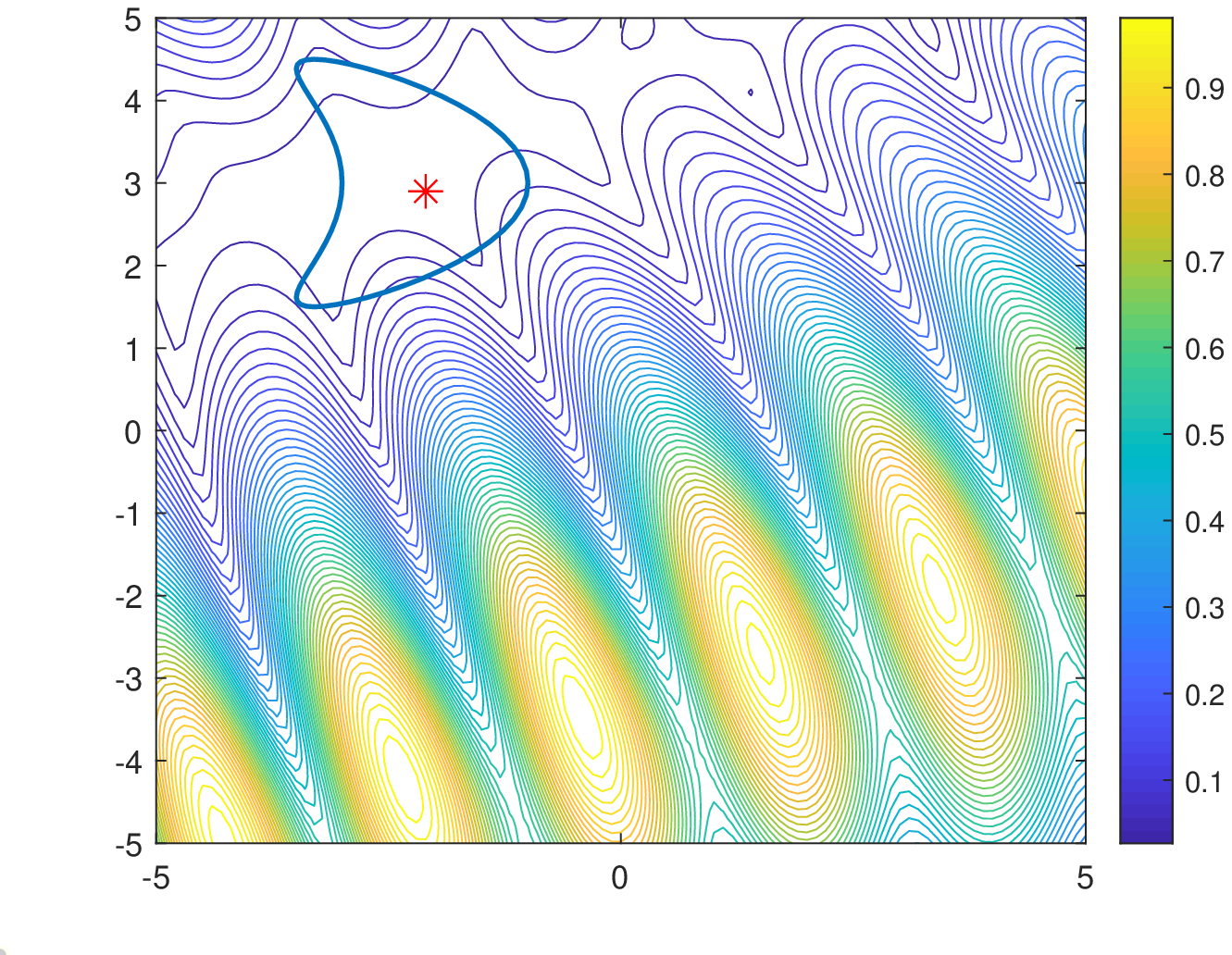}
\end{minipage}
\caption{Contour plots of the indicator function $I_{ESM}(\bm{z})$. The asterisk `*' indicates the reconstructed location by the ESM. Top row: $\gamma^i=\gamma_1^i$, from left to right: $\gamma^o=\gamma_1^o$, $\gamma_2^o$, $\gamma_3^o$. Bottom row: $\gamma^i=\gamma_2^i$, from left to right: $\gamma^o=\gamma_3^o$, $\gamma_4^o$, $\gamma_5^o$.}\label{ESM_kite}
\end{figure}

\begin{figure}[!h]
\centering
\begin{minipage}{0.32\linewidth}
\centering
\includegraphics[height=40mm]{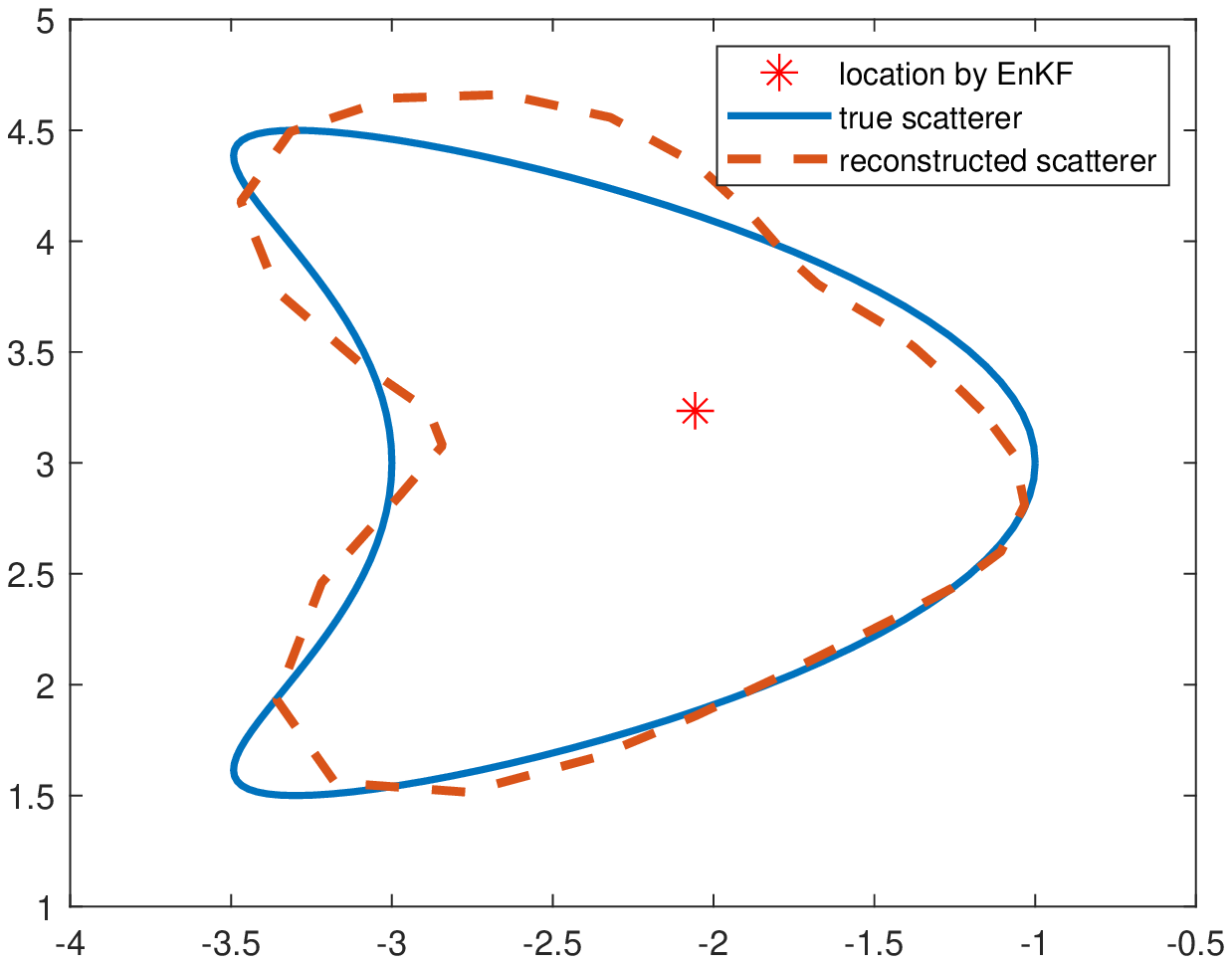}
\end{minipage}
\hfill
\begin{minipage}{0.32\linewidth}
\centering
\includegraphics[height=40mm]{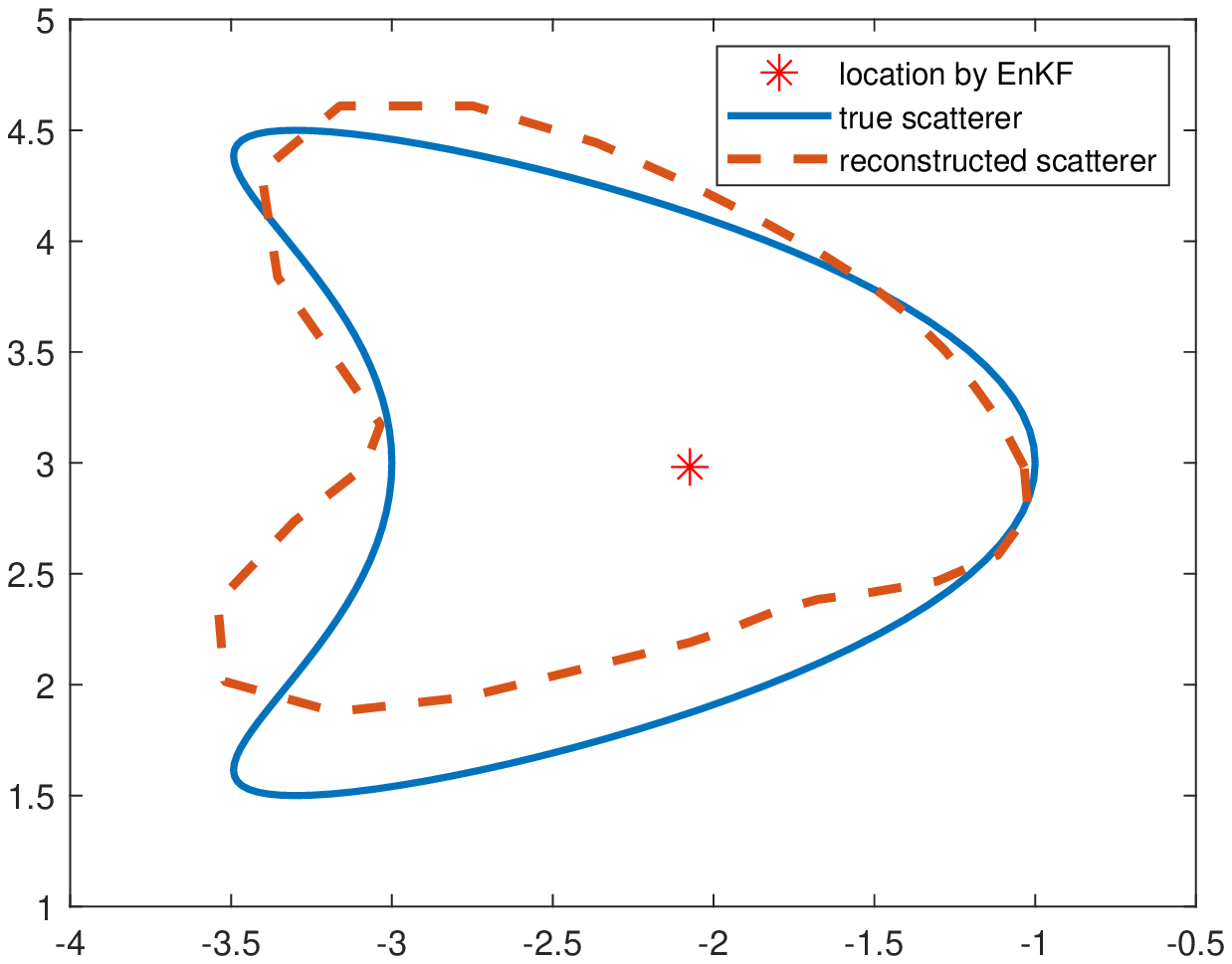}
\end{minipage}
\hfill
\begin{minipage}{0.32\linewidth}
\centering
\includegraphics[height=40mm]{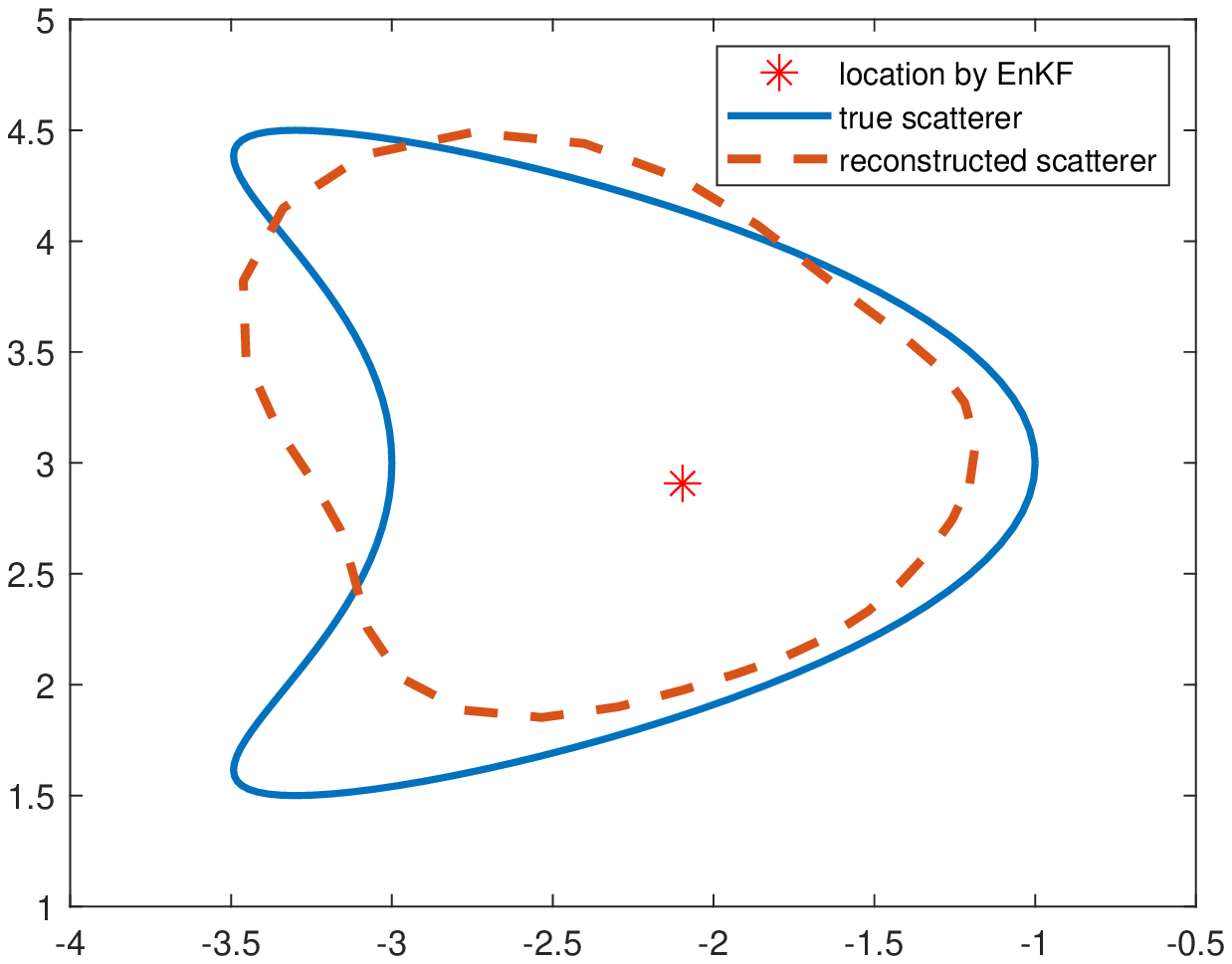}
\end{minipage}

\begin{minipage}{0.32\linewidth}
\centering
\includegraphics[height=40mm]{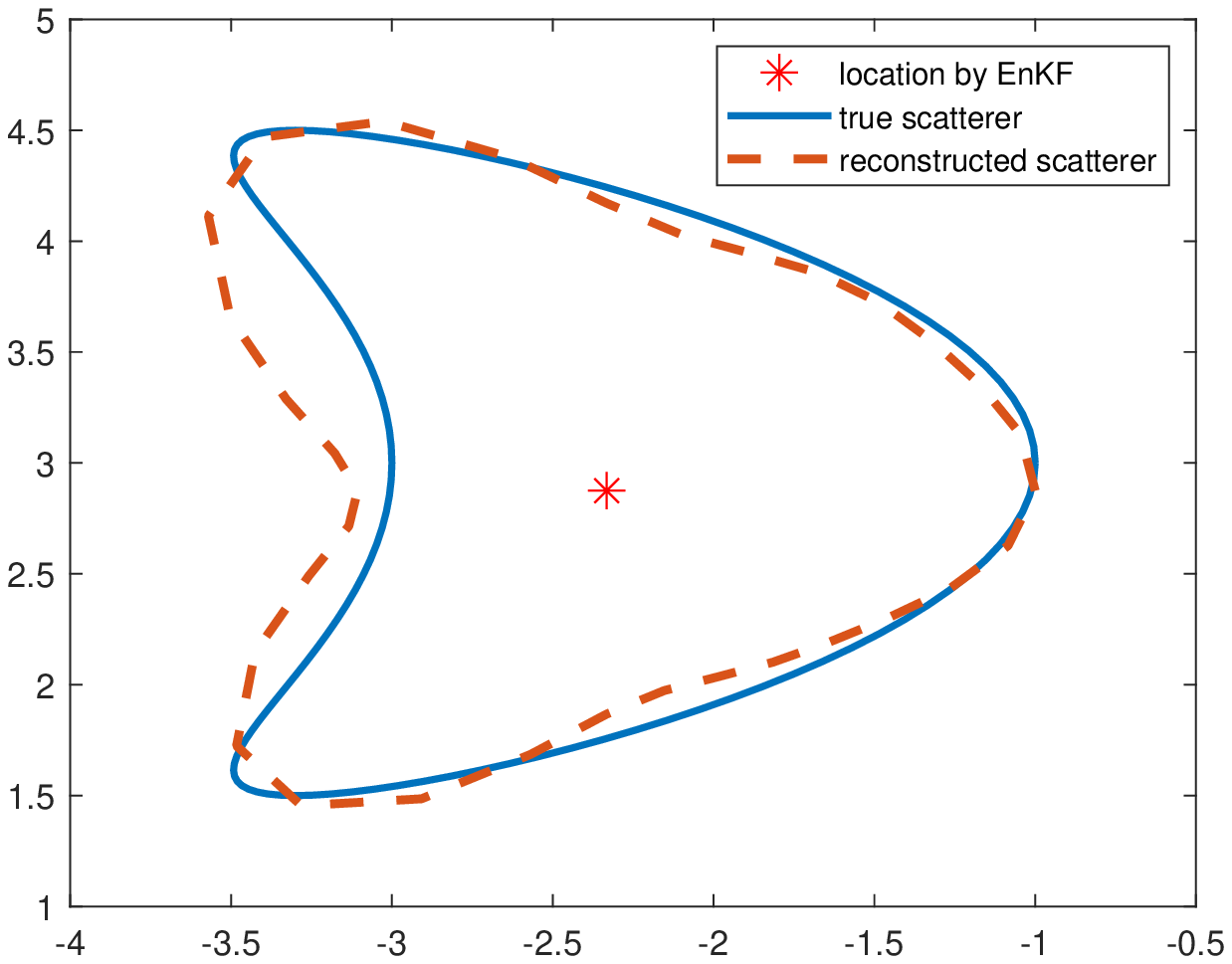}
\end{minipage}
\hfill
\begin{minipage}{0.32\linewidth}
\centering
\includegraphics[height=40mm]{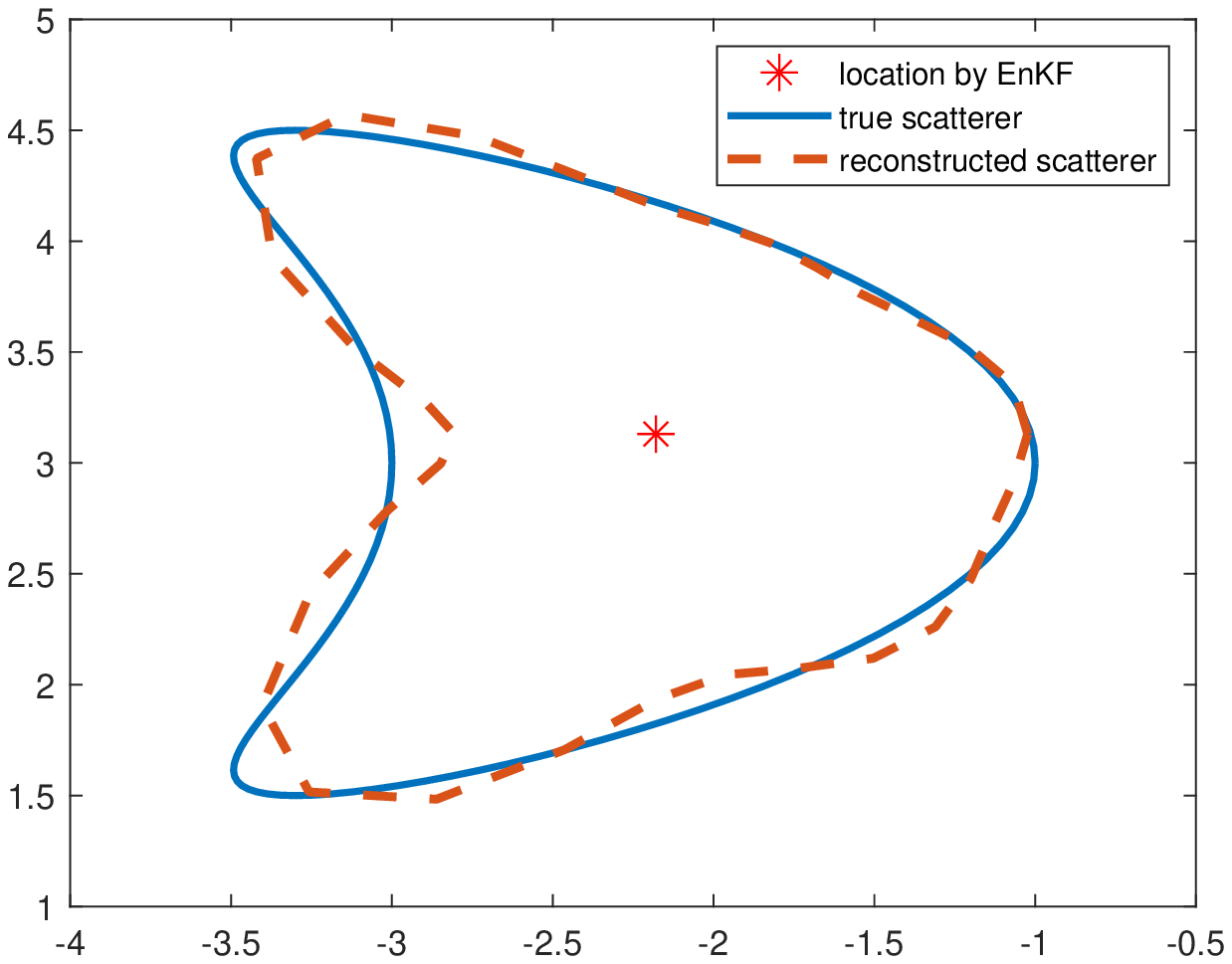}
\end{minipage}
\hfill
\begin{minipage}{0.32\linewidth}
\centering
\includegraphics[height=40mm]{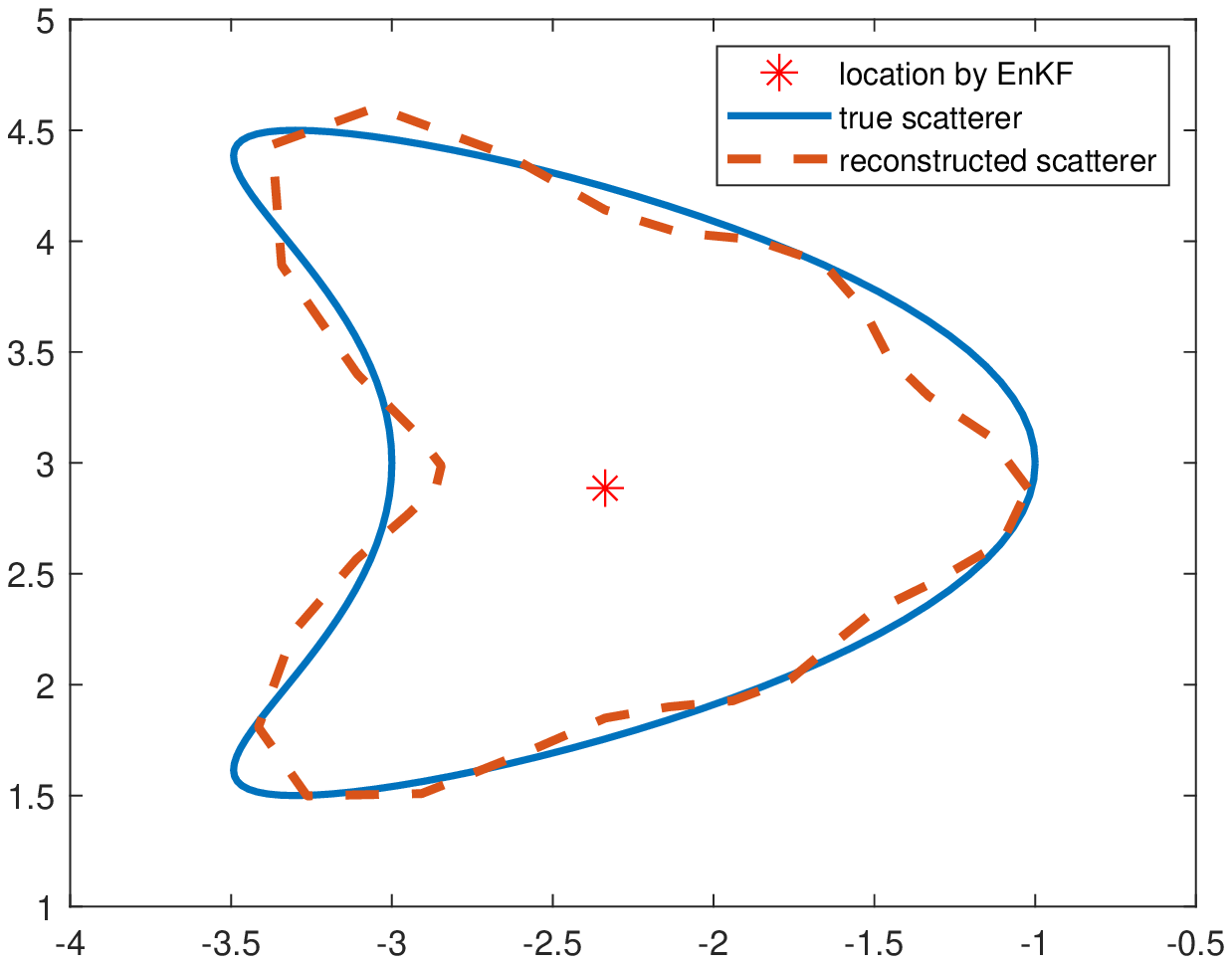}
\end{minipage}
\caption{Boundary reconstructions by the EnKF. The asterisk `*' represents the refined location by the EnKF. Top row: $\gamma^i=\gamma_1^i$, from left to right: $\gamma^o=\gamma_1^o$, $\gamma_2^o$, $\gamma_3^o$. Bottom row: $\gamma^i=\gamma_2^i$, from left to right: $\gamma^o=\gamma_3^o$, $\gamma_4^o$, $\gamma_5^o$.}\label{EnKF_kite}
\end{figure}

\begin{figure}[!h]
\centering
\begin{minipage}{0.45\linewidth}
\centering
\includegraphics[height=50mm]{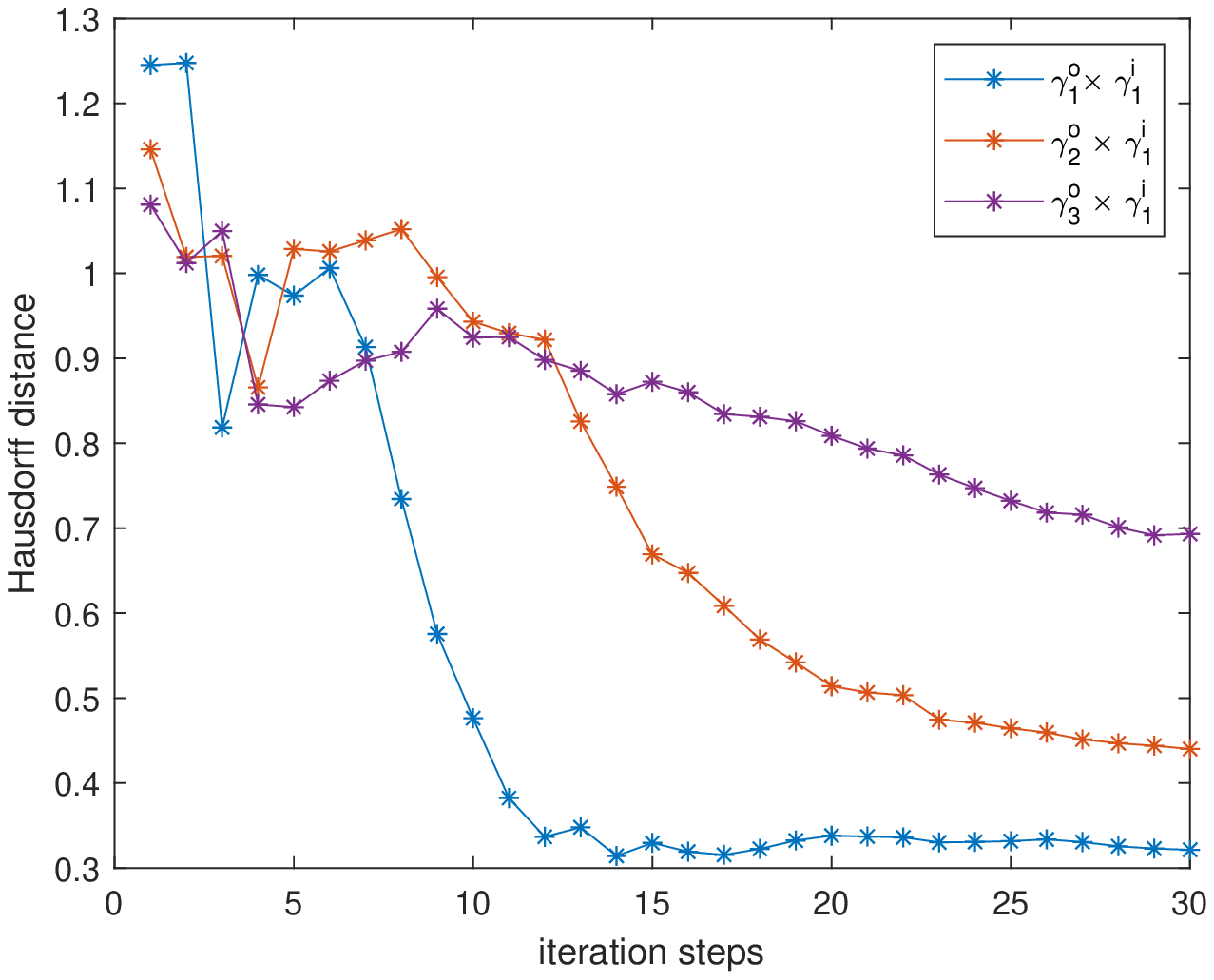}
\end{minipage}
\begin{minipage}{0.45\linewidth}
\centering
\includegraphics[height=50mm]{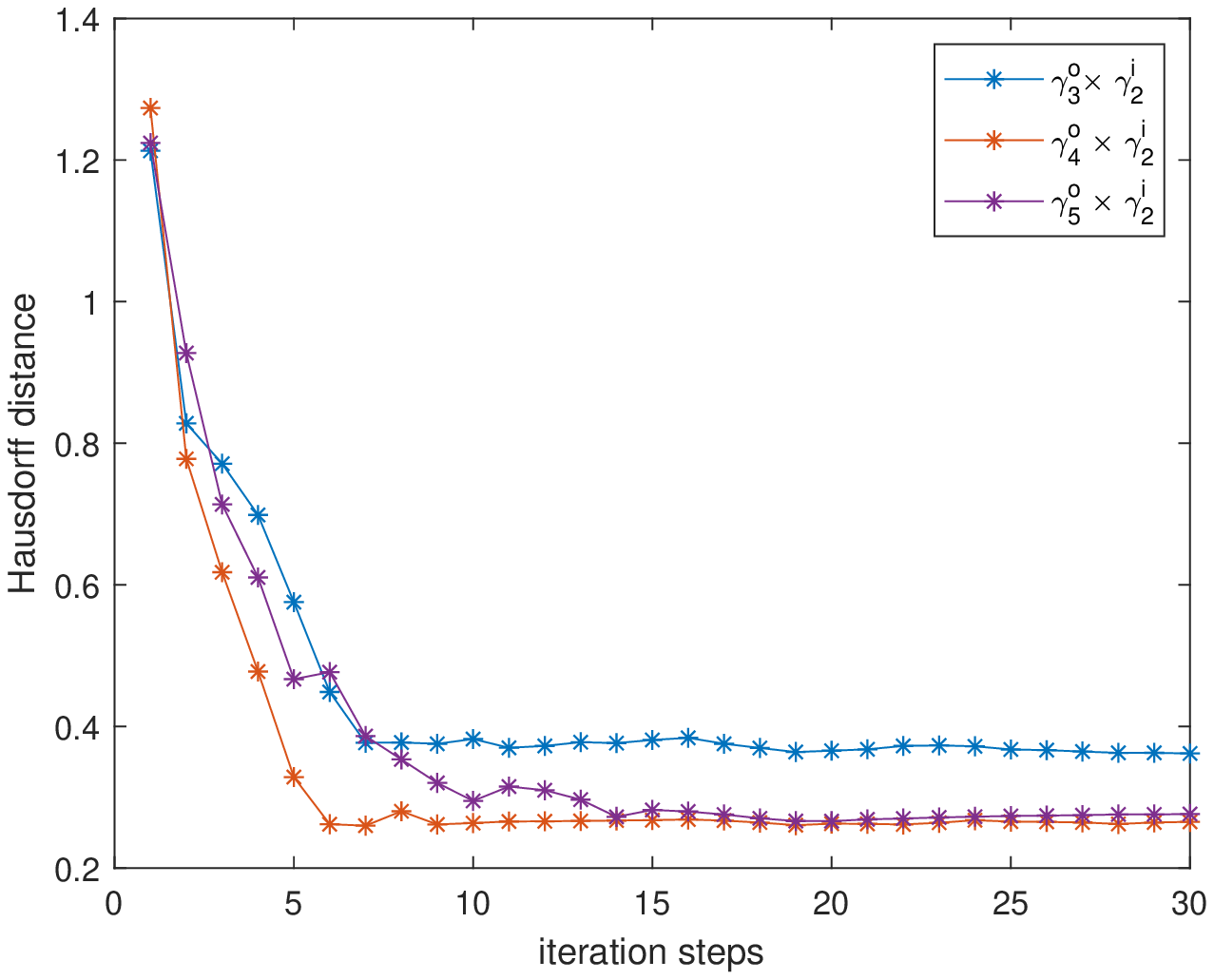}
\end{minipage}
\caption{Hausdorff distance $d_H(\partial \Omega_{\text{exact}},\partial \Omega_{\text{inv}})$ with respect to the iteration steps. Left: $\gamma^i_1$. Right: $\gamma_2^i$.
}
\label{dH_kite}
\end{figure}

\begin{figure}[!h]
\centering
\begin{minipage}{0.45\linewidth}
\centering
\includegraphics[height=50mm]{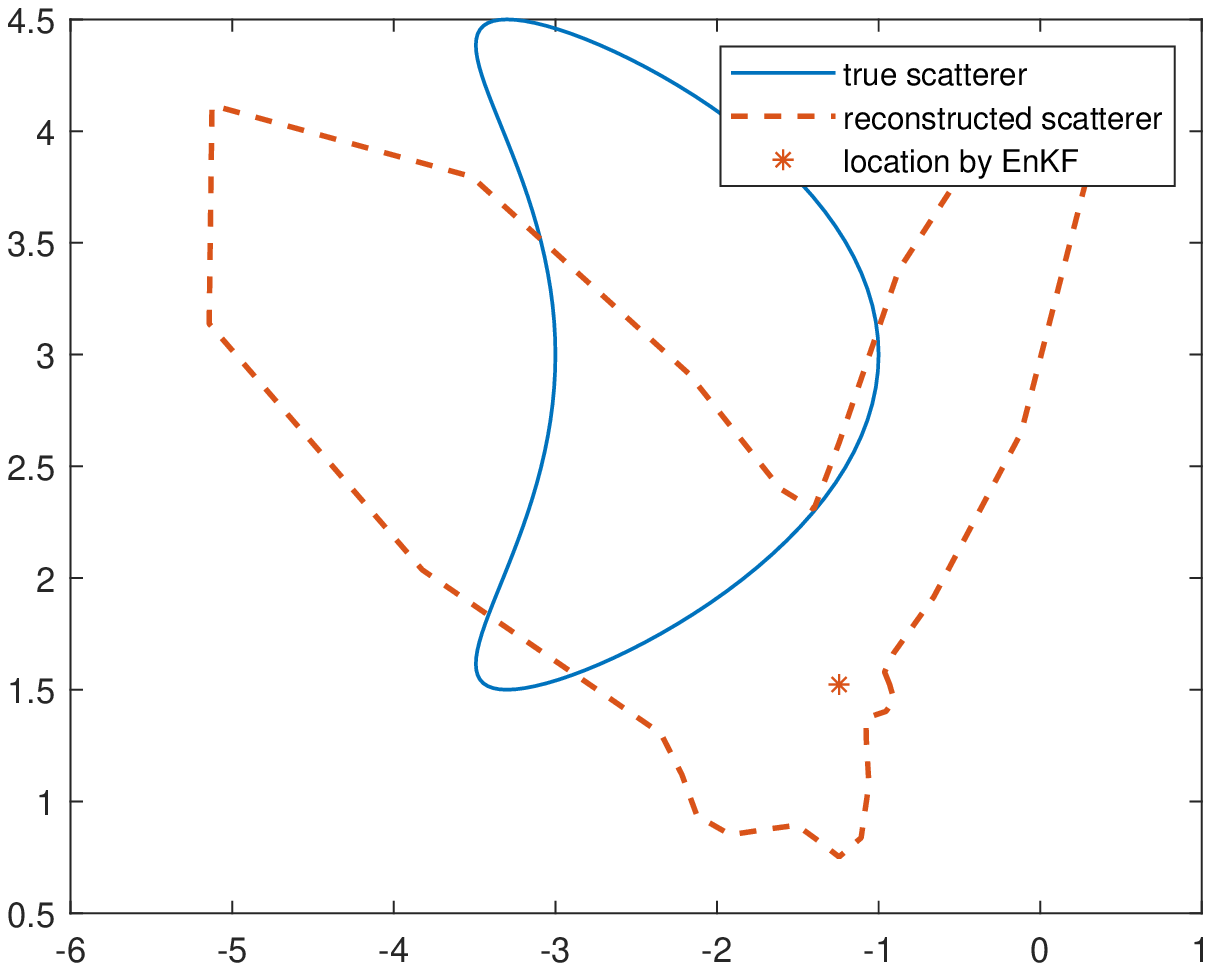}
\end{minipage}
\begin{minipage}{0.45\linewidth}
\centering
\includegraphics[height=50mm]{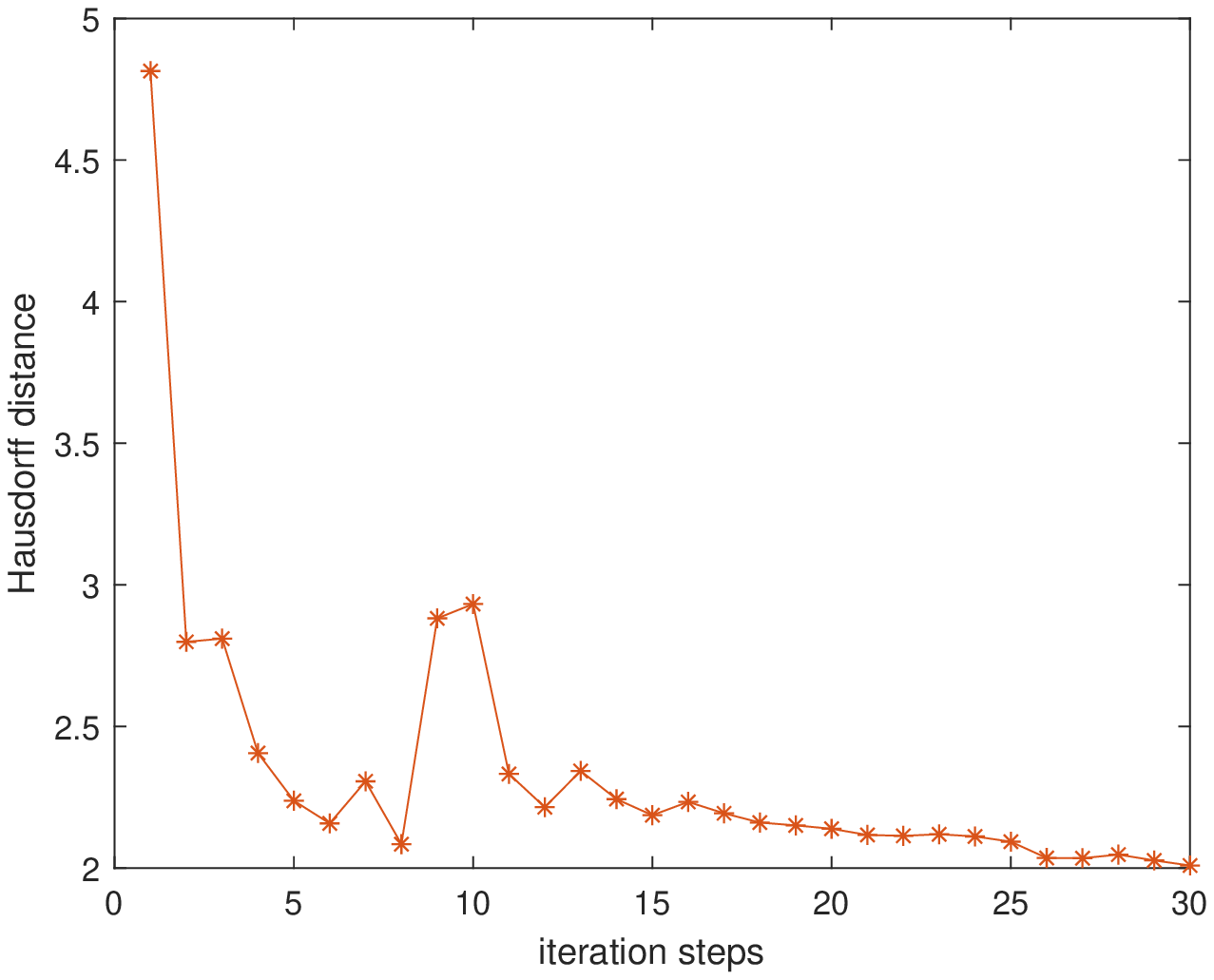}
\end{minipage}
\caption{Poor initial location for the EnKF. Left: reconstructed boundary. Right: Hausdorff distance.}\label{poorguess_kite}
\end{figure}

\begin{figure}[!h]
\centering
\begin{minipage}{0.32\linewidth}
\centering
\includegraphics[height=40mm]{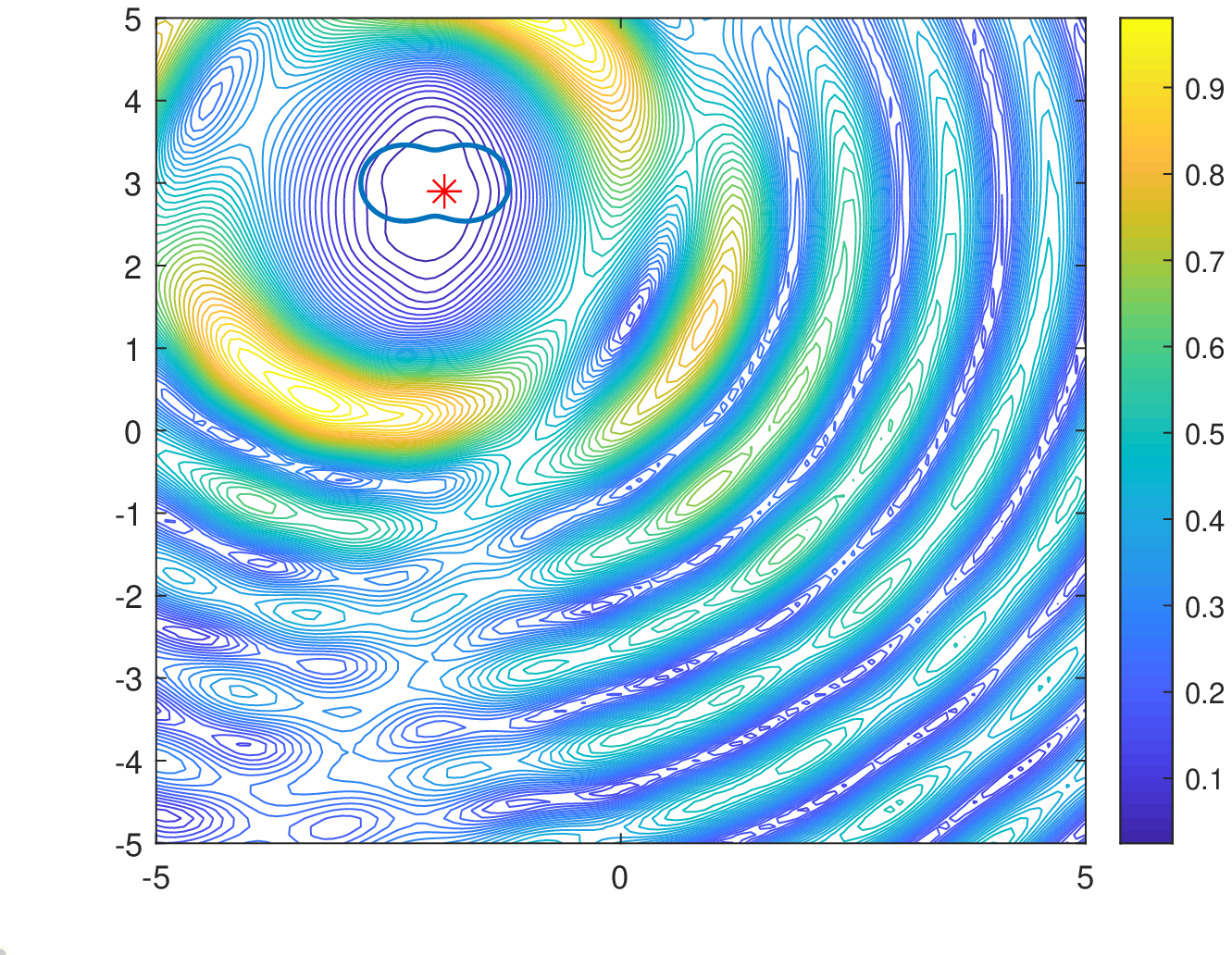}
\end{minipage}
\hfill
\begin{minipage}{0.32\linewidth}
\centering
\includegraphics[height=40mm]{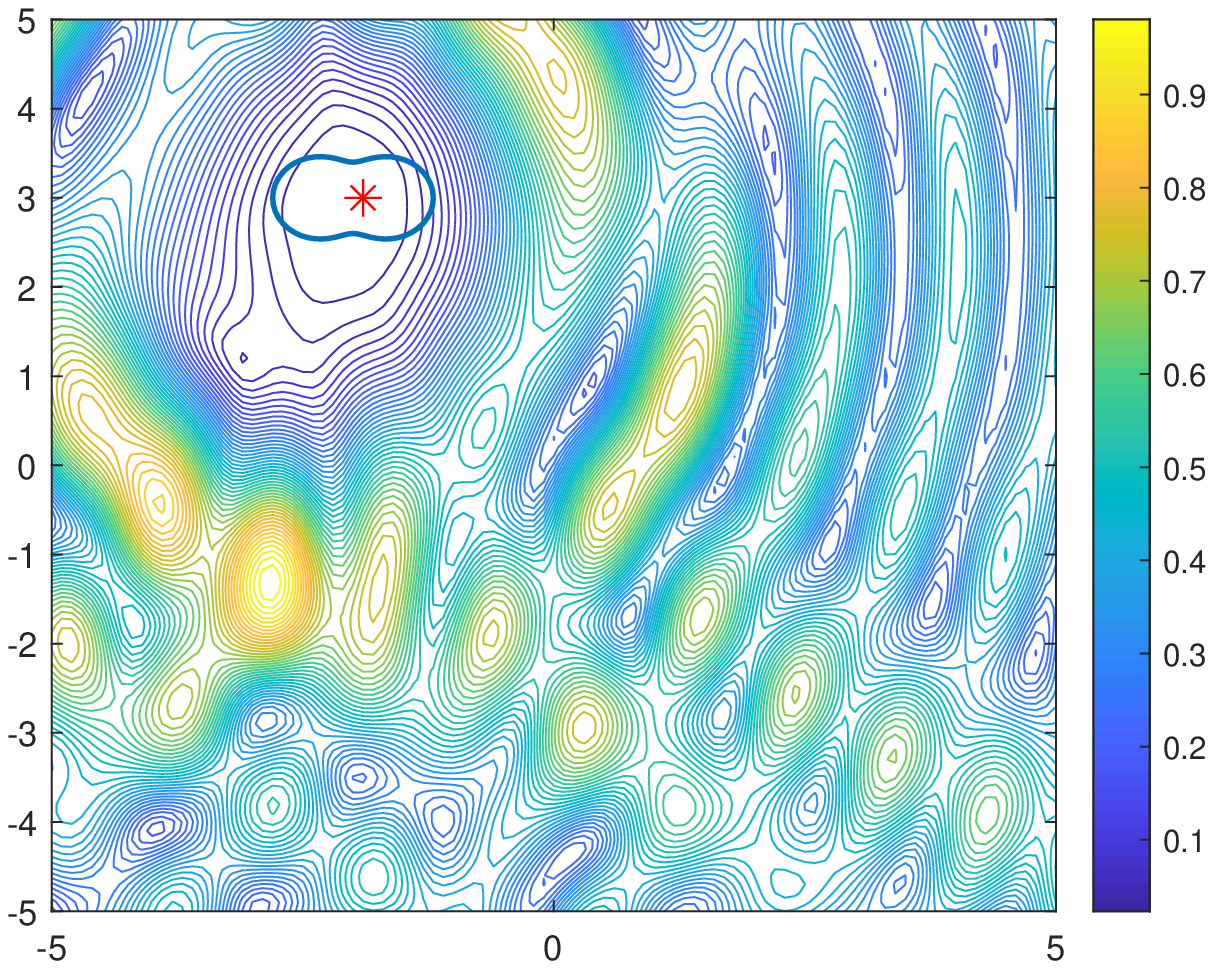}
\end{minipage}
\hfill
\begin{minipage}{0.32\linewidth}
\centering
\includegraphics[height=40mm]{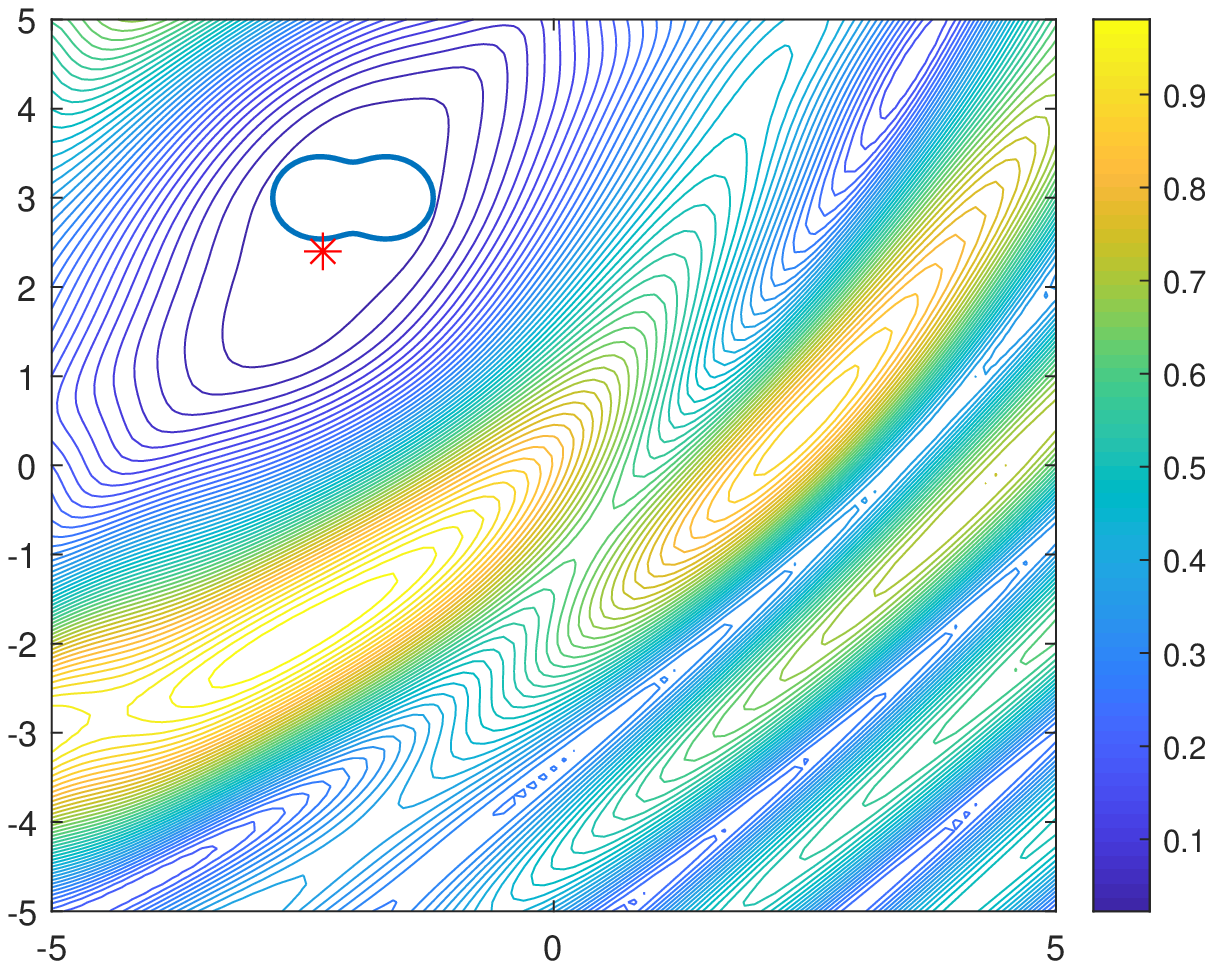}
\end{minipage}

\begin{minipage}{0.32\linewidth}
\centering
\includegraphics[height=40mm]{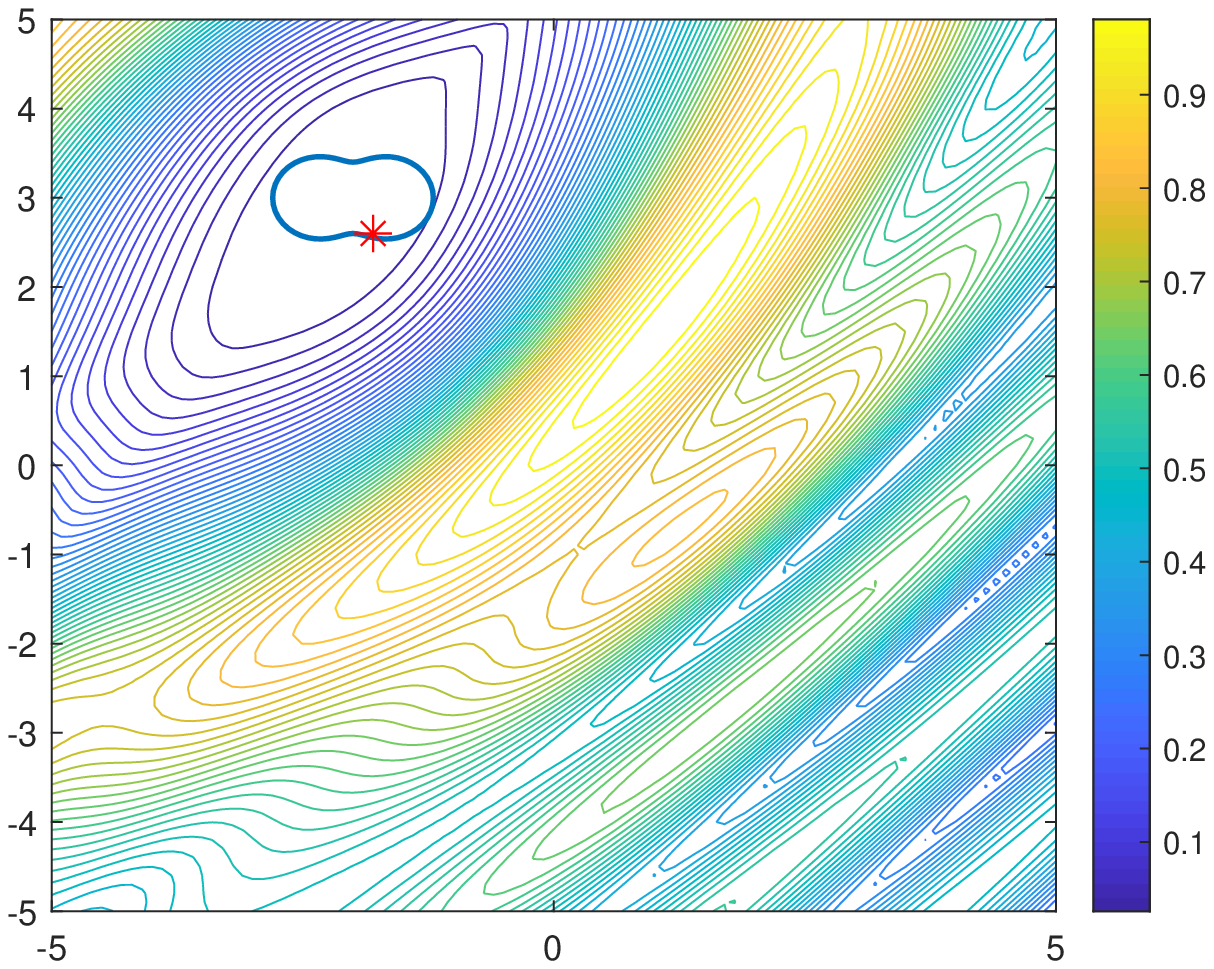}
\end{minipage}
\hfill
\begin{minipage}{0.32\linewidth}
\centering
\includegraphics[height=40mm]{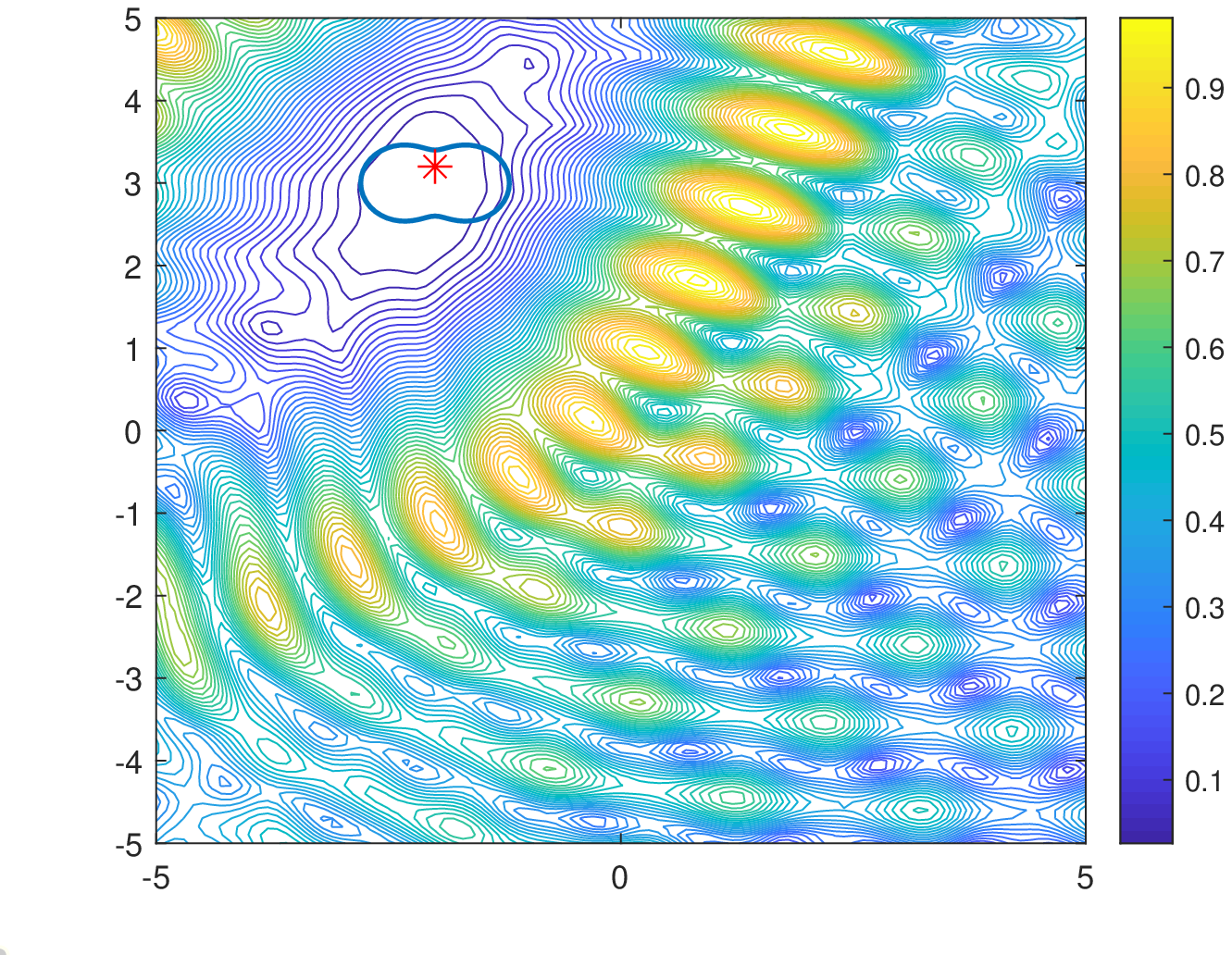}
\end{minipage}
\hfill
\begin{minipage}{0.32\linewidth}
\centering
\includegraphics[height=40mm]{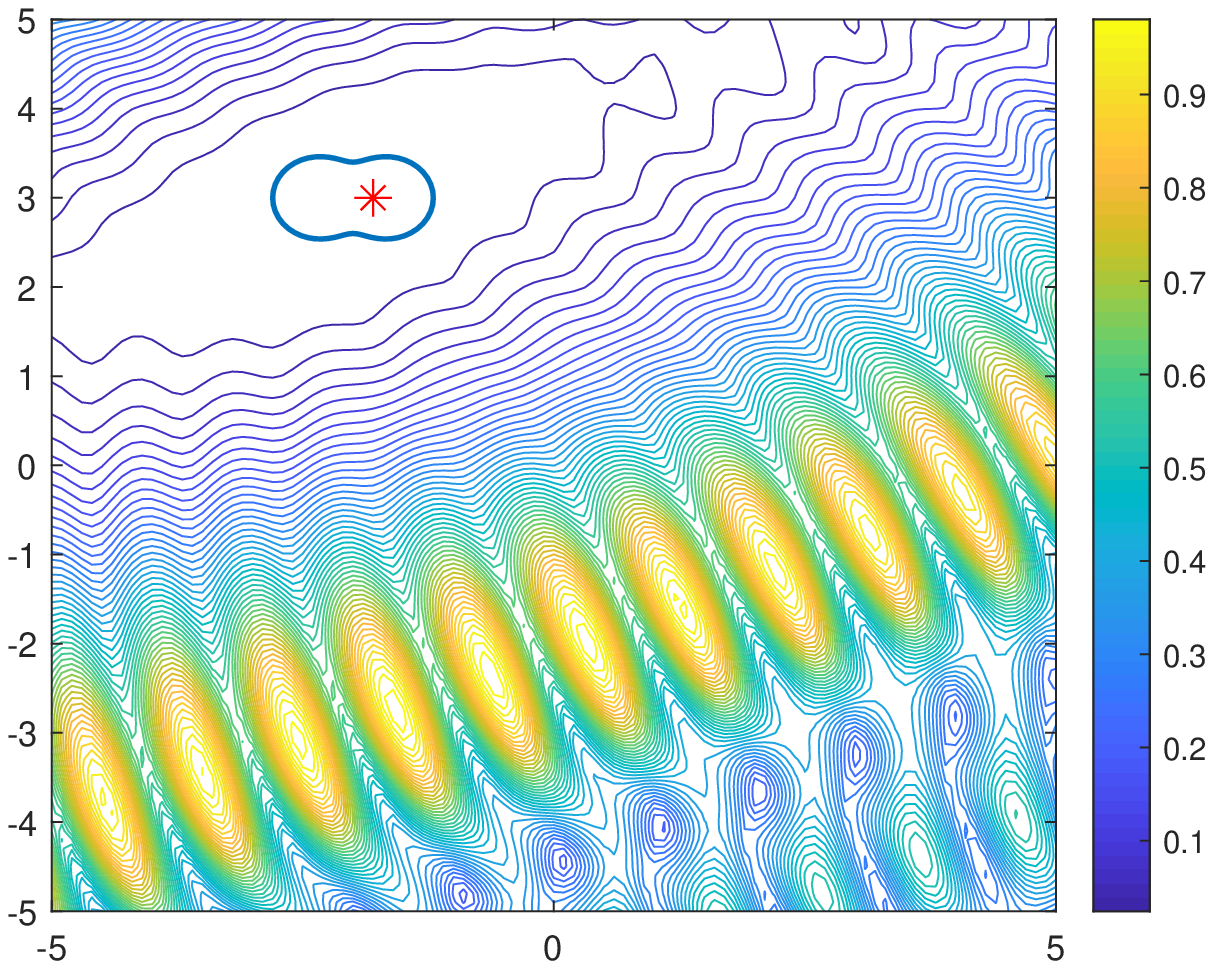}
\end{minipage}
\caption{Contour plots of the indicator function $I_{ESM}(\bm{z})$. Top row: $\gamma_1^i$, from left to right: $\gamma_1^o$, $\gamma_2^o$, $\gamma_3^o$. Bottom row: $\gamma_2^i$, from left to right: $\gamma_3^o$, $\gamma_4^o$, $\gamma_5^o$.}\label{ESM_peanut}
\end{figure}

\begin{figure}[!h]
\centering
\begin{minipage}{0.32\linewidth}
\centering
\includegraphics[height=40mm]{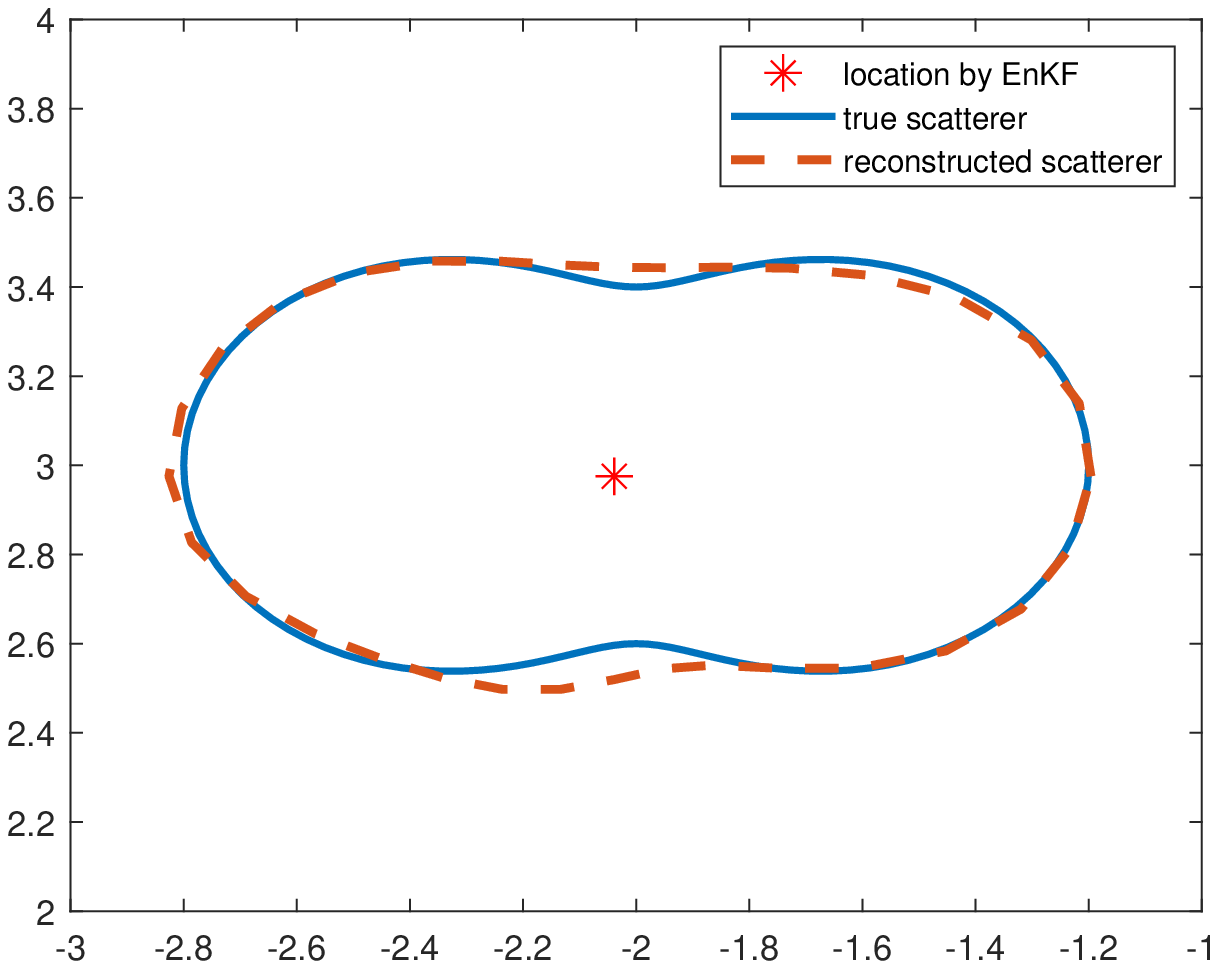}
\end{minipage}
\hfill
\begin{minipage}{0.32\linewidth}
\centering
\includegraphics[height=40mm]{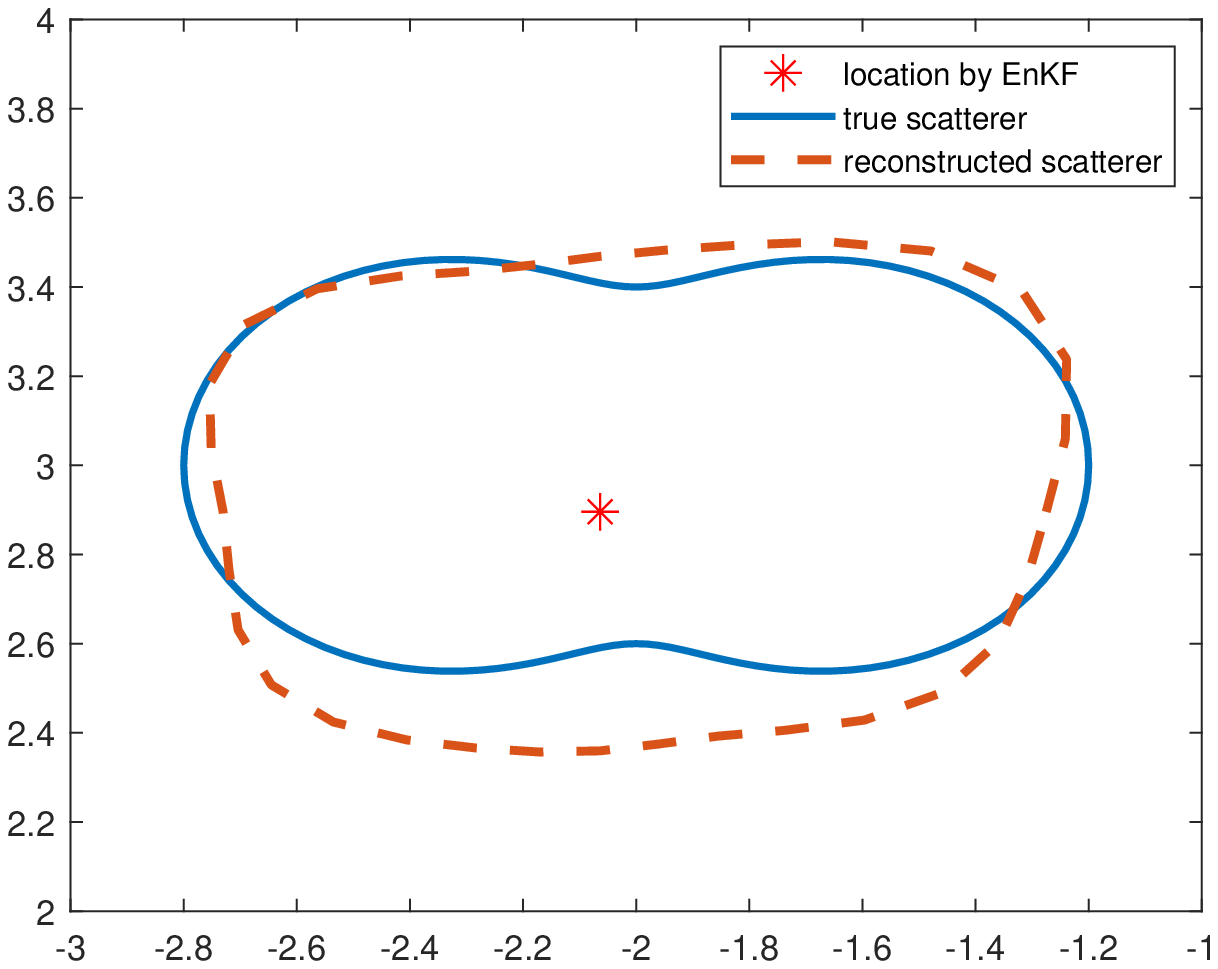}
\end{minipage}
\hfill
\begin{minipage}{0.32\linewidth}
\centering
\includegraphics[height=40mm]{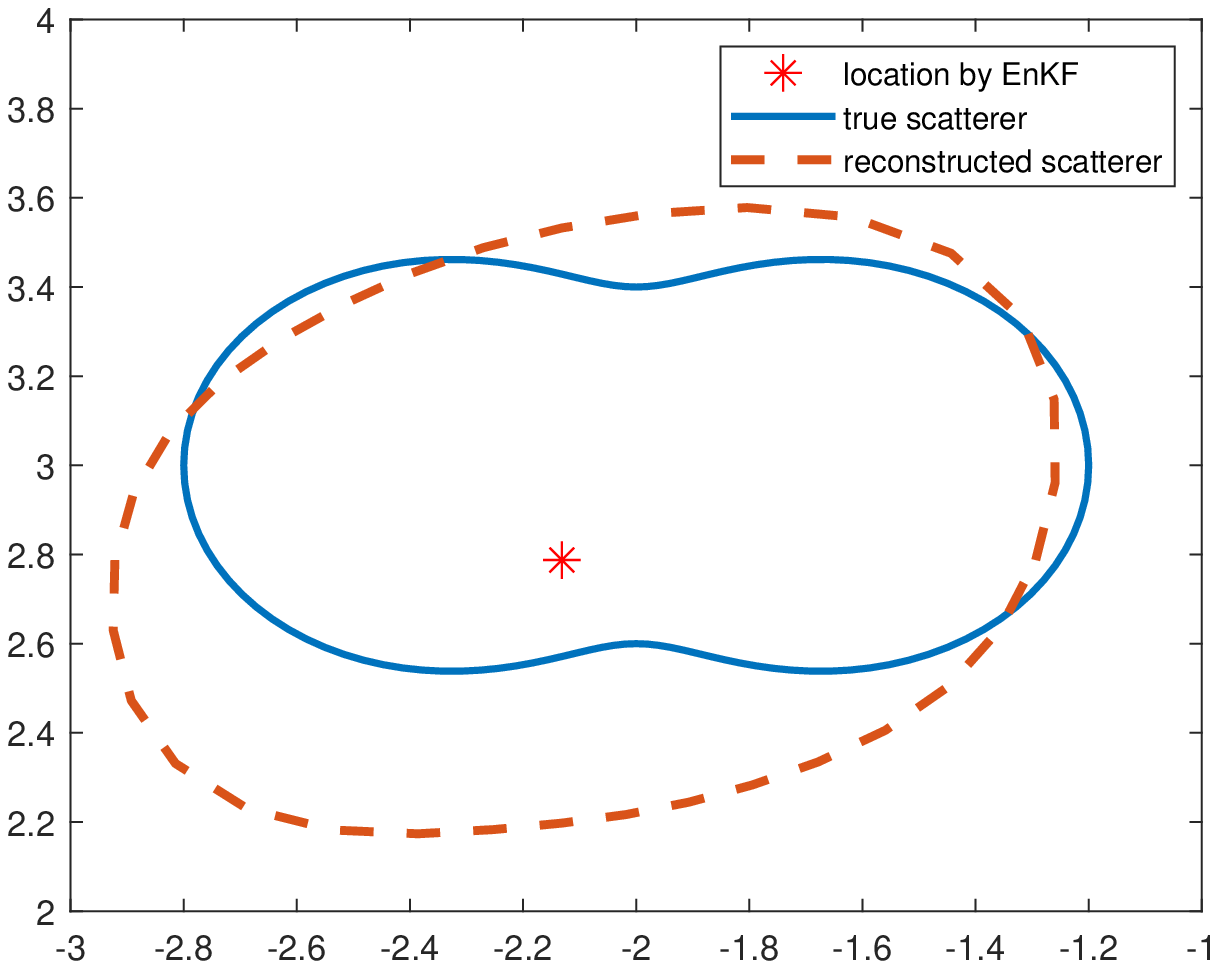}
\end{minipage}

\begin{minipage}{0.32\linewidth}
\centering
\includegraphics[height=40mm]{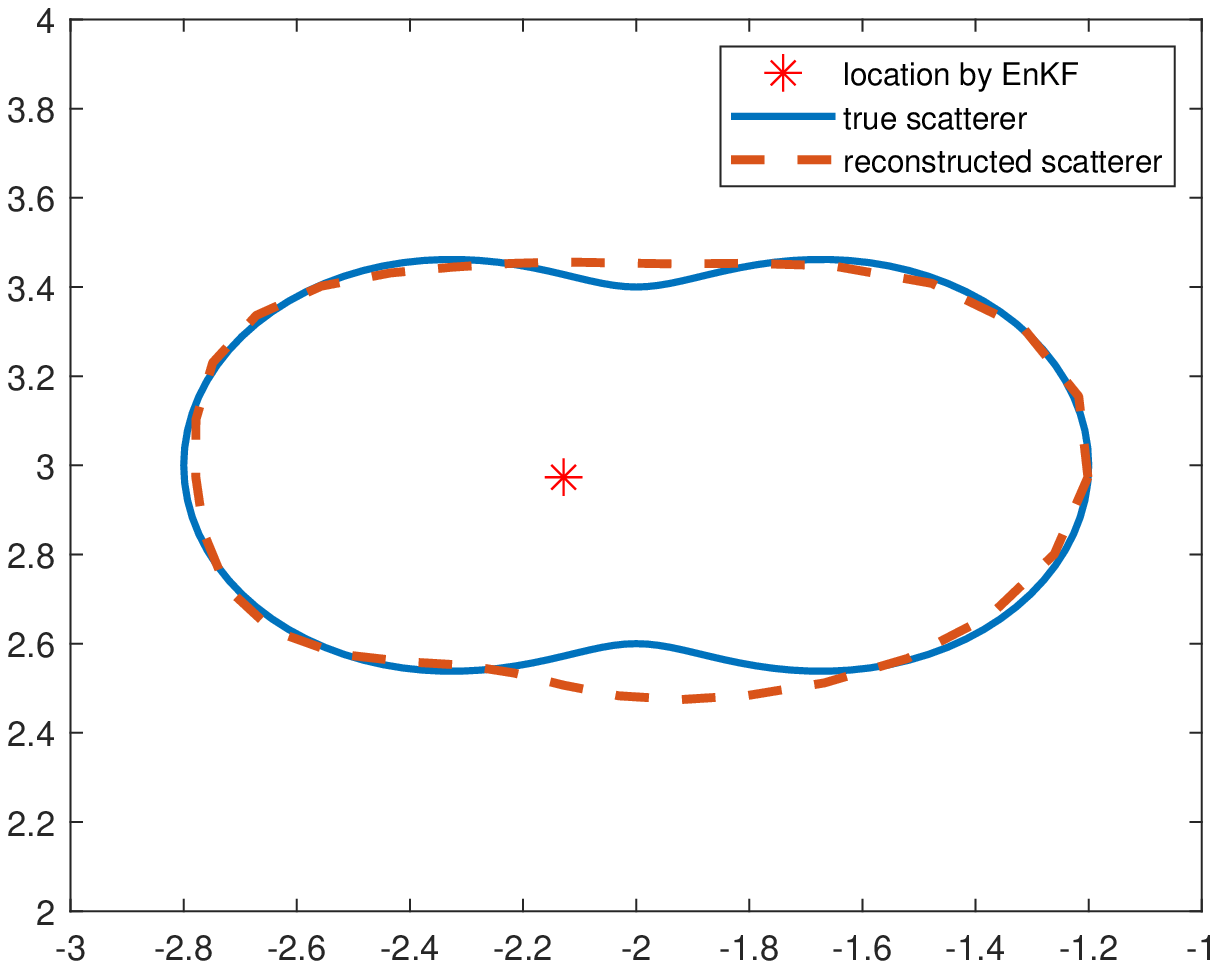}
\end{minipage}
\hfill
\begin{minipage}{0.32\linewidth}
\centering
\includegraphics[height=40mm]{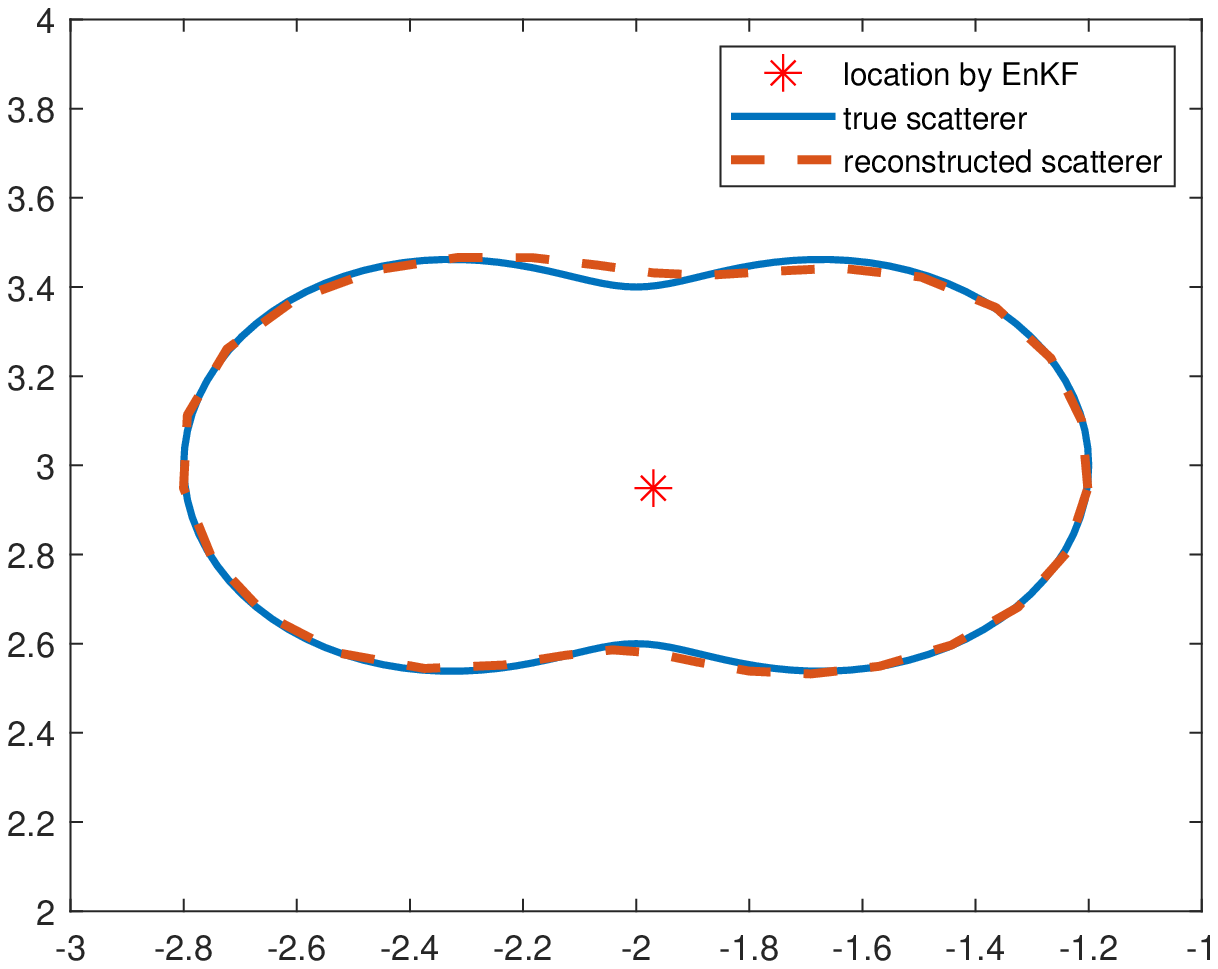}
\end{minipage}
\hfill
\begin{minipage}{0.32\linewidth}
\centering
\includegraphics[height=40mm]{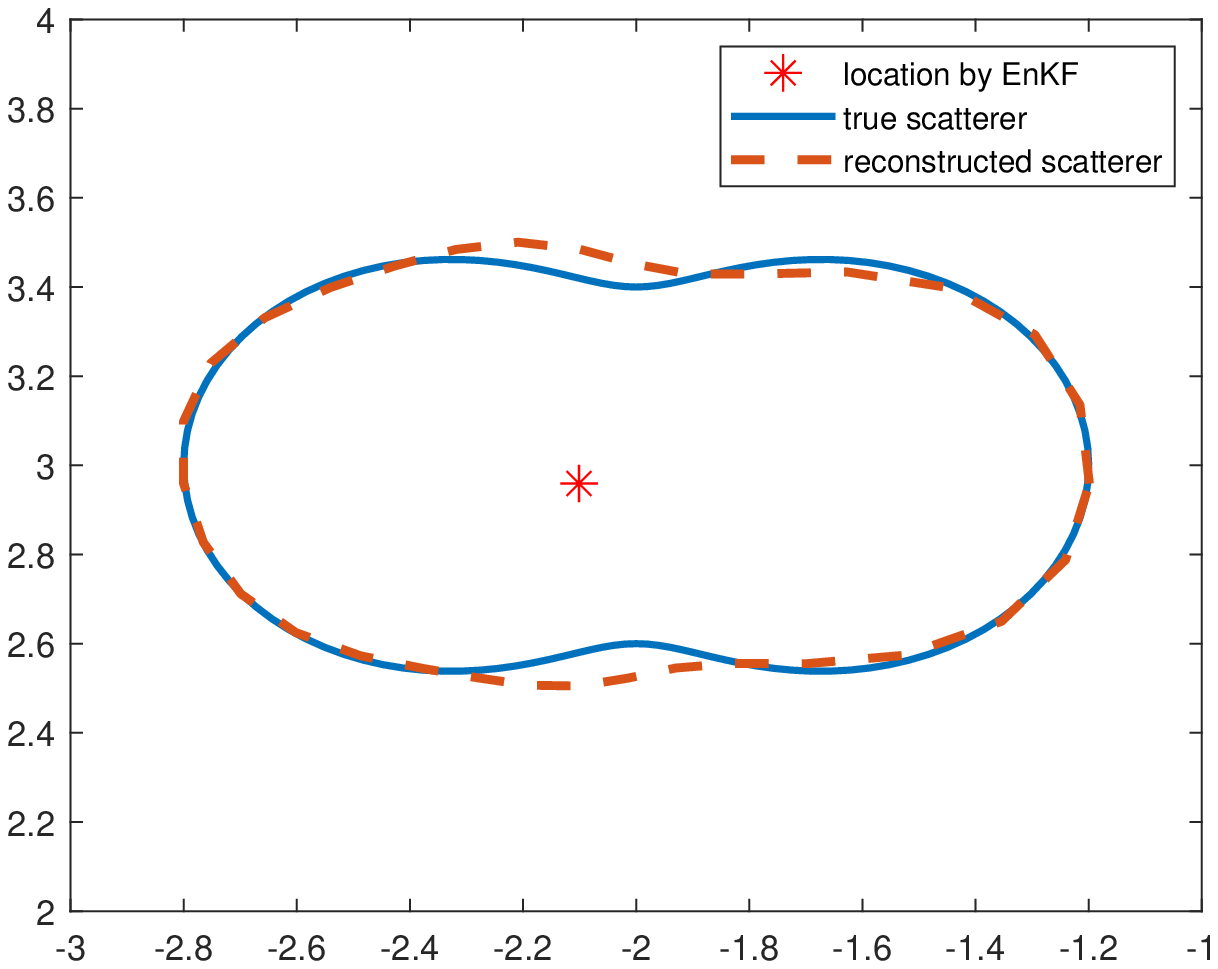}
\end{minipage}
\caption{Boundary reconstructions by EnKF algorithm. Top row: $\gamma_1^i$, from left to right: $\gamma_1^o$, $\gamma_2^o$, $\gamma_3^o$. Bottom row: $\gamma_2^i$, from left to right: $\gamma_3^o$, $\gamma_4^o$, $\gamma_5^o$.}\label{EnKF_peanut}
\end{figure}

\begin{figure}[!h]
\centering
\begin{minipage}{0.45\linewidth}
\centering
\includegraphics[height=50mm]{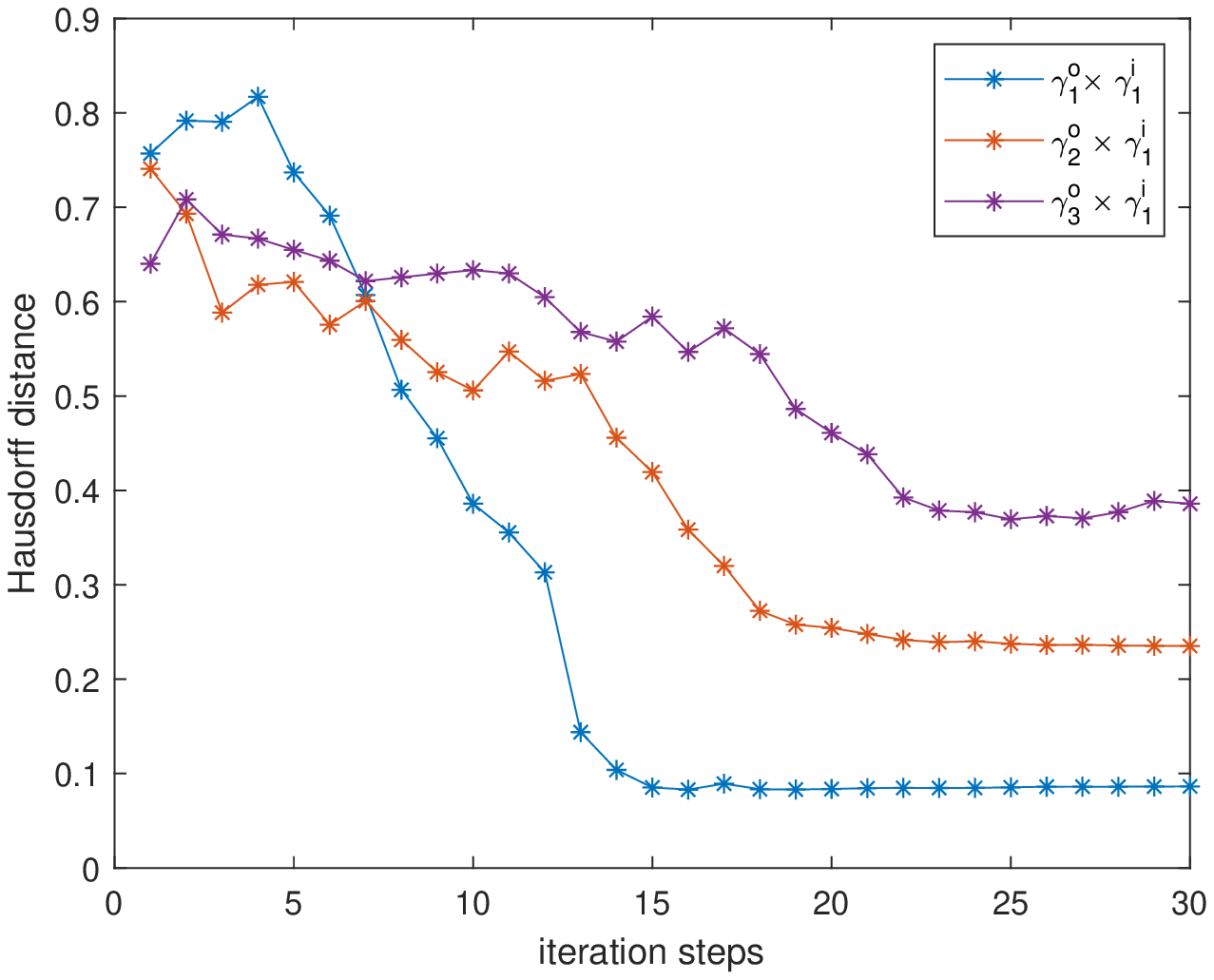}
\end{minipage}
\begin{minipage}{0.45\linewidth}
\centering
\includegraphics[height=50mm]{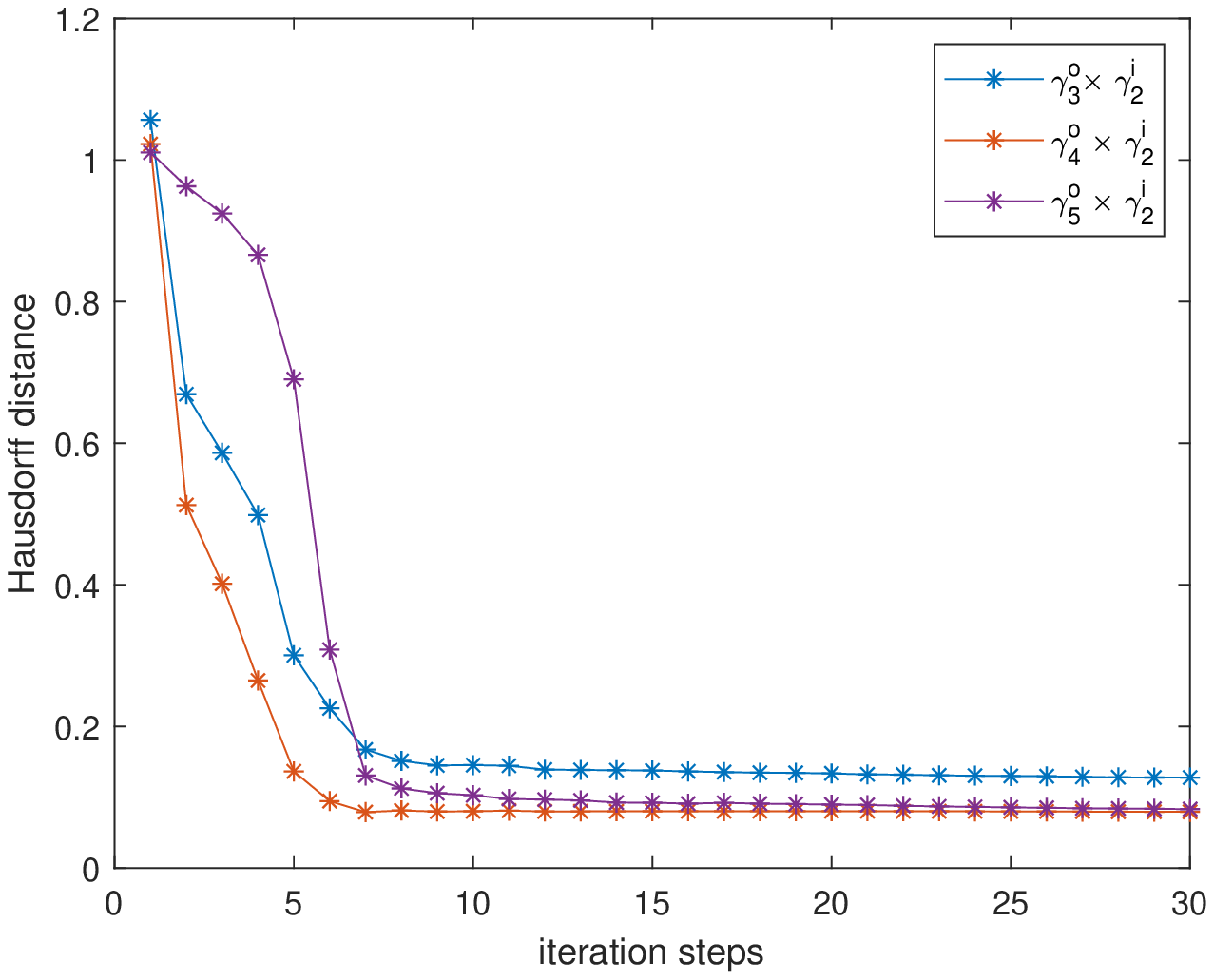}
\end{minipage}
\caption{Hausdorff distance $d_H(\partial \Omega_{\text{exact}},\partial \Omega_{\text{inv}})$ with respect to the iteration steps. Left: $\gamma^i=\gamma^i_1$. Right: $\gamma^i=\gamma_2^i$.}\label{dH_peanut}
\end{figure}

\begin{figure}[!h]
\centering
\begin{minipage}{0.32\linewidth}
\centering
\includegraphics[height=40mm]{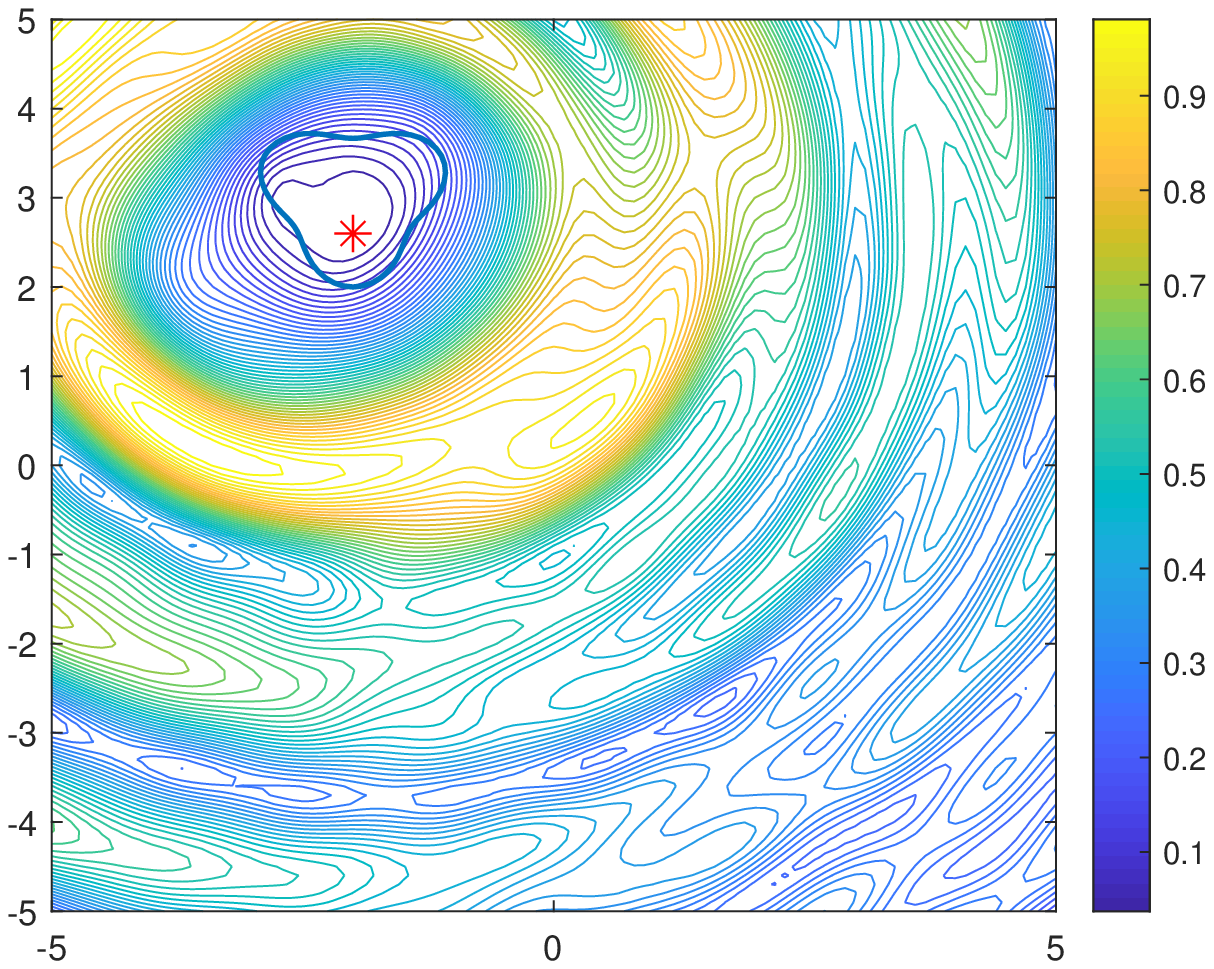}
\end{minipage}
\hfill
\begin{minipage}{0.32\linewidth}
\centering
\includegraphics[height=40mm]{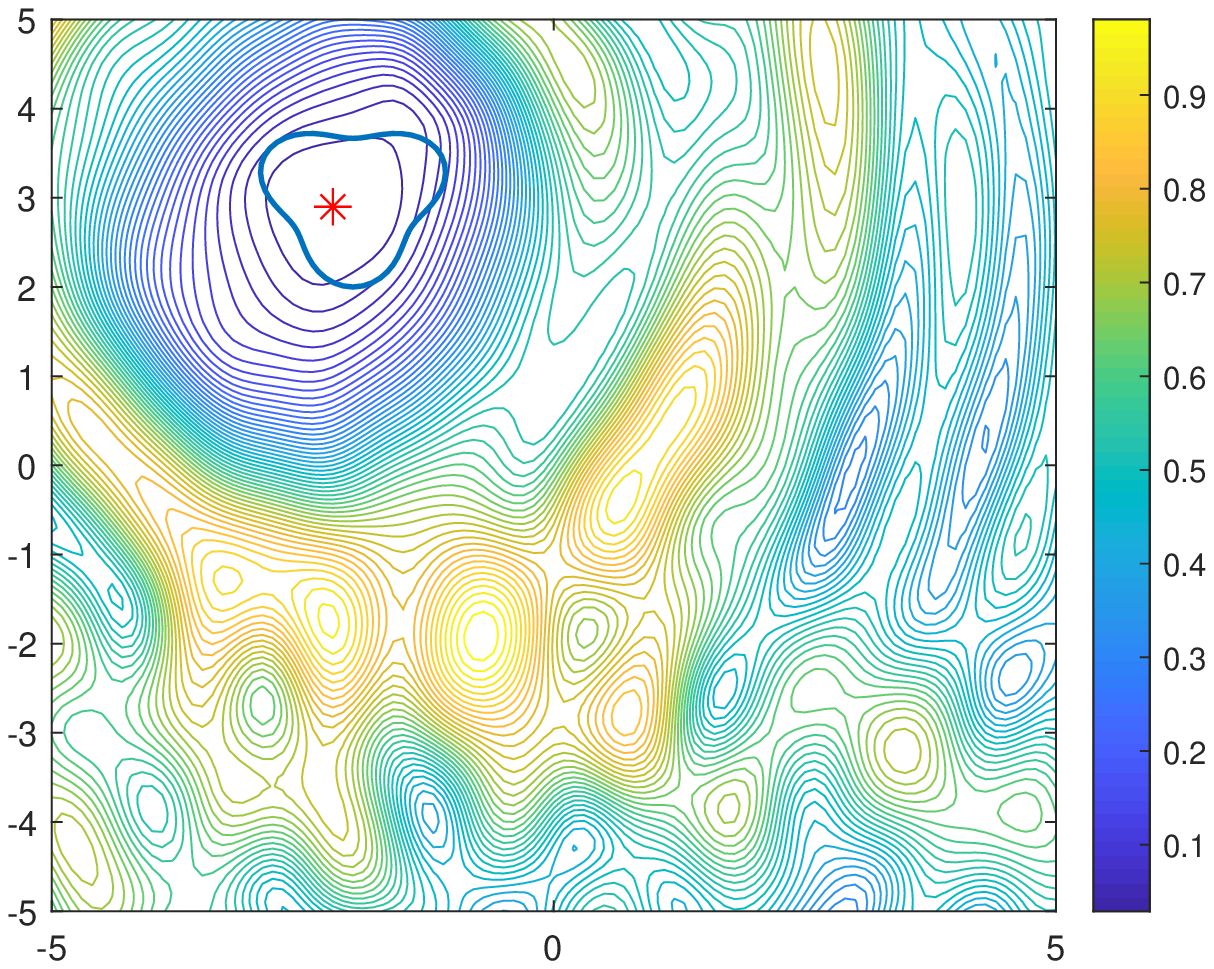}
\end{minipage}
\hfill
\begin{minipage}{0.32\linewidth}
\centering
\includegraphics[height=40mm]{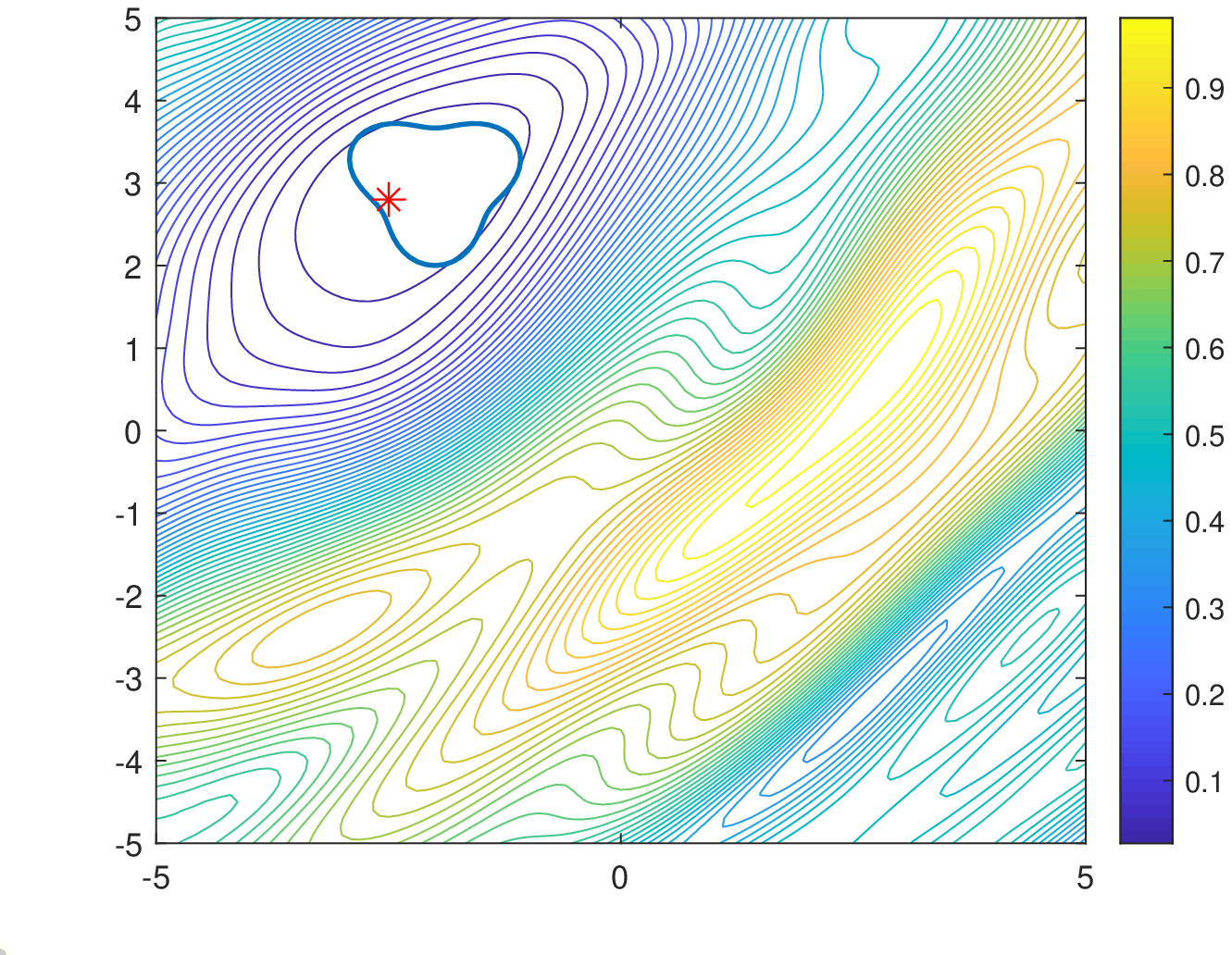}
\end{minipage}

\begin{minipage}{0.32\linewidth}
\centering
\includegraphics[height=40mm]{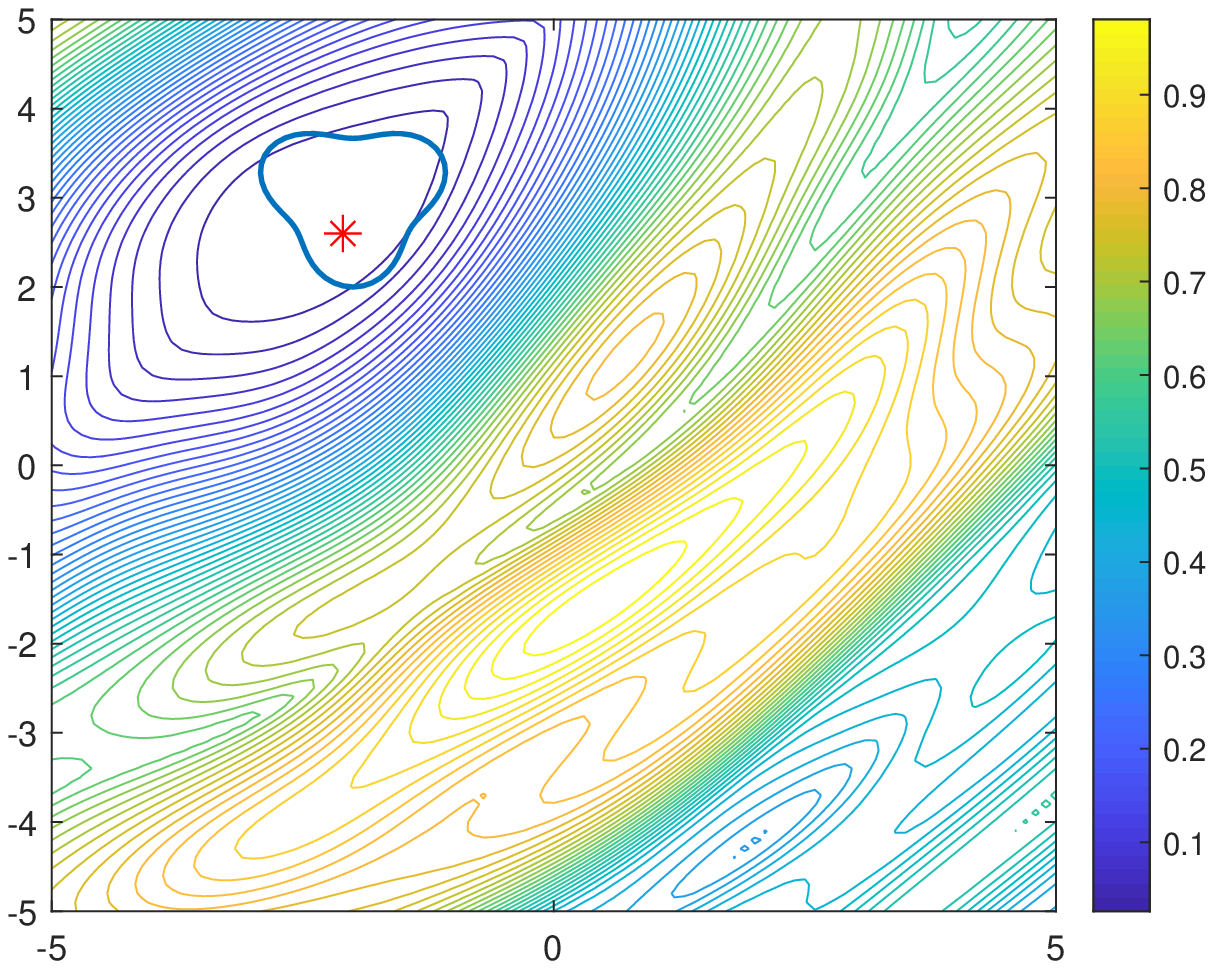}
\end{minipage}
\hfill
\begin{minipage}{0.32\linewidth}
\centering
\includegraphics[height=40mm]{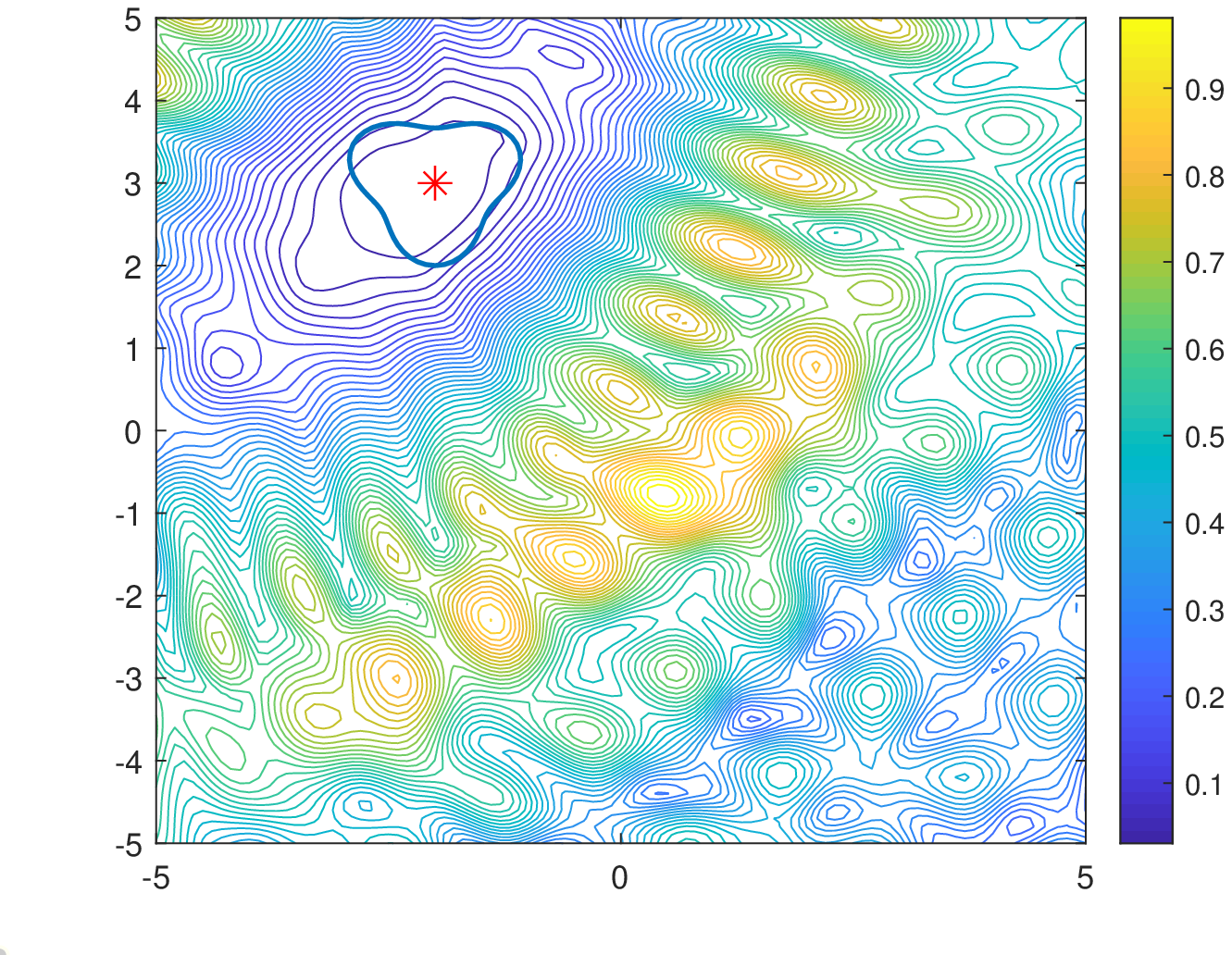}
\end{minipage}
\hfill
\begin{minipage}{0.32\linewidth}
\centering
\includegraphics[height=40mm]{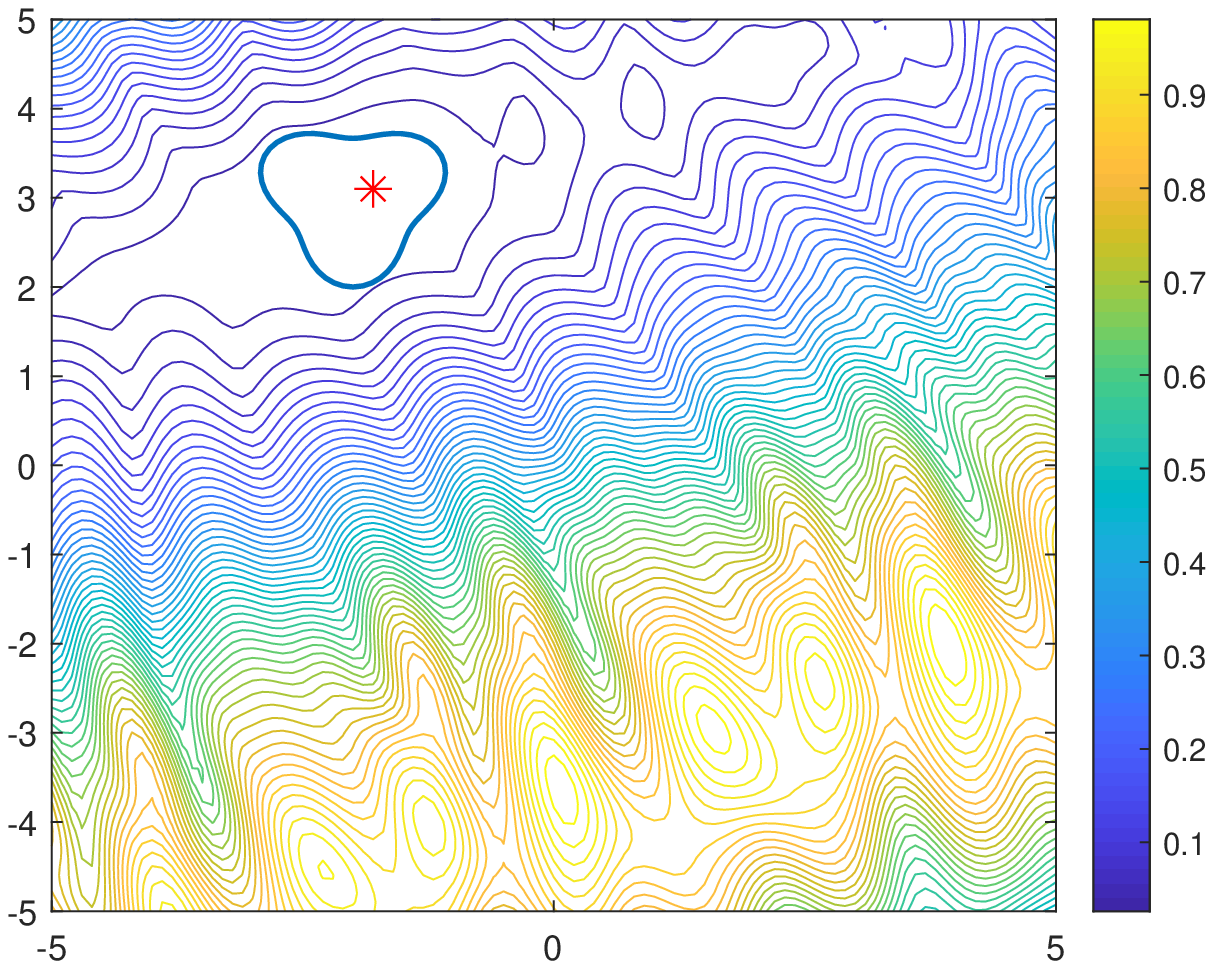}
\end{minipage}
\caption{Contour plots of the indicator function $I_{ESM}(\bm{z})$. Top row: $\gamma_1^i$, from left to right: $\gamma_1^o$, $\gamma_2^o$, $\gamma_3^o$. Bottom row: $\gamma_2^i$, from left to right: $\gamma_3^o$, $\gamma_4^o$, $\gamma_5^o$.}\label{ESM_pear}
\end{figure}

\begin{figure}[!h]
\centering
\begin{minipage}{0.32\linewidth}
\centering
\includegraphics[height=40mm]{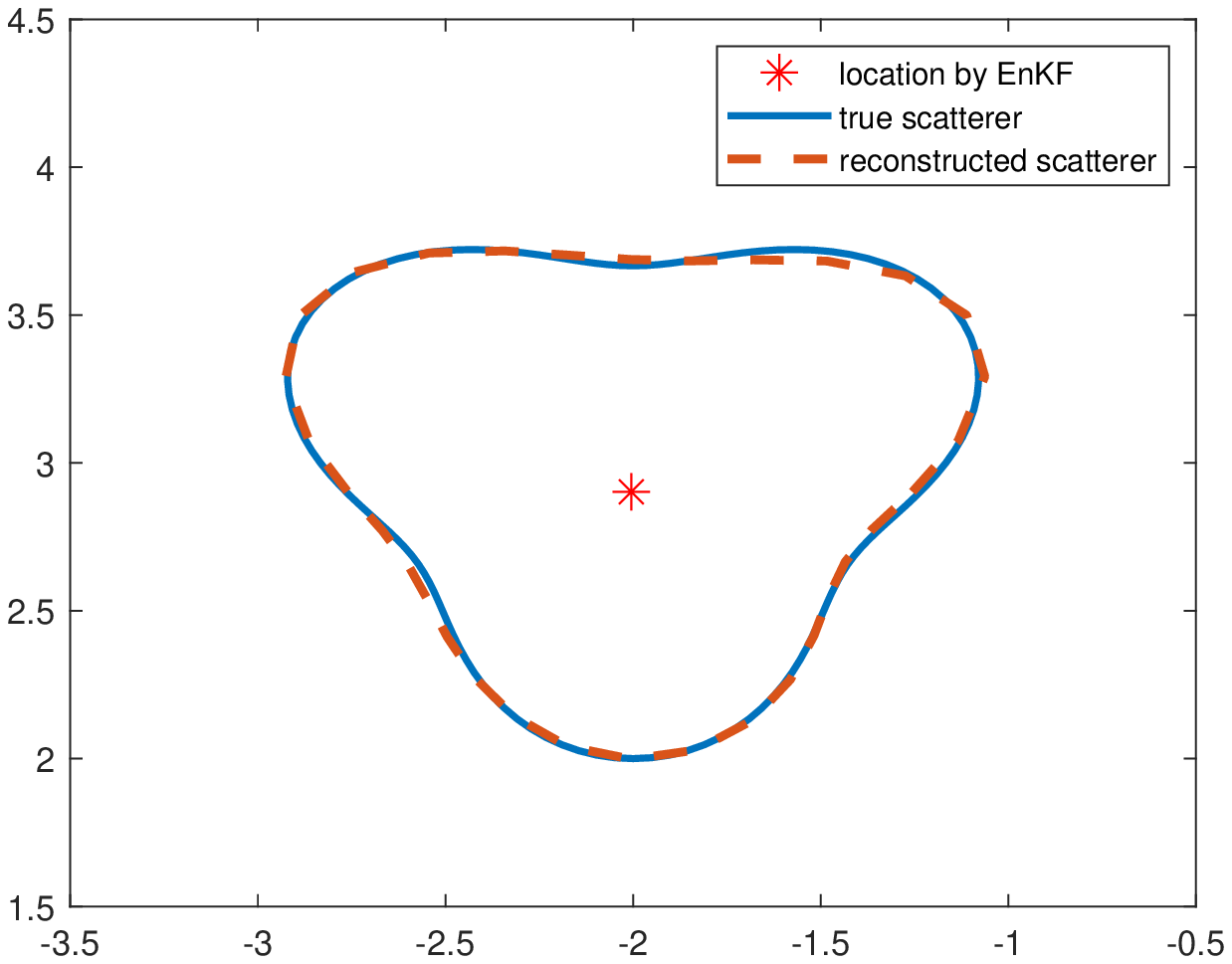}
\end{minipage}
\hfill
\begin{minipage}{0.32\linewidth}
\centering
\includegraphics[height=40mm]{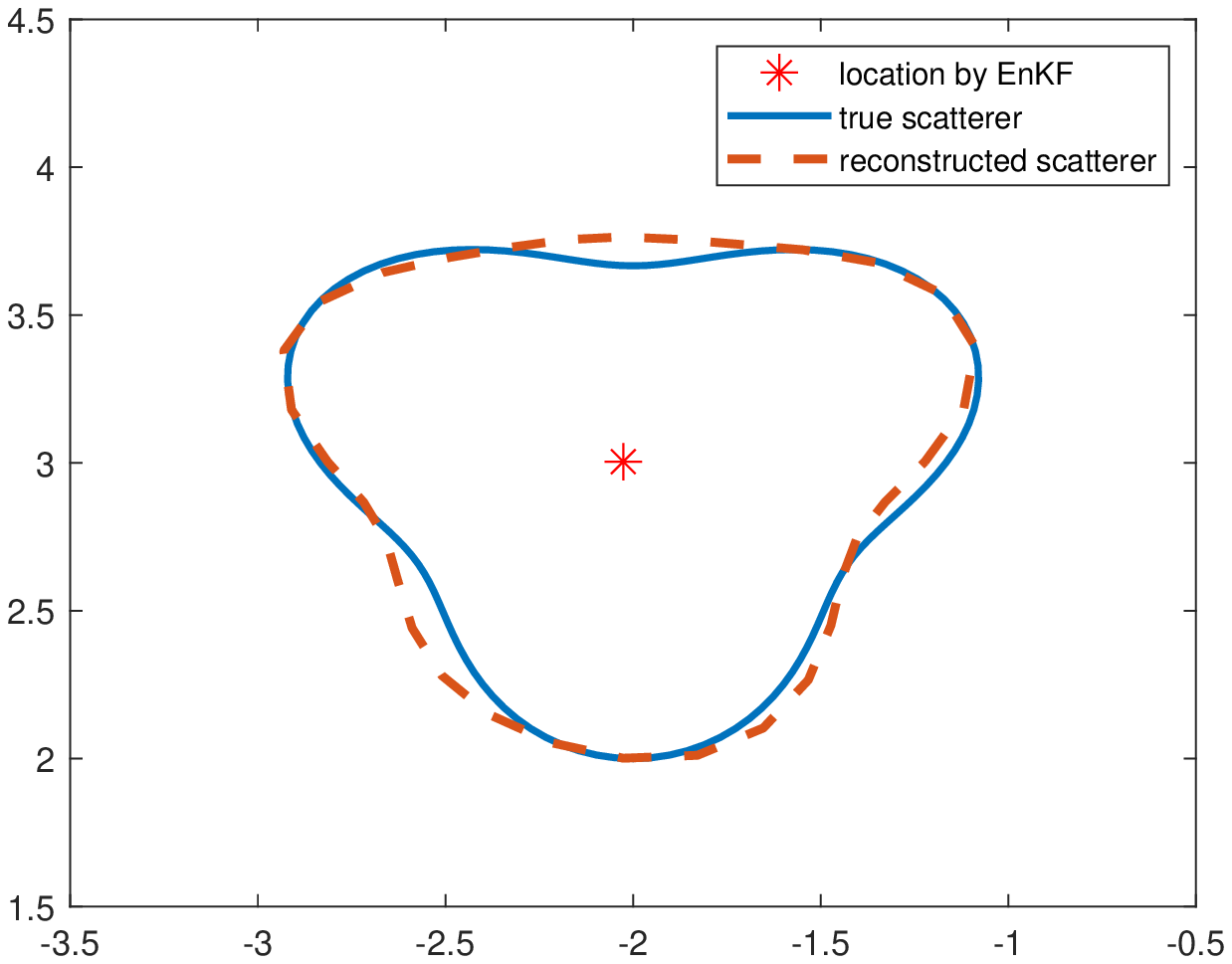}
\end{minipage}
\hfill
\begin{minipage}{0.32\linewidth}
\centering
\includegraphics[height=40mm]{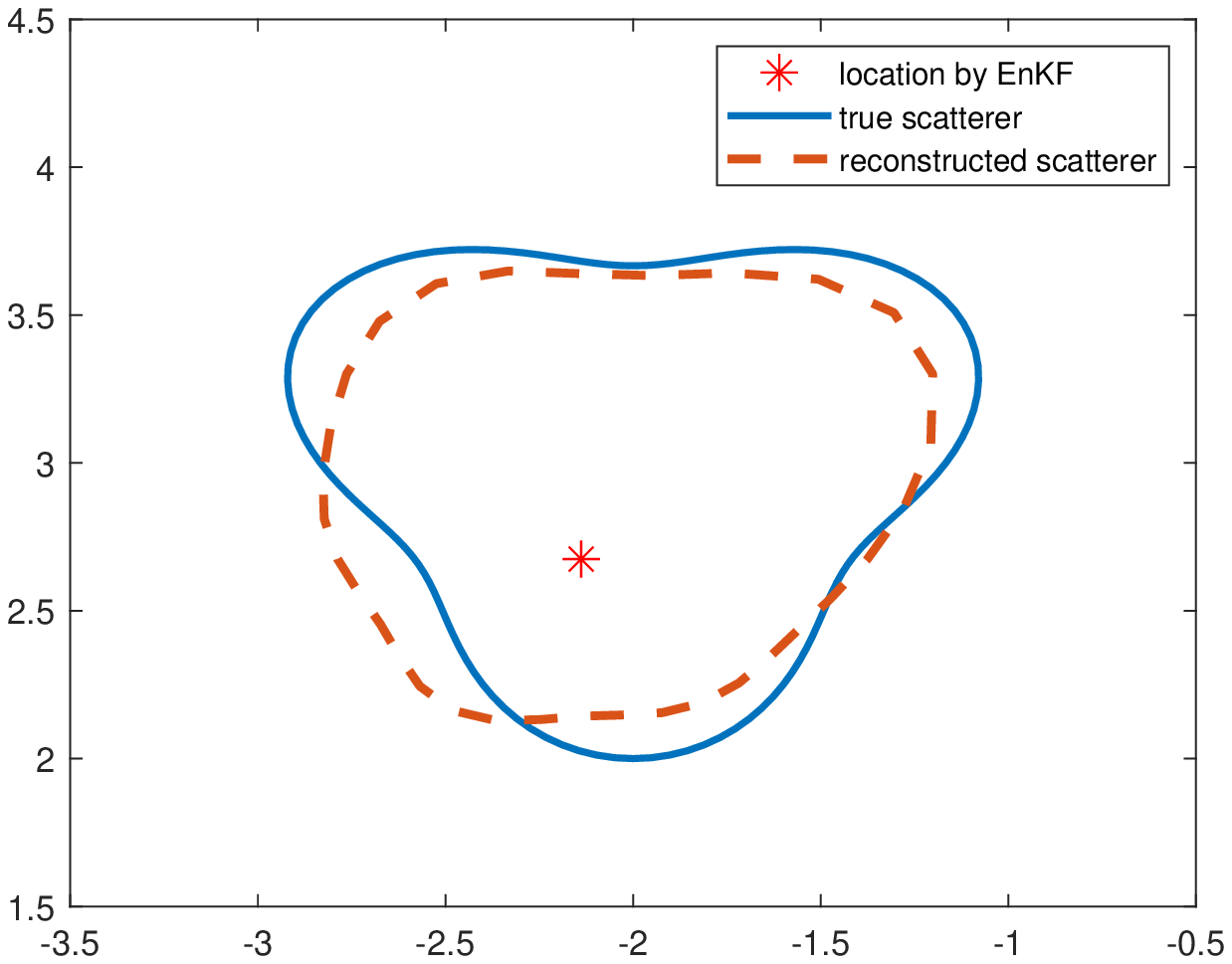}
\end{minipage}

\begin{minipage}{0.32\linewidth}
\centering
\includegraphics[height=40mm]{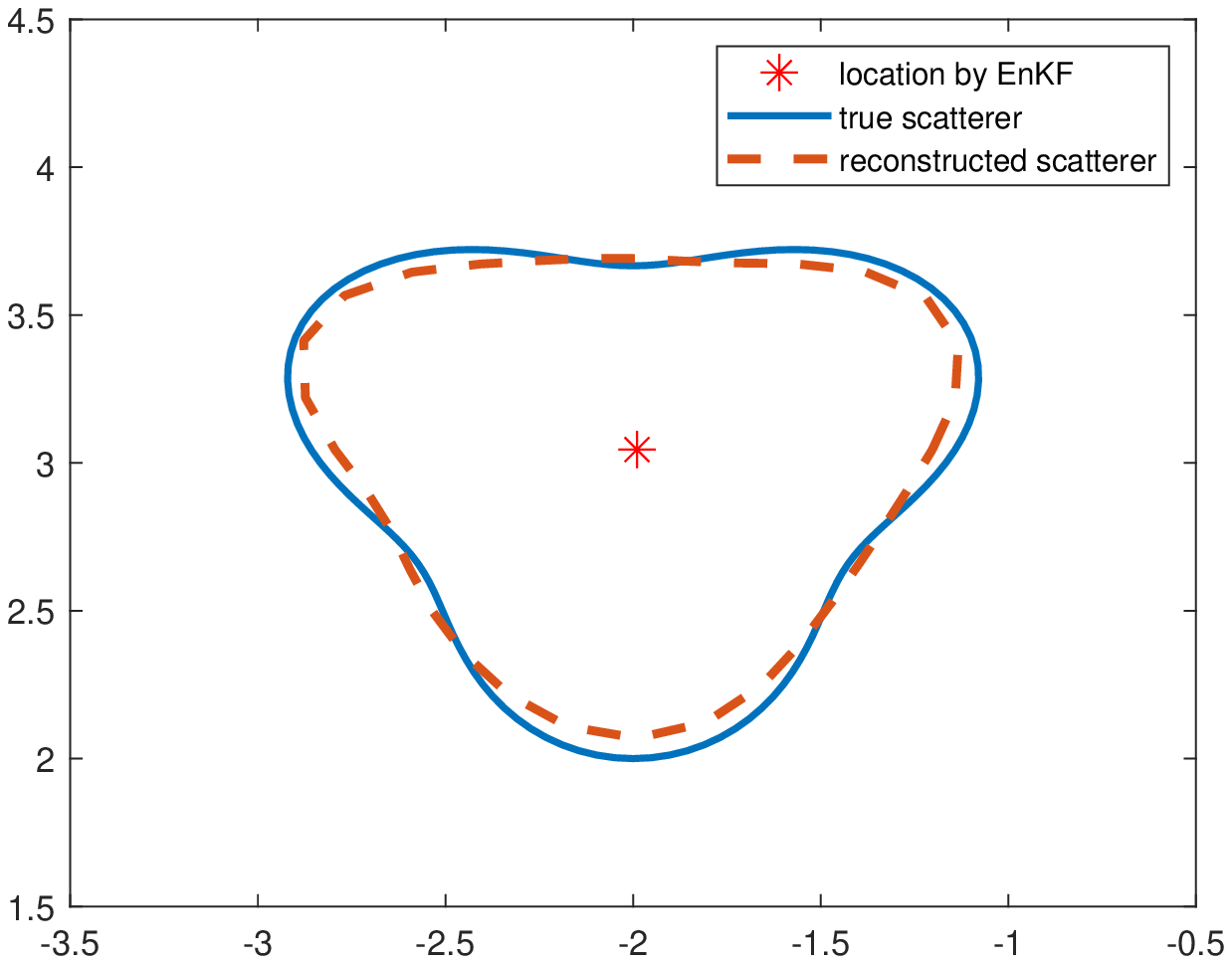}
\end{minipage}
\hfill
\begin{minipage}{0.32\linewidth}
\centering
\includegraphics[height=40mm]{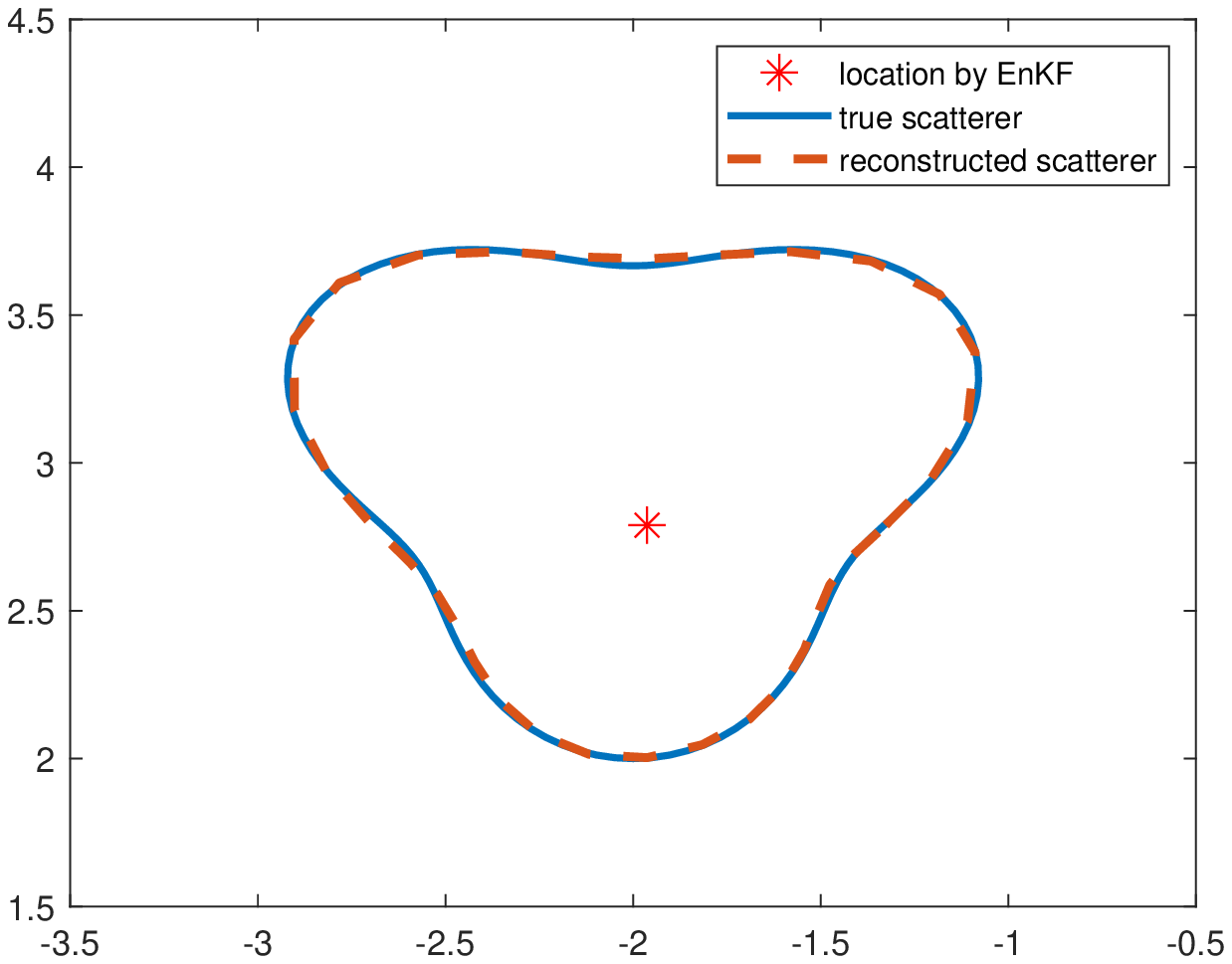}
\end{minipage}
\hfill
\begin{minipage}{0.32\linewidth}
\centering
\includegraphics[height=40mm]{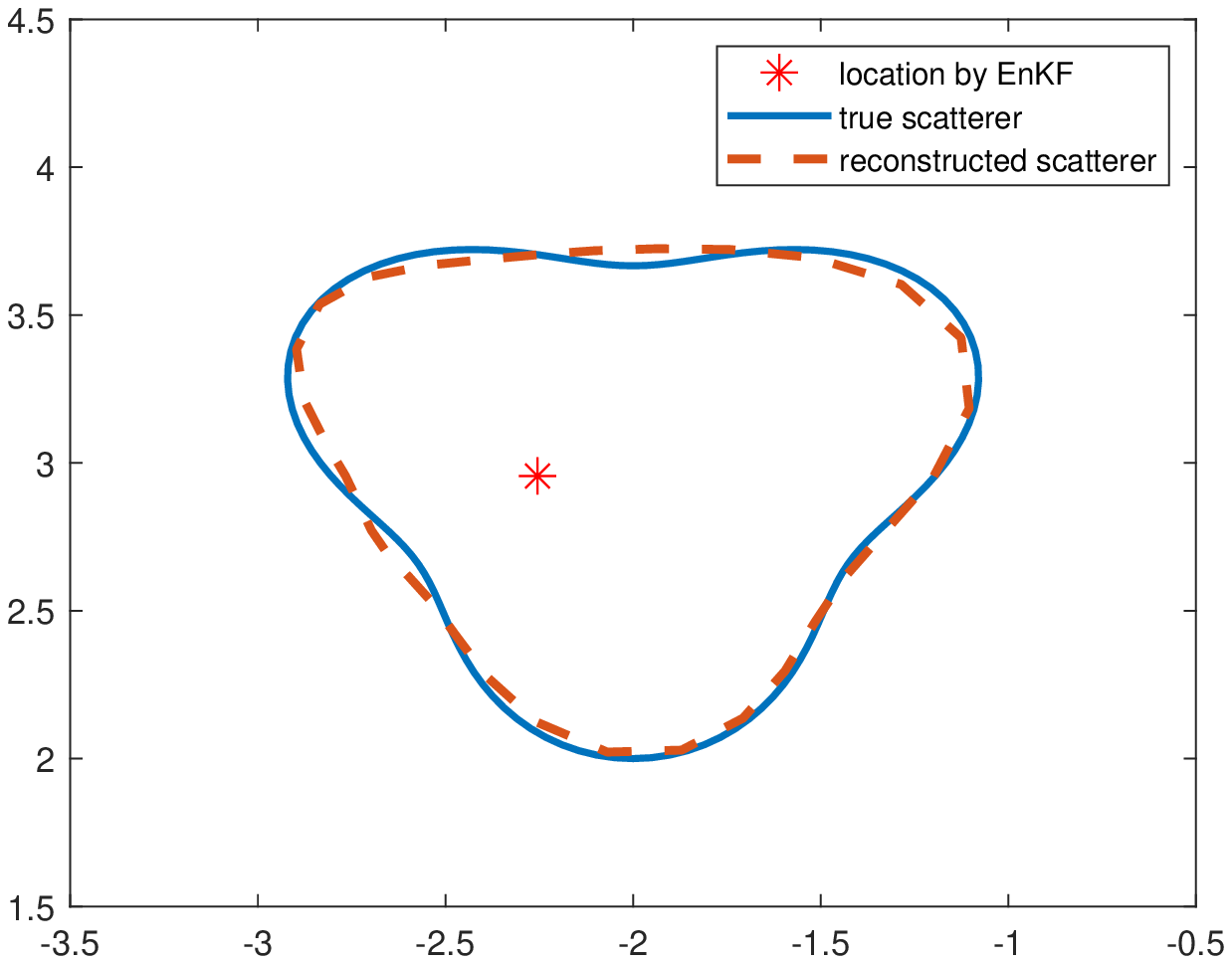}
\end{minipage}
\caption{Boundary reconstructions by EnKF algorithm. Top row: $\gamma_1^i$, from left to right: $\gamma_1^o$, $\gamma_2^o$, $\gamma_3^o$. Bottom row: $\gamma_2^i$, from left to right: $\gamma_3^o$, $\gamma_4^o$, $\gamma_5^o$.}\label{EnKF_pear}
\end{figure}

\begin{figure}[!h]
\centering
\begin{minipage}{0.45\linewidth}
\centering
\includegraphics[height=50mm]{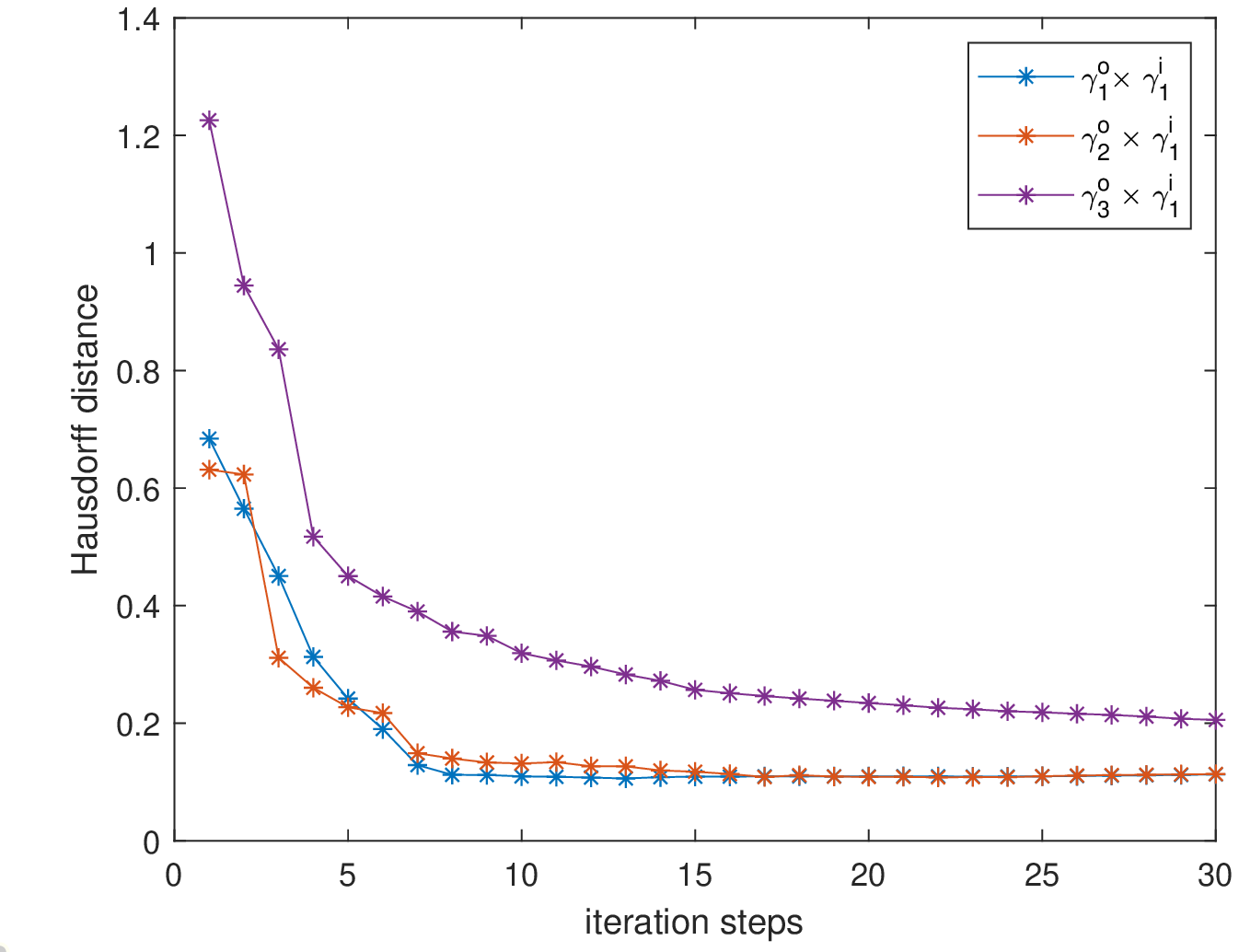}
\end{minipage}
\begin{minipage}{0.45\linewidth}
\centering
\includegraphics[height=50mm]{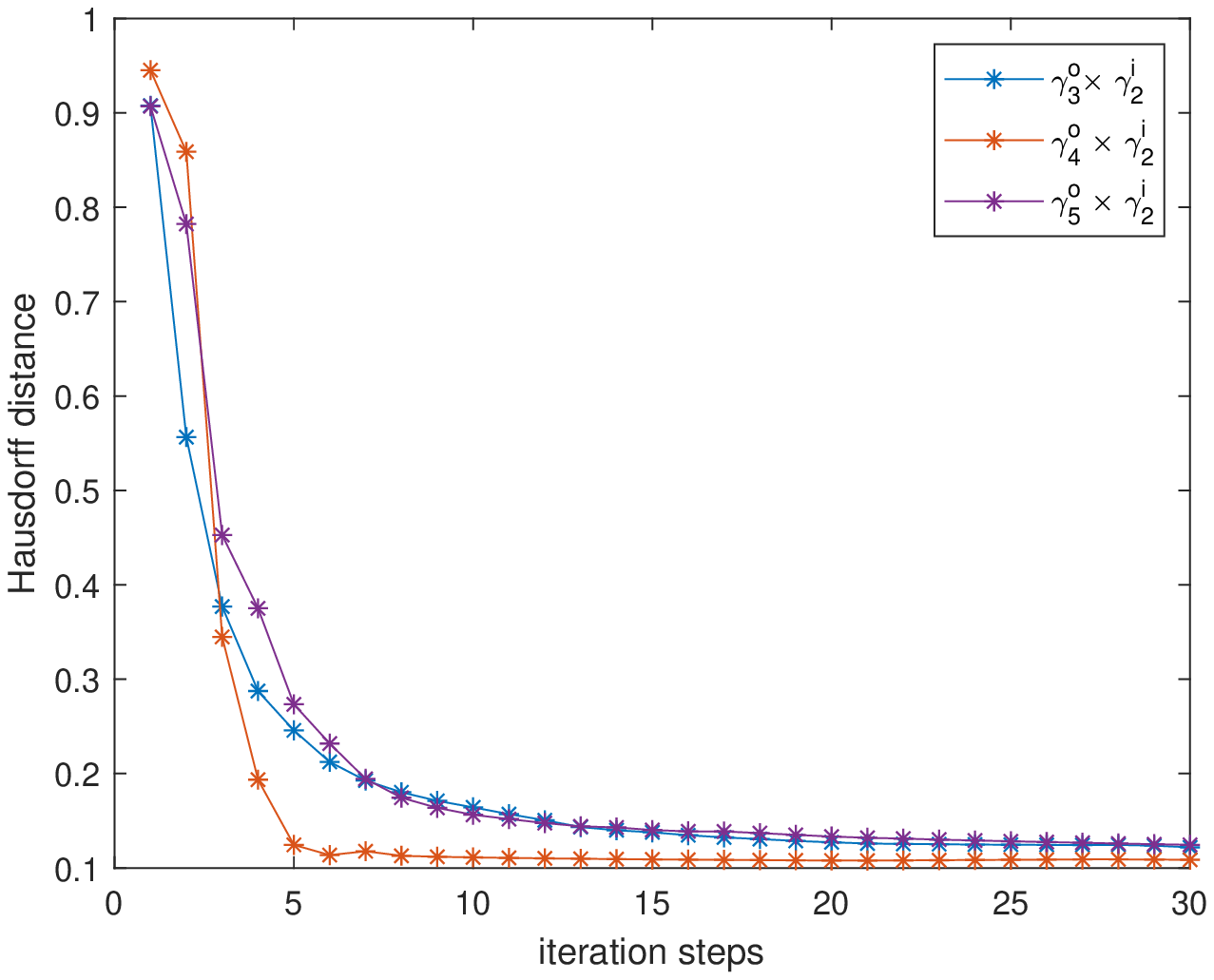}
\end{minipage}
\caption{Hausdorff distance $d_H(\partial \Omega_{\text{exact}},\partial \Omega_{\text{inv}})$ with respect to the iteration steps. Left: $\gamma^i_1$. Right: $\gamma_2^i$.}\label{dH_pear}
\end{figure}

\subsection{Examples for IOSP-P}
Let the measured data be the compressional part of the far-field pattern, i.e.,
$u_p^\infty(\hat{\bm{x}},\bm{d})$, $(\hat{\bm{x}},\bm{d})\in \gamma^o\times\gamma^i$. The obstacle $\Omega$ is a kite with $\partial \Omega$ given by
\begin{equation*}
(\cos \theta+0.65\cos2\theta-0.65, 1.5  \sin \theta)+(-2,3), \;\;\; \theta \in (0, 2\pi].
\end{equation*}

Let $\gamma^i=\gamma_1^i$, i.e.,  one incident direction. The observation apertures are $\gamma^o=\gamma_1^o,\gamma_2^o,\gamma_3^o$.
In Figure~\ref{ESM_kite} (top row), we show the contour plots of the indicator function ${I}_{ESM}(\bm{z})={I_{\bm{z}}}(\bm{z})/{\max_{{\bm{z}}\in T} I_{\bm{z}}(\bm{z})}$,
where the asterisk `*' indicates the reconstructed location by the ESM. The solid curve is the exact boundary.
As expected, when the observation aperture becomes smaller, the result is less satisfactory.
The location reconstructed by the ESM, either inside or outside $\Omega$, is close enough and provides a good initial input for the EnKF.
In Figure \ref{EnKF_kite} (top row), we show the boundary reconstructed by the EnKF in the second step.
The solid line is the exact boundary, the dashed line is the reconstructed boundary, and the asterisk `*' is the refined location generated by the EnKF.
The reconstructions becomes less satisfactory as the observation aperture decreases.
Nonetheless, the reconstruction is very good considering the fact that there is only one incident direction.

Next, we consider the case of $\gamma^i=\gamma_2^i$, $\gamma^o=\gamma_3^o,\gamma_4^o,\gamma_5^o$.
The contour plots of the indicator function ${I}_{ESM}(\bm{z})$ are shown in Figure \ref{ESM_kite} (bottom row).
In Figure \ref{EnKF_kite} (bottom row), we show the boundary reconstructed by the EnKF.
Satisfactory reconstruction can be achieved with quite limited observation data.
In Figure \ref{dH_kite}, the Hausdorff distance $d_H(\partial\Omega_{\text{exact}},\partial\Omega_{\text{inv}})$
between the exact boundary $\partial\Omega_{\text{exact}}$ and the reconstructed boundary $\partial\Omega_{\text{inv}}$ is plotted against to the iteration steps.

The location obtained in the first step using the ESM is critical to the success of the proposed method.
We demonstrate this using a simple example.
Assume that the approximate location of the obstacle is $(z_1,z_2)=(0,0)$. The initial particles are drawn from $\mathcal{N}(z_i, 1)$, $i=1,2$.
In Figure \ref{poorguess_kite}, we display the reconstructions of the boundary and the Hausdorff distance $d_H(\partial\Omega_{\text{exact}},\partial\Omega_{\text{inv}})$
for measured data on $\gamma^o_1\times\gamma^i_1$.
The Hausdorff distance does not become small after reasonable number of iterations and the reconstructed boundary is nowhere close to $\partial \Omega$.

\subsection{Examples for IOSP-S}
We consider the shear part of the far-field pattern,
i.e., $u_s^\infty(\hat{\bm{x}},\bm{d})$, $(\hat{\bm{x}},\bm{d})\in \gamma^o\times\gamma^i$.
The obstacle $\Omega$ is a peanut shape domain with $\partial \Omega$ given by
\begin{equation*}
0.4\sqrt{4\cos^2\theta+\sin^2\theta} (\cos\theta,\sin\theta) +(-2,3), \;\;\; \theta \in (0, 2\pi].
\end{equation*}
In Figure \ref{ESM_peanut}, we plot the contours of the indicator function ${I}_{ESM}(\bm{z})$.
The approximate locations by ESM are marked with asterisks.
In Figure \ref{EnKF_peanut}, we show the reconstructions by the EnKF.
The Hausdorff distance $d_H(\partial\Omega_{\text{exact}},\partial\Omega_{\text{inv}})$ with respect to the number of iterations is shown in Figure \ref{dH_peanut}.

\subsection{Examples for IOSP-F}
Finally, we consider the full far-field pattern, i.e., $\bm{u}^\infty(\hat{\bm{x}},\bm{d})=(u_p^\infty;u_s^\infty)$,
$(\hat{\bm{x}},\bm{d})\in \gamma^o\times\gamma^i$. The obstacle $\Omega$ is a pear shape domain with $\partial \Omega$ given by
\begin{equation*}
\left(\frac{5+\sin3\theta}{6}\cos\theta,\frac{5+\sin3\theta}{6}\sin\theta \right) +(-2,3), \;\;\; \theta \in (0, 2\pi].
\end{equation*}
In Figure \ref{ESM_pear}, we show the contours plots of the indicator function ${I}_{ESM}(\bm{z})$.
The reconstructions by the EnKF are shown in Figure \ref{EnKF_pear}.
In Figure \ref{dH_pear}, we plot the Hausdorff distance $d_H(\partial\Omega_{\text{exact}},\partial\Omega_{\text{inv}})$ with respect to the iteration numbers.

\section{Conclusions}
This paper continues our investigation of the combined deterministic-statistical approach for partial data inverse scattering problems \cite{Li2020, LiEtal2020}.
We propose a two step approach to reconstruct an elastic rigid obstacle with partial data.
In the first step, the approximate location of the unknown obstacle is obtained by the extended sampling method.
In the second step, using the location obtained previously, the ensemble Kalman filter is employed to construct the shape of the obstacle.
Both steps use the same physical model and the same set of measured data.

This approach inherits the merits of the two methods.
Numerical examples show that the proposed method is effective for the inverse elastic scattering problem with partial data.
Demonstrated by the example in Section 5.1, the reconstructed location by the ESM is critical to the success of the ensemble Kalman filter,
which is consistent with the discussions in \cite{Iglesias2013} (Theorem 2.1) and \cite{Iglesias2016} (Proposition 3.1) that the initial ensemble
is a  crucial design parameter. The readers are encouraged to compare the results in this paper with those obtained using the sampling methods
with the same set of measured data \cite{Liu2018, Liu2019}.

\section*{Disclosure statement}

No potential conflict of interest was reported by the author(s).

\bibliographystyle{plain}

%

\end{document}